\documentclass[11pt]{amsart}
\usepackage{amssymb}
\usepackage{verbatim} 
\newtheorem{thm}[equation]{Theorem}
\newtheorem{pro}[equation]{Proposition}

\newtheorem{lem}[equation]{Lemma}

\newtheorem{con}[equation]{Convention}
\newtheorem{DEF}[equation]{Definition}
\newtheorem{rem}[equation]{Remark}

\numberwithin{equation}{section}

\def\a1s{a_1,\cdots, a_s}
\def\a{\alpha}

\def\aa{\mathcal A}

\def\fa{\mathfrak{a}}
\def\fA{\mathfrak{A}}

\def\andd{\quad\hbox{and}\quad}

\def\ad{\hbox{ad}}

\def\fb{\frak{b}}
\def\b{\beta}

\def\bl4{B_{\ell\geq4}}
\def\cc{{\mathcal C}}
\def\fc{\mathfrak{C}}
\def\bbbc{{\mathbb C}}

\def\d{\delta}
\def\D{\Delta}

\def\dd{\mathcal D}
\def\fD{\mathfrak D}

\def\es{e_\sg}

\def\e{\epsilon}
\def\ve{\varepsilon}

\def\bbbf{\mathbb{F}}

\def\gg{{\mathcal G}}

\def\gs{\gg_\sg}

\def\fg{\mathfrak{g}}

\def\hh{{\mathcal H}}

\def\fh{\mathfrak{h}}

\def\fcj{\mathfrak{J}}
\def\fj{\mathfrak{j}}
\def\jj{\mathcal{J}}

\def\kk{\mathcal{K}}

\def\lam{\lambda}

\def\LL{\mathcal{L}}
\def\Lt{\tilde{\mathcal L}}

\def\fm{(\cdot,\cdot)}

\def\bbbn{\mathbb{N}}

\def\bbbr{{\mathbb R}}

\def\1k{\frac{1}{k}}
\def\op{\oplus}
\def\ot{\otimes}

\def\la{\langle}
\def\ra{\rangle}

\def\vp{\varphi}

\def\qed{\hfill$\Box$}
\def\sub{\subseteq}
\def\sg{\sigma}

\def\rcross{R^{\times}}

\def\pf{\noindent{\bf Proof. }}

\def\kk{\mathcal K}

\def\sfs{\mathfrak{s}}

\def\bv{\bold v}
\def\vrs{v_\sg^r}

\def\v{{\mathcal V}}

\def\w{{\mathcal W}}

\def\zn{{\mathbb Z}^{\nu}}
\def\z{{\mathcal Z}}

\def\bbbz{{\mathbb Z}}

\def\1il{1\leq i\leq\ell}

\def\sgn{\hbox{sgn}}

\begin{document}\markboth{}{}
\title{Universal coverings of Lie tori\\ (A finite presentation) }

\author{Saeid Azam, Hiroyuki Yamane, Malihe Yousofzadeh}
\address
{Department of Mathematics\\ University of Isfahan\\Isfahan, Iran,
P.O.Box 81745-163} \email{azam@sci.ui.ac.ir, saeidazam@yahoo.com.}

\address{Department of Pure and Applied Mathematics,
Graduate School of Information Science and Technology,
Osaka University,
Toyonaka, Osaka, 560-0043, Japan}
\email{yamane@ist.osaka-u.ac.jp}

\address
{Department of Mathematics\\ University of Isfahan\\Isfahan, Iran,
P.O.Box 81745-163} \email{ma.yousofzadeh@sci.ui.ac.ir.}
\thanks{This work is partially supported by the Center of Excellence for Mathematics, University of Isfahan.}

\begin{abstract}
Using the well-known recognition and structural theorem(s) for root-graded Lie algebras and their universal coverings, we give a finite presentation for the universal covering algebra of a centerless Lie torus of type $X\not=A,C,BC$. We follow a unified approach for the types under consideration.
\end{abstract}

\maketitle
\markboth{A presentation of Lie tori}{S. Azam, H. Yamane, M. Yousofzadeh}

\setcounter{section}{-1}
\section{ Introduction}
The notion of a Lie torus arises in the study of extended affine
Lie algebras which are natural generalizations of finite
dimensional simple Lie algebras and affine Lie algebras.
Extended affine Lie algebras and Lie tori have been under
intensive investigation in recent years.

It is well understood that the study of an extended affine Lie algebra, in any aspect, somehow relates to the study of
its
core modulo the center, called the centerless core. Here is the
place where the notion of Lie tori arises; the centerless core of
an extended affine Lie algebra is a centerless Lie torus and
conversely any centerless Lie torus is the centerless core of an
extended affine Lie algebra (\cite{Yos}). Lie tori as well as
extended affine Lie algebras are defined axiomatically
(\cite{Yos}, \cite{AABGP}), and their structures are well studied. In particular, it is known that a Lie torus is a {\it root graded Lie algebra}, a notion which turns out to be essential in the structure theory of extended affine Lie algebras and Lie tori.

The structure of root graded Lie algebras of reduced types are
recognized up to central isogeny by \cite{BM} and \cite{BZ}.
Namely a root graded Lie algebra is centrally isogenous to a Lie
algebra having a prescribed structure coming from some known
algebra constructions, including the so called generalized Tits
construction (see Section 1). We refer to these results as
Recognition Theorem(s) (see Theorem \ref{recog}). The Recognition
Theorem together with  several theorems from \cite{ABG},
\cite{BZ}, \cite{AG} and \cite{BGK} enables us to decompose the universal
covering algebra $\mathfrak{A}$ of the centerless core of an
extended affine Lie algebra of type $X\not=A_\ell, C_\ell,
BC_\ell$ as
$$\mathfrak{A}=(\gg\ot \mathcal{A})\op(\v\ot \mathcal{A}^m)\op \dd\eqno{(*)}$$
where $\gg$ is a finite dimensional simple Lie algebra of type
$X$, $\v$ is an irreducible $\gg$-module whose highest weight is
the highest short root of $\gg$, $\mathcal{A}$ is the algebra of
Laurent polynomials in several variables, $\aa^m$ is the direct
sum of $m$ copies of $\aa$ and $\dd$ is a known subalgebra of
$\mathfrak{A}$ related to the inner derivations of the so called
coordinate algebra of $\mathfrak{A}$.

In this work, we use decomposition $(*)$ to give a nullity-free finite presentation for the universal covering  of the centerless core of a Lie torus of type $X$, $X\not=A_\ell,\; C_\ell,\; BC_\ell$. Since the
ingredients appearing in $(*)$, and their  algebra
structures are completely known to us (see Theorems \ref{v1} and \ref{uca}), we are able to select a (finite) set of generators for
$\mathfrak{A}$ and deduce a certain (finite) set of relations among them. This then motivates the generators and defining relations of our presented Lie algebra.

There had been several attempts in presenting a Lie torus (up to center) by a (finite) set of generators and relations. Historically, we may name the works of \cite{Ka} and \cite{MRY} for toroidal Lie algebras; \cite{SY} and
\cite{Yam} for elliptic Lie algebras (2-extended affine Lie algebras); and \cite{You} for Lie tori of type $B_\ell$, $(\ell\geq 3)$. It is worth mentioning that the universal covering of the centerless core of a finite dimensional simple Lie algebra or an affine Kac--Moody Lie algebra is just itself or its derived algebra, respectively. Therefore,
for  finite or affine case the well-known Serre's type presentation is in fact given for the universal covering of the centerless core. This might justify why we are considering the universal covering instead of a Lie torus itself.

The paper is arranged as follows. In Section 1, we record  some necessary backgrounds and results needed for describing the universal covering of the centerless core of a Lie torus in the form $(*)$. In brief, we establish some minor results regarding finite root systems and irreducible modules of  finite dimensional simple Lie algebras (Lemmas \ref{rev7}-\ref{rev1}). The generalized Tits construction and basic properties of root graded Lie algebras are reviewed. The Recognition Theorem(s) from \cite{BM} and \cite{BZ} regarding the structure of root graded Lie algebras are recalled, and some results from \cite{BGK}, \cite{AG} and \cite{ABG} regarding the universal covering algebras of root graded Lie algebras and Lie tori are restated. These all together, enable us to have decomposition $(*)$. The section is concluded with obtaining from $(*)$ a finite set of generators and a certain set of relations for $\mathfrak{A}$.

In Section 2, we introduce the generators and defining relations of our presented Lie algebra $\LL$, associated to a Lie torus of type $X\not=A_\ell, C_\ell, BC_\ell$, and  deduce ceratin immediate consequences from them which will be essential for the proof of our main theorem (Theorem \ref{main}).
For the convenience of the reader, we have shifted the proofs of some of these results containing complicated technicalities to Section 4. Section 3 is devoted to the proof of our Main Theorem (Theorem \ref{main}); the universal covering algebra of the centerless Lie torus under consideration is centrally isogenous to its corresponding presented Lie algebra $\LL$. The proof is proceeded in steps. In step 1, we prove that $\LL$ is root graded and we obtain some information regarding its center. In the remaining steps we complete the proof for simply laced types, types $F_4$ and $G_2$, respectively. For type $B_\ell$ ($\ell\geq 3$), we refer the reader to \cite{You}, though the proof for this type can also be deduced from our setting. This work continues the line of research appeared in \cite{You} for type $B_\ell$ $(\ell\geq 3)$ where here we have enlarged our setting in order to have a unified approach for the types under consideration. Needless to say, in both works the defining generators and relations of the considered presented Lie algebras are inspired from the
algebra structures among the ingredients involved in decomposition $(*)$.

For non-simply laced types, we have used certain set of data (see
(\ref{str})) which are obtained using several Mathematica
programming.

\section{Terminology and
prerequisites}\label{Terminology} Throughout this paper we suppose
that  $R$ is  an irreducible reduced finite root system  of rank
greater than 3 which is not of type $A,C$. By $\rcross,$ we mean
$R \setminus\{0\}.$ All vector spaces and tensor products are
taken over the field of  complex numbers $\bbbc.$ For an algebra
$A$ which is not a Lie algebra, by $[a,b],$ $a,b\in A,$ we mean
$ab-ba.$  If $\a,\b$ are two elements of a unital   algebra $\aa$
which is  commutative associative, alternative or Jordan, we
denote by $D_{\a,\b},$ the inner derivation of $\aa$ based on
$\a,\b$, we also use $D_{\aa,\aa}$ to denote the set of all inner
derivations of $\aa$ \cite{Sc}.
 For a positive integer $\nu,$ we denote by $A_{[\nu]},$ the
commutative associative algebra of Laurent polynomials in $\nu$
variables $t_1,\ldots,t_{\nu}$ and we use $A_{[\nu]}^m,$ for a
positive integer $m,$ to denote the direct sum of $m$ copies
of $A_{[\nu]}.$ Also by $A_{[\nu]}^0,$ we mean the trivial  vector
space. For $\sg=(n_1,\ldots,n_\nu)\in \bbbz^\nu$ we denote
$t_1^{n_1}\ldots t_\nu^{n_{\nu}}$ by $t^\sg$. We also make a
convention that for elements  $x,x_1,\ldots,x_m,y$ of a Lie
algebra, by an expression of the form $[x_1,\ldots,x_m,y]$ we
always mean $[x_1,\cdots[x_{m-1}, [x_m,y]]\cdots].$   In this case
if $m=0$ we interpret $[x_1,\ldots,x_m,y]$ as $y,$ also by
$[x^m,y],$ we mean $[x,\ldots , x,y]$ in which $x$ appears
$m$-times. In a vector space, we define  ${\sum_{i=1}^m\cdots}$ to
be zero if $m=0.$ For the sake of simplicity of notation we denote
the set $\{1,\ldots,k\}$ be $J_k$, $k$ a positive integer. For a
Lie algebra $\LL$,  by $Z(\LL),$ we mean the center of $\LL$ and
by a  centerless Lie algebra, we mean a Lie algebra with trivial
center. We start the paper by recalling some (in general)
non-associative algebra constructions and some related theorems
which are of our use throughout this work.
\subsection{General facts}

Suppose that   $\gg$ is  a finite dimensional simple Lie algebra
over $\bbbc$ of type $X\neq A,C$ and of  rank $\ell, $  also
suppose that $R$ is the corresponding irreducible finite root
system with a root base $\D$.  Let $R^+$, $R^+_{sh}(R_{sh})$ and
$R^+_{lg}(R_{lg})$ be the set of positive roots, the set of
positive short roots (the set of short roots), and the set of
positive long roots (the set of long roots), respectively. If $X$
is simply laced we assume by convention that $R^+_{lg}=\emptyset$.
Also assume that $n$ is the number of positive roots. We arrange
the set of positive roots as follows:
\begin{equation}\label{fin2}
\stackrel{R^+_{lg}\cap\D}{\overbrace{\a_1,\ldots,\a_{n_\ell}}},\stackrel{R^+_{sh}}
{\overbrace{\a_{n_\ell+1},\ldots,\a_\ell,\ldots,
\a_{n_s}}},\stackrel{R^+_{lg}\setminus\D}{\overbrace{\a_{{n_s}+1},\ldots,\a_n}}.
\end{equation}
Clearly if there is no long root ($n_\ell=0$), then $n=n_s.$ We
may assume that $\a_{n_s}$ is the highest short root and $\a_n$ is
the highest root. For $t\in J_n,$ we fix  $e_t\in\gg_{\a_t}$ and
$f_t\in\gg_{-\a_t}$ such that $(e_t,h_t:=[e_t,f_t],f_t)$ is an
$\frak{sl_2}-$triple.


\smallskip

\smallskip
Using  the finite dimensional theory, one can prove  the following
lemma:
\begin{lem}\label{rev7}
 (i) If $R_{lg}\not=\emptyset$, $R_{sh}\sub R_{lg}+ R_{sh}.$

\smallskip

(ii) For $\a\in R_{sh}^+,$ there is $\b\in R_{sh}^+$ such that
$\a+\b\in R.$
\end{lem}

\begin{lem}
\label{rev1} Let $\w$ be an irreducible finite dimensional $\gg$-module.

(i) Suppose that highest weight $\w$ is $\lam\in R$ and take $\Pi$ to be the
set of weights of $\w.$ Then for $\mu\in \Pi\setminus\{0\},$ there
are $\a_{i_1},\ldots,\a_{i_p}\in R^+$ such that
$\w_\mu=\gg_{-\a_{i_p}}\cdot \cdots\cdot
\gg_{-\a_{i_1}}\cdot\w_\lam$ with $\displaystyle{\a_{i_k}\neq
\pm(\lam-\sum_{r=1}^{k-1}\a_{i_r})}.$

(ii) Suppose that  the
highest weight of $\w$ is a short root  and consider  the weight space decomposition
$\w=\w_0\op\sum_{i=n_\ell+1}^{n_s}\bbbc w_{\pm i}.$
for $\w$ in which
for $n_\ell+1\leq i\leq n_s,$ $w_{\pm i}$ is a weight vector
of weight $\pm\a_i.$ Then for $n_\ell+1\leq i\leq n_s,$
$$e_i\cdot f_i\cdot w_i=2w_i,\;\; e_i\cdot f_i\cdot f_i\cdot
w_i=2f_i\cdot w_i,\;\; \bbbc e_i\cdot w_{-i}=\bbbc f_i\cdot w_i$$
and
 $\w_0$ is spanned by $f_i\cdot w_i,$ $n_\ell+1\leq i\leq \ell.$

\end{lem}
\pf $(i)$ Using a simple   argument from  the finite dimensional  theory,  we are done.

$(ii)$ We first note that if $\b\in\Pi\setminus\{0\},$ then
$\dim(\w_\b)=1,$ also if $\b\in \Pi$ and $\a\in R$ are  such that
$\a+\b\in \Pi,$ then $\gg_\a\cdot\w_\b=\w_{\a+\b}.$ Now let
${\mathcal U}(N^-)$ be the universal enveloping algebra of
$N^-:=\displaystyle{\sum_{\a<0}\gg_{\a}}.$ Then if $w$ is a
highest vector, we have ${\mathcal U}(N^-)\cdot w=\w$. Also for
each $\b\in R^+$ we have
\begin{equation}
\gg_{-\b}=[\gg_{-\b_s},\ldots,\gg_{-\b_1}] \hbox{ for some
$\b_i\in\D$}. \label{rev2}
\end{equation}
 Now combining these facts and using
module action, the required expression for $\w_{\mu}$ with $\mu\in
R^+$ follows. Note that the claim concerning $\a_{i_k}$'s is
automatically satisfied.

Next let $\mu\in -R^+,$ then there is $\lam\in R^+$ such that
either $\mu-\lam\in R^+$ or $\mu+\lam\in R^+.$ Moreover if $\mu\in
R^+_{sh},$  $\lam$  can be chosen as an element of $R^+_{sh}$ such
that $\mu+\lam\in R.$
\begin{eqnarray*}
\w_{\mu}=\left\{\begin{array}{ll} \gg_{-(\lam+\mu)}\cdot
\w_{-\mu}&\hbox{if } \lam+\mu\in
R^+, \\
\gg_{\lam-\mu}\cdot\gg_{\mu}\cdot\w_{\mu-\lam}&\hbox{if }
\lam-\mu\in R^+.
\end{array}\right.
\end{eqnarray*}
Now we may replace $\w_{\mu}$ or $\w_{\mu-\lam}$ by an expression
obtained from the previous step to get the required expression for
$\w_{\mu}$. Note that from the way we found this expression for
$\w_{\mu}$, it becomes clear that the claim concerning
$\a_{j_k}$'s holds. \qed

\subsubsection{Generalized Tits construction} Let $A$ be a unital commutative associative
algebra  and assume that $X$ is a unital algebra over $A$. A {\it
normalized trace} on $X$ is an $A$-linear map $T:X\longrightarrow
A$ satisfying, for $x,x',x''\in X$,
$$T(1)=1,\quad T(xx')=T(x'x),\quad T\big((xx')x''\big)=T\big(x(x'x'')\big).$$
If $T$ is a normalized trace, the maps $t$ and $*$ defined by
\begin{equation}\label{v5}
\begin{array}{c}
t:X\times X\longrightarrow A;\;(x,y)\mapsto T(xy),\vspace{2mm}\\
*:X\times X\longrightarrow X;\; (x,y)\mapsto xy-t(x,y)1,
\end{array}
\end{equation}
are called, respectively, the {\it trace form} and the {\it
$*$-operator} with respect to $T.$ We use the same symbols $T$,
$t$ and $*$ to denote the normalized trace, the trace form and the
corresponding $*$-operator for different algebras.
 We also  have
$X=A1\op X_0,$ where $X_0:=\{x\in X\mid T(x)=0\}$. Let
$\hbox{Der}^0_A(X)$ be the Lie subalgebra of the $A$-derivations
of $X$ which send $X_0$ to $X_0.$ Let $D$ be a Lie subalgebra of
$\hbox{Der}^0_A(X)$ and assume there is an $A$-bilinear
transformation $\a:X_0\times X_0\longrightarrow D$ which is
skew-symmetric.

Suppose now that $\mathfrak{A}$ is another unital commutative
associative algebra over $\bbbc$ and $Y$, $Y_0$, $D'$ are
similarly chosen over $\mathfrak{A}.$ Assume $\b:Y_0\times
Y_0\longrightarrow D'$ is a $\mathfrak{A}$-bilinear transformation
which is skew symmetric. Let the following relations
hold:$$[d,\a(x,x')]=\a(dx,x')+\a(x,dx')\andd
[d',\b(y,y')]=\b(d'y,y')+\b(y,d'y')$$ for $x,x'\in X_0$, $y,y'\in
Y_0,$ $d\in D$ and $d'\in D'$. Then the vector space
$${\mathcal T}(X/A,Y/\mathfrak{A}):=(D\ot\mathfrak{A})\op(X_0\ot Y_0)\op(A\ot D'),$$
 provides an algebra over $\bbbc$ with the anticommutative multiplication given by
\begin{equation}
\begin{array}{l}
\ [d_1\ot b,a\ot d'_1]=0,\;\; [d_1\ot b,d_2\ot b']=[d_1,d_2]\ot bb',\\
\ [a\ot d'_1,a'\ot d_2']=aa'\ot [d'_1,d'_2],\\
\ [d_1\ot b,x\ot y]=d_1x\ot by=-[x\ot y,d_1\ot b],\\
\ [a\ot d'_1,x\ot y]=ax\ot d'_1y=-[x\ot y,a\ot d'_1],\\
\ [x\ot y,x'\ot y']=\a(x,x')\ot t(y,y')+(x*x')\ot(y* y')+t(x,x')\ot \b(y,y')\\
\end{array}
\label{m9}
\end{equation}
for $d_1,d_2\in D$, $b,b'\in\mathfrak{A},$ $a,a'\in A$,
$d'_1,d'_2\in D'$, $x,x'\in X_0$ and $y,y'\in Y_0$. This is called
a {\it  generalized Tits construction}. If $X,Y$ are suitably
chosen, then ${\mathcal T}(X/A,Y/\mathfrak{A})$ will be a Lie
algebra \cite[ Proposition 3.9]{BZ}.

\subsubsection{Root-graded Lie algebras}
The main purpose  of this work is to give a finite presentation
for universal covering algebra of a Lie torus. Lie tori are
centerless cores of extended affine Lie algebras. In \cite{AG} and
\cite{BGK}, the authors classify the Lie tori under consideration,
using the so called recognition theorems for root graded Lie
algebras. We recall here these theorems and some related topics
which will be of use in the the proof of our main theorem. In what
follows we denote the $8$-dimensional Octonion (Cayley) algebra by
${\mathfrak C}.$ We consider the  usual normalized trace on
${\mathfrak C}$ and denote its subspace of trace zero elements
with ${\mathfrak C}_0$. One knows that the algebra of inner
derivations of ${\mathfrak C}$ is a finite dimensional simple Lie
algebra of type $G_2$. Also we denote by ${\mathfrak J}$ the
exceptional simple Jordan algebra whose inner derivations form a
finite dimensional simple algebra of type $F_4$. In fact if
${\mathfrak C}_{3\times 3}$ is the algebra of $3\times 3$ matrices
with entries from the Octonion algebra ${\mathfrak C}$, then
${\mathfrak J}$ is its subspace of self-adjoint elements, under
transpose-conjugate involution $x\mapsto \bar{x}^t$, with the
product $x\cdot y:=(xy+yx)/2$. We also consider the usual
normalized trace on $\jj$ and denote the subspace of trace zero
elements of ${\mathfrak J}$ by ${\mathfrak J}_0$ (see Section \ref{app} for details).

\begin{DEF} {\rm Suppose that  $\gg$ is  a finite dimensional  simple Lie
algebra over $\bbbc$ with a Cartan subalgebra $\hh$ and root
system $R$ so that $\gg$ has a root space decomposition
$\gg=\op_{\mu\in R}\gg_{\mu}$ with $\hh=\gg_0.$
 An
{\it $R$-graded Lie algebra} $\LL$ over $\bbbc$ with {\it grading
pair} $(\gg,\hh)$ is  a Lie algebra satisfying the following
conditions:

(i) $\LL$ contains $\gg$ as a subalgebra,

(ii) $\LL=\op_{\mu\in R}\LL_{\mu}, \hbox{\;where\;}
\LL_{\mu}:=\{x\in\LL\mid [h,x]=\mu(h)x\hbox{\;for
all\;}h\in\hh\},$

(iii) $\LL_0=\sum_{\mu\in R^\times}[\LL_{\mu},\LL_{-\mu}]. $

For a positive integer $\nu$, an $R$-graded Lie algebra  $\LL$
with grading pair $(\gg,\hh)$ is called {\it $(R,\zn)$-graded }if
$\LL=\op_{\sg\in \zn}\LL^\sg$ is a $\zn$-graded Lie algebra such
that $\gg\sub \LL^0$ and supp$(\LL):=\{\sg\in\zn\mid
\LL^\sg\neq\{0\}\}$ generates $\zn$. Since $\gg\subseteq\LL^0,$
$\LL^\sg$ is an $\hh$-module for $\sg\in \zn$ and so we have
$\LL=\op_{\mu\in R}\op_{\sg\in \zn}\LL^\sg_\mu$ where
$\LL_\mu^\sg:=\LL^\sg\cap\LL_\mu$ for $\sg\in \zn$ and $\mu\in R.$
An $(R,\zn)$-graded Lie algebra  $\LL$   is called {\it division
$(R,\zn)$-graded} if for each $\mu\in R^\times$, $\sg\in \zn$ and
$0\neq x\in\LL^\sg_\mu$, there exists $y\in\LL^{-\sg}_{-\mu}$ such
that modulo $Z(\LL),$  $[x,y]$ equals to the unique element of
$\hh$ representing $\mu$ through the induced form on the dual of
$\hh.$ A division $(R,\bbbz^\nu)$-graded Lie algebra $\LL$ with
$\dim_\bbbc(\LL_\mu^\sg)\leq1$ for all $\sg\in\bbbz^\nu$ and
$\mu\in R^\times$ is called a {\it Lie $\nu$-torus}  or  simply a
{\it Lie torus.} In this case the set $\{\a+\sg\mid \a\in
R,\sg\in\zn, \LL_\a^\sg\neq\{0\}\}\sub \hbox{span}_\bbbr R\op\zn$
is called the {\it root system } of $\LL.$ For $\mu\in R^\times,$
define $S_\mu:=\{\sg\in\zn\mid \LL_{\mu}^\sg\neq\{0\}\},$ by
\cite[Theorem 1.5]{Yos}, $S_\mu=S_\nu$ if $\mu$ and
$\nu$ have the same length. If two root lengths occur, we set $S:=
S_\mu$ for any choice of a short root $\mu$ and $L:=S_\nu$ for any
choice of a long root $\nu$ and  call $(S,L)$ the {\it
corresponding pair} of $\LL.$}
\end{DEF}

\begin{DEF} {\rm Let $B$ be a unital commutative associative algebra,
$\w$ be a $B$-module and $g:\w\times\w\longrightarrow B$ be a
symmetric $B$-bilinear form on $\w$. Then $J(\w):=B1\op\w$
 with the  multiplication, for $w,w'\in \w,\;b,b'\in B$,
$$(b1+w)\cdot(b'1+w')=bb'1+g(w,w')1+bw'+b'w$$ is a Jordan algebra called the
{\it Clifford Jordan algebra of $g$}.}
\end{DEF}

\begin{DEF}
Two perfect Lie algebras are said to be {\it centrally isogenous}
if they have the same universal covering algebra, up to
isomorphism.
\end{DEF}
\begin{thm}[\bf Recognition Theorem]\label{recog}
Let $\LL$ be an $R$-graded Lie algebra with grading pair
($\gg,\hh)$.

\medskip

(i) \cite{BM} If $R$ is simply laced, there is a commutative
associative unital  algebra $A$ such that $\LL$ is centrally
isogenous  with $\gg\ot A.$

\medskip

(ii) \cite{BZ} If  $R$ is of type $B_\ell,$ $\ell\geq3,$ there
exists a unital commutative associative algebra $A$ and a unital
$A$-module $B$ with  a symmetric $A$-bilinear  form $\fm: B\times
B\longrightarrow  A$ such that $\LL$ is centrally isogenous  with
$${\mathcal T}\big(J(\v)/\bbbc,J(B)/A\big)=(\gg\ot A)\op(\v\ot B)\op D_{J(B),J(B)}$$
where $\v$ is the  $(2\ell+1)$-dimensional vector space equipped
with a non-degenerate  symmetric bilinear form with respect to
which skew symmetric endomorphisms of $\v$ is isomorphic to $\gg.$

\medskip

(iii) \cite{BZ} If $R$ is of type $G_2,$ then there is a unital
commutative associative algebra $A$ and a unital Jordan algebra
$J$ over $A$ having a normalized trace $T$ satisfying the identity
{\small $$ch_3(y):=y^3 -3T(y)y^2 +(\frac{9}{2}T(y)^2-\frac{3}{2}
T(y^2))y+(T(y^3)-\frac{9}{2}T(y^2)T(y) + \frac{9}{2}T(y)^3)1=0$$}
such that $\LL$ is centrally isogenous with
$${\mathcal T}(\fc/\bbbc,J/A)=(\gg\ot A)\op(\fc_0\ot J_0)\op D_{J,J}$$
where $\fc_0$ and $J_0$ are the trace zero elements of the
Octonion algebra $\fc$ and $J$ respectively.

\medskip
(iv) \cite{BZ} If $R$ is of type $F_4,$ there exists a unital
commutative associative algebra $A$ and a unital  alternative
algebra $\aa$ over $A$ having  a normalized trace $T$ satisfying
the identity
$$ch_2(y):=y^2 -2T(y)y + (2T(y)^2-T(y^2))1 = 0$$ such that $\LL$ is
centrally isogenous with
$${\mathcal T}(\fcj/\bbbc,\aa/A)=(\gg\ot A)\op(\fcj_0\ot \aa_0)\op D_{\aa,\aa},$$
where $\aa_0$ is the subspace of trace zero elements of $\aa$.
\end{thm}

One knows  that any $R$-graded Lie algebra $\LL$ has a
decomposition as
$$(\gg\ot A)\op (\v\ot B)\op D$$ where $\gg$ is the finite dimensional  simple Lie
algebra of type $R,$ $\v$ is the irreducible finite dimensional
$\gg$-module whose highest weight is the highest short root, equipped with  a normalized trace  whose highest
weight is a short root and $\dd$ is a subalgebra of $\LL.$ Also
$A,B$ are two vector spaces with, $B={0}$ in simply laced cases,
such that $\fa:=A\op B$ is an algebra which we refer to as the
coordinate algebra of $\LL$ (see \cite{BGKN},\cite{BM}).

\begin{thm}[\hbox{\cite[Theorem 4.13]{ABG}}]  Let $\LL=(\gg\ot
A)\op(\v\ot B)\op D_{\fa,\fa}$ be a centerless $R$-graded Lie
algebra where $\fa=A\op B$ is the coordinate algebra of $\LL.$
 Take $\sfs$ to be the subspace of $\fa\ot \fa $ spanned by the elements of the form
$$(\a\ot \b)+(\b\ot \a),\; (\a\b\ot \gamma)+(\b\gamma\ot \a)+(\gamma\a\ot\b),\;a\ot b$$ for $a\in A,b\in B,\a,\b,\gamma\in\fa.$
Consider the factor space $$\{\fa,\fa\}:=(\fa\ot\fa)/\sfs$$ and
for $\a,\b\in \fa,$ let $\{\a,\b\}$ denote $(\a\ot \b)+\sfs$ in
$\{\fa,\fa\}.$ Set $\Hat\LL:=(\gg\ot A)\op(\v\ot
B)\op\{\fa,\fa\}.$ Define a multiplication on $\Hat\LL$ by
\begin{equation}\label{ubracket}
\begin{array}{l}
\ [x\ot a, x'\ot a']=[x, x']\ot aa'+\kappa(x,x')\{a,a'\},\\
\ [x\ot a, v\ot b]=xv\ot ab=-[ v\ot b,x\ot a],\\
\ [x\ot a, \{\a,\a'\}]=0=-[\{\a,\a'\},x\ot a],\\
\ [v\ot b, v'\ot b']=D_{v, v'}\ot t(b, b')+(v*v')\ot (b*b')+t(v, v')\{b, b'\},\\
\  [\{\a,\a'\},  v\ot b]=v\ot D_{\a,\a'}b=-[v\ot b,\{\a,\a'\}],\\
\
[\{\a,\a'\},\{\b,\b'\}]=\{D_{\a,\a'}\b,\b'\}+\{\b,D_{\a,\a'}\b'\},
\end{array}
\end{equation}
for $x,x'\in\gg,\; a,a'\in A, \;v,v'\in \v,\; b,b'\in
B\hbox{\;and\;}\a,\a',\b,\b'\in\fa$ where $\kappa$ denotes the
Killing form of $\gg.$ Also consider the map
$\Hat\pi:\Hat\LL\longrightarrow\LL$ given by $ x\ot a\mapsto x\ot
a;\;u\ot b\mapsto u\ot b;\;\{\a,\a'\}\mapsto D_{\a,\a'}$. Then
$(\Hat\LL,\Hat\pi)$ is the universal covering algebra of
$\LL.$\label{uca}
\end{thm}

\begin{thm}
\label{v1}  Let $\LL$ be a Lie $\nu$-torus of type $R$ with
grading pair ($\gg,\hh)$. Then $\LL$ is centrally isogenous with
one of the followings

\smallskip

(i) {\rm \hbox{\cite{BGK}}} $\gg\ot A_{[\nu]}$,  if $R$ is simply
laced.

\medskip

(ii) {\rm \hbox{\cite{AG}}} $ {\mathcal
T}(J(\v)/\bbbc,J(A_{[\nu]}^{m})/A_{[\nu]})$, if $R$ is of type
$B_\ell$ ($\ell\geq 3$), where $\v$ is the $2\ell+1$-dimensional
vector space having a symmetric non-degenerate bilinear form with
respect to which the set of  skew symmetric endomorphisms of $\v$
is isomorphic to $\gg$, $m\geq 2$, and $J(A_{[\nu]}^{m})$ is the
Clifford Jordan algebra with respect to the symmetric
$A_{[\nu]}$-bilinear form on $A_{[\nu]}^{m}$ given by
\begin{equation}
\begin{array}{c}
g:A_{[\nu]}^{m}\times A_{[\nu]}^{m}\longrightarrow A_{[\nu]}\\
g\big(\sum_{r=1}^{m}a_rw_r,\sum_{r=1}^{m}b_rw_r)=\sum_{r=1}^{m}a_rb_rt^{\tau_r}
\end{array}
\label{a15}
\end{equation}
where $\tau_0,\ldots,\tau_{m}\in\zn$ satisfy $\tau_0=0$ and
$\tau_r{\not\equiv}\tau_s$ (mod $2\zn$) for $0\leq s\neq r\leq m.$
Moreover the root system of $\LL$ is of the form
$(S+S)\cup(R_{sh}+S)\cup(R_{lg}+2\zn)$ where
$S=\cup_{j=0}^{m}(2\zn+\tau_j)$. (Here $\{w_1,\ldots,w_{m}\}$ is
the standard basis of $A_{[\nu]}^{m}$ as a free
$A_{[\nu]}$-module.)

\medskip

(iii) {\rm \hbox{\cite{AG}}} ${\mathcal
T}(\fcj/\bbbc,\cc/A_{[\nu]})$ if $R$ is of type $F_4,$ where
$\cc=\aa_p$ (see Subsection \ref{app3}) for some $0\leq p\leq 3.$ Moreover, if ${\mathcal
C}=\aa_p$, $0\leq p\leq 3$, then the root system of $\LL$ is
$\bbbz^\nu\cup(R_{sh}+\bbbz^\nu)\cup(R_{lg}+(2\bbbz^p\oplus\bbbz^{\nu-p}))$.

\medskip

(iv) {\rm \hbox{\cite{AG}}} ${\mathcal
T}(\fc/\bbbc,\jj/A_{[\nu]})$ if $R$ is of type $G_2,$ where
$ \jj=\jj_p$  (see Subsection \ref{app3}) for some $0\leq p\leq 3.$ Moreover, if ${\mathcal
J}=\jj_p,$ $0\leq p\leq 3$, then the root system of $\LL$ is
$\bbbz^\nu\cup(R_{sh}+\bbbz^\nu)\cup(R_{lg}+(3\bbbz^p\oplus\bbbz^{\nu-p}))$.
\end{thm}

\begin{rem}
\label{jadid} {\rm Using the same notation as  before, if
$\fa=A\op B$ is the coordinate algebra of a Lie torus $\LL$ with
the universal covering algebra $\fA,$ one gets from \cite{BZ} that
$A$ is a subset of the {\it associative center} of $\fa$ and so
(\ref{ubracket}) implies that $\{A,A\}=\{\{a,a'\}\mid a,a'\in
A\}\sub Z(\fA).$}
\end{rem}
\subsection{Induced relations} As before we assume that $R$ is a finite irreducible
reduced root system of type $X\not= A, C$. Let $\mathfrak{A}$ be
the universal covering algebra of a centerless Lie torus $\LL$ of type
$X$. Using Theorems \ref{v1} and \ref{uca}, we may write
\begin{equation}\label{structure}
\mathfrak{A}=(\gg\ot A_{[\nu]})\op(\v\ot A_{[\nu]}^m)\op \dd
 \end{equation}
for some nonnegative integer $m,$ where $\gg$ is a finite
dimensional simple Lie algebra of type $X$, $\v$ is an irreducible
$\gg$-module whose highest weight is the highest short root of
$\gg$, and $\dd$ is a known subalgebra of $\mathfrak{A}.$ We
remark that if $X$ is simply laced, we have by convention $m=0$
and so the middle part in (\ref{structure}) vanishes. Since the
ingredients appearing in (\ref{structure}), and their  algebra
structures   are completely known to us (see \cite{BGK},
\cite{AG} and \cite{ABG}), we are able to select a set of generators for
$\mathfrak{A}$ and deduce certain relations among them which in
turn motivate the generators and defining relations of a
presented Lie algebra (isomorphic to $\mathfrak{A}$) which we
define in the next section.

 If $\LL=\sum_{\a\in R}\LL_\a$, $\gg=\sum_{\a\in R}\gg_\a$ and $\v=\sum_{\a\in R_{sh}}\v_\a\oplus\v_0$, are the
 corresponding weight space decompositions of $\LL$, $\gg$ and $\v$ respectively, then
\begin{equation}\label{decom1}
\LL_\a=\left\{\begin{array}{ll} \gg_\a\ot A_{[\nu]}&\hbox{if
$\a\in R_{lg}$}\\
(\gg_\a\ot A_{[\nu]})\op(\v_\a\ot A_{[\nu]}^m)&\hbox{if
$\a\in R_{sh}$}\\
\displaystyle{\sum_{\b\in R^\times}[\LL_\b},\LL_{-\b}]&\hbox{if
$\a=0.$}
\end{array}
\right.
\end{equation}

One knows that if $X$ is not simply laced, $\gg$ is the set of
inner derivations of an algebra equipped with a normalized trace
$T$, and that $\v$ is the set of zeros of this trace. We fix a
highest weight vector $\bv:=v_{n_s}$ of $\v.$ Next using  Proposition
\ref{rev1}($i$),   we  fix
$j_1^t,\ldots,j_{n_t}^t,k_1^t,\ldots,k_{n'_t}^t\in J_\ell,$
$n_\ell+1\leq t\leq n_s,$ such that
\begin{equation}\label{my1}
v_t:=[f_{j^t_1},\ldots,f_{j^t_{n_t}},\bv]\andd
v_{-t}=[f_{k^t_1},\ldots,f_{k^t_{n'_t}},\bv]
\end{equation}
are nonzero elements of $\v_{\a_t}$ and $\v_{-\a_t}$ respectively.
Then
$$\v=(\bigoplus_{t=n_\ell+1}^{n_s}\bbbc v_{\pm t})\oplus\v_0,$$
where $\v_0$ is the corresponding zero weight space. We recall the
trace form $t$ as in (\ref{v5}) and ask the reader to check that
\begin{equation}\label{fin}
t(v_{n_s},v_{\pm i})=0;\;\;\;\;(n_\ell+1\leq i\leq n_s-1).
\end{equation}

Now consider the coordinate algebra  $\fa=A_{[\nu]}\op
A_{[\nu]}^m$ of $\LL.$ Using Theorem  \ref{recog}, we get that
$A_{[\nu]}^m$ is the  kernel of a normalized trace $T$ of $\fa.$
 We recall the trace form $t\fm:\fa\times\fa\longrightarrow
A_{[\nu]}$ on $\fa;$ $t(x,y):=T(xy)$ and the operator
$*:\fa\times\fa\longrightarrow A_{[\nu]}^m $ defined by
$x*y:=xy-t(x,y)1.$ Let $\{w_r\mid 1\leq r\leq m\}$ be the standard
basis for $A_{[\nu]}$-module $A_{[\nu]}^m.$ Then one observes that
for $1\leq r,s\leq m,$ there are $a_{r,s},a'_{r,s}\in\bbbc,$
$\sg_{r,s},\sg'_{r,s}\in\zn$ and $t_{r,s}\in J_m$ with
\begin{equation}\label{fin8}
\parbox{4.3in}{$\sg_{r,s}=\sg_{s,r},$  $t_{r,s}=t_{s,r},$ $a_{r,s}=-a_{s,r}$  for type $F_4,$ $a_{r,s}=a_{s,r}$  for type
$B,G_2,$   $(1-\d_{r,s})a'_{r,s}+\d_{r,s}a_{r,s}=0$ and
$(1-\d_{r,s})a_{r,s}+\d_{r,s}a'_{r,s}\neq0.$ } \end{equation} such
that
\begin{equation} t(w_r,w_s)=a'_{r,s}t^{\sg'_{r,s}}\andd
w_r*w_s=a_{r,s}t^{\sg_{r,s}}w_{t_{r,s}}. \label{newm1}
\end{equation}

We set  $$\fD_c:=\{a_{r,s},a'_{r,s},\sg_{r,s},\sg'_{r,s},t_{r,s}\mid r,s\in  J_m\}.$$ We consider $\fD_c$ as a data-set
for the coordinate algebra $\fa $ of $\LL.$ In fact, as  $\{w_i\mid i\in J_m\}$ is an $A_{[\nu]}$-basis for $A_{[\nu]}^m$,
also $*$-operator and $t\fm$ are $A_{[\nu]}$-bilinear, the data in $\fD_c$ can  completely  describe  the structure of  $\fa=A_{[\nu]}\op
A_{[\nu]}^m.$

Next for $n_\ell+1\leq i\leq n_s-1,$ define $t'_{\pm i}\in J_n,$
$n_\ell+1\leq t_{ i}\leq n_s$ to be as follow. If
$\a_{n_s}\pm\a_i$ is a root, take $t'_{\pm i}$ to be such that
$\a_{n_s}\pm\a_i=\a_{t'_{\pm i}}$ and $n,$ otherwise. Also  if
$\a_{n_s}-\a_i$ is a short root, take $t_{ i}$ to be such that
$\a_{n_s}-\a_i=\a_{t_{ i}}$ and $n_s,$ otherwise. One observes
that there are  $m_{\pm i},m'_{\pm i}\in\bbbc$ with $m'_{\pm
i'}=0$ and $m_{\pm i}=0$ if $\a_{n_s}\pm\a_i$ is not a root or a
short root respectively such that
\begin{equation}
d_{v_{n_s},v_{\pm i}}=m'_{\pm i}e_{t'_{\pm i}}\andd v_{n_s}*v_{\pm
i}=m_{\pm i}v_{ t_{ i}}.\label{newm2}
\end{equation}
We set $$\fD_m:=\{t'_{\pm i},t_i,m_{\pm i},m'_{\pm i}\mid n_\ell+1\leq i\leq n_s-1\}.$$ We drew the attention of the
 reader to the fact that the  data appearing in  $\fD_m$ are derived  from the structure of the $\gg$-module  $\v.$     We next take $C$ to be the Cartan matrix of $R$ with respect to $\D$ and call
 \begin{equation}
 \label{str}
 \fD:=\{C\}\cup \fD_m\cup \fD_c
 \end{equation}
  the {\it structural data} associated to  $\LL.$

Now let $x_1,\ldots,x_p\in\hh$ and $y_1,\ldots,y_q$ be  some
 non-zero root vectors of $\gg$. Let $j_1,\ldots,j_p\in J_\ell$ and
set $t^\sg=t_{j_1}^\pm\cdots t_{j_p}^\pm$. Now it follows from
(\ref{ubracket}) that, for $1\leq i\leq m$,
\begin{equation}\label{findgen1}
\begin{array}{l}
\hbox{\small$[y_1\ot 1,\ldots,y_q\ot 1,x_1\ot
t_{j_1}^\pm,\ldots,x_p\ot t_{j_p}^\pm, e_n\ot
1]=[y_1,\ldots,y_q,x_1,\ldots,x_p,e_n]\ot t^\sg,$}\\
\hbox{\small$[y_1\ot 1,\ldots,y_q\ot 1,x_1\ot
t_{j_1}^\pm,\ldots,x_p\ot t_{j_p}^\pm,v_{n_s}\ot w_i]=y_1\ldots
y_q x_1 \ldots x_p v_{n_s}\ot t^\sg w_i.$}
\end{array}
\end{equation}

Next suppose that $\{e_i,f_i,h_i\}$ is a set of Chevalley
generators for $\gg.$ Then using  (\ref{findgen1}) together with
(\ref{decom1}), we see that $\mathfrak{A}$ as a Lie algebra is
generated by
\begin{equation}\label{genset}
\{e_i\ot 1,f_i\ot 1,h_i\ot 1,h_i\ot t_j^\pm,v\ot w_t\mid i\in
J_\ell,\;j\in J_\nu,\;t\in J_m \}.
\end{equation}
For  $n_\ell+1\leq i\leq n_s-1$ and $ r,s\in J_m,$ take
$\sg_{r,s},\sg'_{r,s}\in\zn,$ $m_{\pm i},m'_{\pm
i},a_{r,s},a'_{r,s}\in\bbbc,$ $ t'_{\pm i}\in J_n,$ $t_{r,s}\in
J_m,$ $n_\ell+1\leq t_i\leq n_s$ and $m_{\pm i},m'_{\pm
i}\in\bbbc$ to be  defined as in (\ref{newm1}) and  (\ref{newm2}).
Using(\ref{ubracket}) and   (\ref{fin}), one sees that
\begin{equation}\label{mrel1}
[v_{n_s}\ot w_r,v_{\pm i}\ot w_s]=(m'_{\pm i} a'_{r,s}e_{t'_{\pm
i}}\ot t^{\sg'_{r,s}})+(m_{\pm i} a_{r,s}v_{t_i}\ot
t^{\sg_{r,s}}w_{t_{r,s}}).
\end{equation}
Following  \cite{You}, we call  this kind of relations {\it basic
short part relations.}

Next let $\a\in R^\times$, $\b\in R_{sh}$, $x_\a\in\gg_\a$ and
$v_\b\in\v_\b$. Then (\ref{ubracket}) implies  that for $ i,r\in
J_\ell,$ $ j\in J_\nu$ and $ s\in J_m,$ the following relations
are satisfied in  $\mathfrak{A},$ the universal covering algebra
of $\LL:$
$$
\begin{array}{l}
[\a(h_r)(h_i\ot t_j^{\pm1})-\a(h_i)(h_r\ot t_j^{\pm1}),x_\a\ot1]=0,\vspace{1mm}\\
\hspace{0cm}[\b(h_r)(h_i\ot t_j^{\pm1})-\b(h_i)(h_r\ot
t_j^{\pm1}),v_\b\ot w_s]=0,
\end{array}$$  to which we
refer  as  {\it  quasi-diagonal relations}, and also
$$
\begin{array}{l}
\;[{h}_r\ot t_j,h_{i}\ot t_j^{-1},x\ot1]=\a(h_i)\a(h_r)x\ot1,\vspace{1mm}\\
\;[{h}_r\ot t_j,h_{i}\ot t_j^{-1},y\ot w_s]=\b(h_i)\b (h_r)y\ot
w_s,
\end{array}$$ that we refer to as {\it  cancelling
relations.}

\section{ A generic presentation}
From now on we fix a centerless Lie torus $\LL$ of type $X\not= A, C, BC$ and we let $\fA$ be its universal covering
algebra. Let $\mathfrak{D}$
be the structural data associated to $\LL$. In this section starting from $\fD$, we introduce a presented Lie algebra
$\Lt:=\Lt(\fD)$, called the {\it presented Lie algebra associated to $\LL$ (or $\fD$)}. The main objective of this work is to show that $\Lt$ is isomorphic to $\fA$.

We use the same notation as in Section \ref{Terminology}. In particular, we let $\gg$ be the finite dimensional
simple Lie
algebra of rank $\ell$ corresponding to the Cartan matrix $C=(c_{i,j})$. Let $\hh$ be a fixed Cartan subalgebra of $\gg$
and
$\fm$ be the killing form $\fm$. We recall that $R$ is the root system
of $\gg$ and $\D$ is a base of $R.$ Also $\v$ is
an irreducible $\gg$-module whose highest weight is a short root.
We keep the same arrangement for roots as in  (\ref{fin2}). Throughout  this work,
by
$\la \a,\b\ra,$ for roots $\a,\b\in R$ with $\b\neq 0,$  we mean $2(\a,\b)/(\b,\b).$
Let $m$ be a
non-negative integer which we take it to be zero if $R$ is
simply-laced. Whenever we use expressions containing a letter with
subscripts going through $\{1,\ldots,m\},$ we understand that we
are in the case $m\neq0.$ Now for a non-negative integer $\nu$,
let $\Lt$ be the Lie algebra defined by $3\ell+m+\ell\nu$
generators
\begin{equation}
\{e_i,f_i,h_i,h_{i,a}^{\pm},v^r\mid i\in J_\ell,\; r\in J_m,\;a\in
J_\nu\},
 \label{3.1'}
\end{equation}
subject to the relations (we  collect our relations, depending on
their natures, in groups (R1)-(R9) below and  give a name to some
of these groups, based on the role which play):

\vspace{.5cm} \noindent\hbox{\bf (R1)}
$
\begin{array}{c} \hbox{\underline{\bf Serre's relations}:}\vspace{1mm}\\
\;[h_i,h_j]=0,\;[e_i,f_j]=\d_{i,j}h_i,\;[h_i,e_j]=c_{j,i}e_j,\;[h_i,f_j]=-c_{j,i}f_j,
\vspace{2mm}\\
(\ad e_i)^{-c_{j,i}+1}(e_j)=0,\;(\ad
f_i)^{-c_{j,i}+1}(f_j)=0,\vspace{2mm}\\
(i,j\in J_\ell).
\end{array}$

\medskip

\noindent$(\hbox{\bf R2})\begin{array}{c} \hbox{\bf
\underline{Highest short weight relations}:}\vspace{1mm}\\
\;[e_i,v^r]=0,\;[h_i,v^r]=\la\a_{n_s},\a_i\ra v^r,
\;[f_i^{\la\a_{n_s},\a_i\ra+1},v^r]=0,\vspace{1mm}\\
 i\in J_\ell,\; r\in j_m.
\end{array}$

Since $\{e_i,f_i,h_i\}_{i=1}^ \ell$ satisfies Serre's relations,
the subalgebra of $\Lt$ generated by $3\ell$ elements
$\{e_i,f_i,h_i\mid i\in J_\ell\}$ is a finite dimensional simple
Lie algebra of the same type of $R$ [H, Theorem 18.3]. So we
identify this subalgebra with the Lie algebra $\gg$ as in the
Subsection 1.1 with Cartan subalgebra $\hh=\oplus_{i\in
J_\ell}\bbbc h_i$, and the corresponding root system $R$. This
then also allows us to identify the $3\ell$ generators $e_i, f_i,
h_i$, $i\in J_\ell$ here with the corresponding elements as in
Subsection 1.1. Next let $\ell+1\leq i\leq n$ and  use  Lemma \ref{rev1}($i$) to fix $j_1,\ldots,j_{n_i}\in J_{n}$ such that
\begin{equation}\label{fix1}
e_i:=[f_{j_1},\ldots,f_{j_{n_i}},e_n]
\end{equation}
is a nonzero element of $\gg_{\a_i}.$ For $a\in J_\nu$ define
\begin{equation}
e_{i,a}^\pm:=(1/2)[f_{j_1},\ldots,f_{j_{n_i}},h_{n,a}^{\pm},e_n].
 \label{*''}
\end{equation}

Next we  set
\begin{equation}\begin{array}{c}
H:=\hbox{span}_{\bbbc}\{h_{i,a}^\pm\mid i\in J_\ell,\;a\in J_\nu\},\vspace{1mm}\\
Z_h:=\hbox{span}_\bbbc\{[h_{i,a}^\pm,h_{j,b}^{\pm}]\mid i,j\in J_\ell,\;a,b\in J_\nu\},\vspace{1mm}\\
S_g:=\hbox{span}_\bbbc\{e_i,f_i,h_i,h_{i,a}^\pm,v^r\mid i\in
J_\ell,\;a\in J_\nu,\;r\in J_m\}.
\end{array}
\end{equation}

We also  note that using (R2), one concludes that, for $r\in J_m,$ the
$\gg-$submodule $\v^r$ of $\Lt$ generated by $v^r$ is an
irreducible $\gg$-module whose highest weight is the highest short
root $\a_{n_s}$ of $R$ (see \cite[Theorem 21.4]{H}). Considering (\ref{my1}), we set
\begin{equation}\label{my2}
v^r_t:=[f_{j^t_1},\ldots,f_{j^t_{n_t}},v^r]\andd
v^r_{-t}=[f_{k^t_1},\ldots,f_{k^t_{n'_t}},v^r].
\end{equation}
Now we have the following weight space decomposition
\begin{equation*}
\v^r=\bigoplus_{t=n_\ell+1}^{n_s} (\bbbc v_{\pm t}^r)\oplus\v^r_0,\label{my6}
\end{equation*}
where  $\v^r_0$ is the corresponding zero weight space.
\medskip

  For $h=\sum_{i=1}^\ell r_ih_i\in\hh$ and $a\in J_\nu,$ take
\begin{equation}
\label{rev3}
\displaystyle{h_a^{\pm}:=\sum_{i=1}^\ell r_ih_{i,a}^\pm}
\end{equation}
and for  $\ell+1\leq t\leq n,$  set
$$h_{t,a}^\pm:=(h_t)_a^\pm.$$

Let   $\sg=(m_1,\ldots,m_\nu)\in\zn$.  We call $|\sg|:=\sum_{t=1}^{\nu}|m_t|$ the {\it norm} of $\sg$.
 Let $\sg\neq 0$ and $1\leq i_1<\cdots<i_p\leq \nu$ be all $i_j\in J_\nu$ for which $m_{i_j}\not=0$.
 Then $|\sg|=\sum_{j=1}^{p}|m_{i_j}|$. For $1\leq t\leq n,$ we set
 \begin{equation}\label{fin3}
 \mathfrak{b}^t_\sg:=(\underbrace{h_{t,i_1}^{\sgn(m_{i_1})},\ldots, h_{t,i_1}^{\sgn(m_{i_1})}}_{|m_{i_1}|},
 \ldots,\underbrace{h_{t,i_p}^{\sgn(m_{i_p})},\ldots,h_{t,i_p}^{\sgn(m_{i_p})}}_{|m_{i_p}|})
\end{equation}
where $sgn(m)$ for $m\in\bbbz$ is the sign of $m.$ In fact
$\mathfrak{b}^t_\sg=(b^t_1,\ldots,b^t_{|\sg|})$ where  for $k\in
J_{|\sg|},$ $\begin{array}{c} b^t_k:=\left\{\begin{array}{ll}
h_{t,i_a}^+ &\hbox{if}\;\;m_{i_a}>0\\
h_{t,i_a}^{-} &\hbox{if}\;\;m_{i_a}<0,\\
\end{array}\right.\end{array}$
in which $a$ is the unique element of $\{1,\ldots,p\}$ with
 $$1+\sum_{j=1}^{a-1}|m_{i_j}|\leq i\leq\sum_{j=1}^{a}|m_{i_j}|.$$
We also  set $\fb^t_0:=(h_{t,1}^-,h_{t,1}^+).$ We call
$\fb_\sg:=\fb_\sg^n$ the {\it norm-tuple} of $ \sg$.

\begin{con} For $\sg\in\zn,$ we denote the norm-tuple of $\sg$ by $\fb_\sg=(\fb^\sg_1,\ldots,\fb^\sg_{|\sg|})$ and by
$(\fb_1,\ldots,\fb_{|\sg|})$ if there is no confusion. For $h\in
\hh,$ if $i\in J_{|\sg|},$ define $h^{i,\sg}$ to be $h_{a}^\pm$
(see (\ref{rev3})) if $\fb_i^\sg=h_{n,a}^\pm$ for some $a\in
J_\nu.$ Also with the abuse of notation we write
$$
[r \fb^t_\sg,x]:=[r \fb^t_1,\ldots,r \fb^t_{|\sg|},x],\qquad (t\in
J_n,x\in \Lt,r\in\bbbc).
$$

\label{my8}
\end{con}

We now introduce some more relations:
\medskip

\noindent $(\hbox{\bf R3})\;\;\;
[e_{i,a}^\pm,f_i]=h_{i,a}^\pm;\;\;i\in J_\ell,\;\;a\in J_\nu.$

\medskip

\noindent $(\hbox{\bf R4})\begin{array}{c}
[H,\hh]=\{0\},\;\;[S_g,Z_h ]=\{0\}.
\end{array}$
\medskip

\noindent $(\hbox{\bf R5})\begin{array}{c} [\la\a_{t},\a_j\ra
h_{i,a}^\pm-\la \a_t,\a_i\ra h_{j,a}^\pm,\bbbc e_t+\bbbc f_t]=0;\;
a\in J_\nu,\;i,j\in J_\ell,\;t\in J_n. \end{array}$

\medskip

\noindent $(\hbox{\bf R6})\begin{array}{c}
\hbox{\bf \underline{Canceling relations}:}\vspace{1mm}\\
\;[{h}_{i,a}^{-},h_{j,a}^+,e_t]=\la \a_t,\a_i\ra\la\a_{t},\a_j\ra
e_t,\;
\;[{h}_{i,a}^{-},h_{j,a}^+,f_t]= \la \a_t,\a_i\ra\la\a_{t},\a_j\ra f_t\vspace{1mm}\\
a\in J_\nu,\;i,j,t\in J_n. \end{array}$

\medskip

We consider (\ref{my2}) to define our next two sets of relations:
\medskip

\noindent $(\hbox{\bf R7})\begin{array}{c} [\la\a_{t},\a_j\ra
h_{i,a}^\pm-\la \a_t,\a_i\ra h_{j,a}^\pm,\bbbc v^r_{\pm t}]=0,\vspace{1mm}\\
 a\in J_\nu,\;r\in J_m,\;i,j\in J_\ell,\; n_\ell+1\leq t\leq n_s.\end{array}$

\medskip

\noindent $(\hbox{\bf R8})\begin{array}{c}
\hbox{\bf \underline{Canceling module relations:}}\vspace{1mm}\\
\;[{h}_{i,a}^{-},h_{j,a}^+,v^r_{\pm t}]=\la \a_t,\a_i\ra \la
\a_t,\a_j\ra v^r_{\pm t},\\ a\in J_\nu,\;i,j\in J_\ell,\; r\in
J_m,\;n_\ell+1\leq t\leq n_s.
\end{array}$

\medskip

Finally we recall (\ref{mrel1}) and  (\ref{fin3}) to state our
last set of relations as follows:

\noindent $(\hbox{\bf R9})$ $\begin{array}{c} \hbox{\bf \underline{Basic short part relations:}}\vspace{1mm}\\
\;[v^r,v^s_{\pm i}]=m'_{\pm i}a'_{r,s}[(1/2)\fb^{t'_{\pm i}}_{\sg'_{r,s}},e_{t'_{\pm i}}]+
m_{\pm i}a_{r,s}[(1/2)\fb^{t_{ i}}_{\sg_{r,s}},v_{t_i}^{t_{r,s}}]\\
r,s\in J_{m},n_\ell+1\leq i\leq n_s-1.
\end{array}$

\begin{rem}
{\rm One knows that $\la \a_i,\a_j\ra=\a_i(h_j)$ for $i,j\in J_n,$ so
we may use  relations  (R1)-(R9) with   $\a_i(h_j)$ instead of $\la
\a_i,\a_j\ra.$}
\end{rem}

\begin{DEF}\label{presentedliealgebra}
{Starting from a centerless Lie torus $\LL$ with the structural data $\fD$, we call the Lie algebra $\Lt$ defined by generators
(\ref{3.1'}) and relations (R1)-(R9), the {\it presented Lie algebra associated to $\LL$ (or $\fD$)}.}
\end{DEF}

We are now ready to state our main theorem.

\noindent
{\bf Main Theorem:} Let $\LL$ be a centerless Lie torus of type $X\not= A, C, BC$ with the universal covering $\fA$, and associated presented Lie algebra $\Lt$. Then $\Lt\cong\fA$.
 In
particular  $\mathfrak{A}$ is a finitely presented Lie algebra.

\vspace{5mm}
In the remaining part of this section we establish several results which are needed prior to the proof of the main theorem.
 The proof of the main theorem will be presented in Section \ref{maintheorem}. An outline of our plan for completing this work is as follows.
We first show that $\Lt$ as a $\gg$-module ($\gg$ is a subalgebra of $\Lt$) is a direct
sum of irreducible $\gg$-modules whose highest weights belong to
$R.$ We then use this to show that  $\Lt$ is an $R$-graded
Lie algebra. Using the fact that both $\Lt$ and $\fA$ are $R$-graded Lie algebras and appealing   the structure of $R$-graded Lie
algebras we prove that $\Lt$ is a central extension of $\fA$ and finally  we prove that $\fA$ is isomorphic to $\Lt.$

 We recall  that
$\mathfrak{A}$ is generated by (\ref{genset}) and so one can
easily get from Subsection 1.2  that there exists a Lie algebra
epimorphism from $\Lt$ to $\mathfrak{A}$ as follows:
\begin{equation}
\begin{array}{c}
\psi:\Lt\longrightarrow\mathfrak{A}\\
e_i\mapsto e_i\ot1,\;f_i\mapsto f_i\ot1,\;h_i\mapsto h_i\ot1,\;
h_{i,a}^\pm\mapsto h_i\ot t_a^{\pm1},\;
 v^r \mapsto v_{n_s}\ot w_r,\\
 a\in J_\nu,i\in J_\ell, r\in J_m.\end{array} \label{epi}
\end{equation}

From now on we fix $\sg=(m_1,\ldots,m_\nu)\in\zn.$  We associate to $\sg$ the following   elements of $\Lt$
\begin{equation}\label{3.3}
\es:=[\frac{1}{2}\fb_\sg,e_n]\andd v_{\sg}^r:=[\fb_\sg,v^r];\;\;
r\in J_m \hbox{ (see Convention \ref{my8})}.
\end{equation}
Using (R6) and (R8), we have
$$
 e_0=e_n\andd v_0^r=v^r.$$

Now considering $\Lt$ as a $\gg$-module, we set
\begin{equation}
\begin{array}{l}
\gg_\sg:=\hbox{the $\gg$-submodule of $\Lt$ generated by $\es$ and}\\\\
\v_\sg^r:=\hbox{ the $\gg$-submodule
of $\Lt$ generated by $v_\sg^r$}.
\end{array}
\label{3.4}
\end{equation}

\begin{pro}
\label{m8} For $\a\in \rcross, $ set
$$\Lt_\a^0:=\left\{\begin{array}{ll}
\gg_{\a}+\displaystyle{\sum_{r=1}^m\v^r_\a}& \hbox{if $\a\in R_{sh}$} \\
\gg_\a & \hbox{if $\a\in R_{lg}.$}
\end{array}
\right.$$ $(i)$ Recall (\ref{rev3}) and let  $a\in J_\nu,$ $x,y\in\hh$ and $t\in J_n,$ then
$$[\a_t(x)y_a^\pm-\a_t(y)x_a^\pm,\Lt^0_{\pm\a_t}]=\{0\},$$ in
particular if $x\in\hh$ is such that $\a_t(x)=0,$ we get
$$[x_a^\pm,\Lt^0_{\pm\a_t}]=\{0\}.$$

$(ii)$ Let $p\in\bbbz^{> 0},$ $a_1,\ldots,a_p\in J_\nu,$ and $t\in J_n,$ suppose that
$x_j,y_j\in \hh$    are  such that $\a_t(y_j)\neq0$ for $j\in J_p.$
   Then for     $e\in \Lt^0_{\pm\a_t},$ we have
 $$\begin{array}{c}
 [(x_1)_{a_1}^\pm,\ldots,(x_p)_{a_p}^\pm,e]=
 \displaystyle{\prod_{j=1}^p(\a_t(x_j)/\a_t(y_j))}
 [(y_1)_{a_1}^\pm,\ldots,(y_p)_{a_p}^\pm,e]
 \end{array}$$
 with the same sign in $(x_j)_{a_j}^\pm$ and $(y_j)_{a_j}^\pm$ for $j\in J_p.$

$(iii)$ Let $p,q\in\bbbz^{> 0},$ $a_1,\ldots,a_p,b_1,\ldots,b_q\in
J_\nu,$ and $t,t'\in J_n$ with $t\not=t'$, suppose that
$x_j,y_j\in \hh,$ $j\in J_p,$  are  such that $\a_t(y_j)\neq0$ and
$\a_{t'}(y_j)=0.$
   Then for   $z_j\in\hh$ with $\a_{t}(z_j)=0,$  $j\in J_q,$ $e\in \Lt^0_{\pm\a_t}$ and
$f\in \Lt^0_{\pm\a_{t'}},$ we have
 $$\begin{array}{c}
 [[(z_1)^\pm_{b_1},\ldots, (z_q)^\pm_{b_q},f],[(x_1)_{a_1}^\pm,\ldots,(x_p)_{a_p}^\pm,e]]=\\
 \displaystyle{\prod_{j=1}^p(\a_t(x_j)/\a_t(y_j))}
 [(z_1)^\pm_{b_1},\ldots, (z_q)^\pm_{b_q},(y_1)_{a_1}^\pm,\ldots,(y_p)_{a_p}^\pm,[f,e]]
 \end{array}$$
 with the same sign in $(x_j)_{a_j}^\pm$ and $(y_j)_{a_j}^\pm$ for $j\in J_p.$
\end{pro}

\pf $(i)$ The first expression    is immediate using  (R5) and
(R7). For the second expression, we note that as $\a_t\neq 0,$ one
finds $y\in \hh$ such that $\a_t(y)\neq 0.$ Now the statement is
fulfilled by considering the first expression.

$(ii)$ Using part ($i$), we have for all $j\in
J_p$,
$$
[(x_j)_{a_j}^\pm,e]=(\a_t(x_j))/(\a_t(y_j))[(y_j)_{a_j}^\pm,e].$$Now we are done using the Jacobi identity together with (R4).

$(iii)$ Using  part $(i),$ we get that  for all $i\in J_p,j\in
J_q$,  $[(y_i)_{a_i}^\pm,f]=0$ and $[(z_j)_{b_j}^\pm,e]=0.$ Now
(R4) together with the Jacobi identity implies that
$[[(z_1)^\pm_{b_1},\ldots, (z_q)^\pm_{b_q},f],(y_i)_{a_i}^\pm]=0,$
$i\in J_p.$ Also using part $(ii),$ we can replace
$[(x_1)_{a_1}^\pm,\ldots,(x_p)_{a_p}^\pm,e]$ with an expression as
in the right hand side of the display appearing in part $(ii).$
Now we are done  using these together with the Jacobi identity.

The following proposition determines, up to isomorphism, the modules $\gg_\sg$ and
$\v^r_\sg$ appearing in (2.13).
\medskip
\begin{pro}\label{y3.1}
(i) $\gs$, is an irreducible finite dimensional $\gg$-module with
highest pair  $(\es,\a_n).$ In fact as $\gg$-modules,
$\gs\cong\gg.$

(ii) For   $r\in J_m,$ $\v_\sg^r$  is an irreducible finite
dimensional $\gg$-module with  highest pair $(v^r_\sg,\a_{n_s}).$
In fact as $\gg$-modules $\v^r_\sg\cong\v.$
\end{pro}

\pf $(i)$ Contemplating (\ref{epi}), one can see that $
\psi(\es)=e_n\ot t^\sg\neq0$. Therefore $\es\neq 0$. Thus by
\cite[Theorem 21.4]{H}, it is enough to show that, for all $i\in
J_\ell$,
\begin{equation}
 \begin{array}{c}[h_i,\es]=\a_n(h_i)\es,\;
 \;\;
 [e_i,\es]=0,\;\;
[f_i^{\a_n(h_i)+1},\es]=0.\\
  \end{array}
\label{con}
  \end{equation}
Fix $i\in J_\ell.$ We first note that the equalities in
(\ref{con}) hold for $\sg=0$ as $e_n $ is a maximal vector in the
$\gg$-module $\gg$. Now let $\sg$ is nonzero with norm-tuple
$\fb_\sg=(\fb_1,\ldots,\fb_{|\sg|}).$ Considering Convention
\ref{my8}, one gets from  (R4) together with (an iterated use of)
the Jacobi identity  that
$$[h_i,\es]=[h_i,[\frac{1}{2}\fb_\sg,e_n]]=-[\frac{1}{2}\fb_\sg,[e_n,h_i]]=\a_n(h_i)e_\sg.
$$

To prove the next two equalities in (\ref{con}), we use induction
on $|\sg|.$  Assume $0\not=\sg$ is nonzero with
$[e_i,e_\sg]=0=[f_i^{\a_n(h_i)+1},e_\sg]$. We are done if we show
that for $a_0\in J_\nu$,
$$
[e_i,[\frac{1}{2}h^{\pm}_{n,a_0},e_\sg]]=0=[f_i^{\a_n(h_i)+1},[\frac{1}{2}h^\pm_{n,a_0},e_\sg]].
$$
Since $\a_n$ and $\a_i$ are not proportional, there exists $x\in
\hh$ such that $\a_n(x)\neq 0$ and $\a_i(x)=0.$ Therefore using
(R4) and Proposition \ref{m8}$(i),$ we have
$$
[\frac{1}{2}h_{n,a_0}^\pm,\es]\in\bbbc[x_{a_0}^\pm,\es]\andd
[x_{a_0}^\pm,e_i]=0.
$$Now using these together with the Jacobi identity and the induction hypothesis, we get
\begin{eqnarray*}
[e_i,[\frac{1}{2}h_{n,a_0}^\pm,\es]]\in\bbbc[e_i,[x_{a_0}^\pm,\es]]=\bbbc
[[x_{a_0}^\pm,e_i],\es]=0=[[x_{a_0}^\pm,f_i],\es].
\end{eqnarray*}
and
\begin{eqnarray*}
[f_i^{\a_n(h_i)+1},\frac{1}{2}h_{n,a_0}^\pm,\es]\in\bbbc[f_i^{\a_n(h_i)+1},x_{a_0}^\pm,\es]=\bbbc
[x_{a_0}^\pm,f_i^{\a_n(h_i)+1},\es]=0.
\end{eqnarray*}

$(ii)$ Using an  argument analogous to the one in part $(i),$ we are done. \qed
\bigskip

Using Proposition  \ref{y3.1}($i$), one concludes that there exists a
$\gg$-module isomorphism $\vp_\sg:\gg\longrightarrow\gs$  mapping
$e_n$ to $\es$. Therefore   $\gg_\sg$ admits a weight space decomposition
$\gg_\sg=\sum_{\gamma\in R}(\gg_\sg)_\gamma,$ where
$(\gg_\sg)_\gamma=\varphi_\sg(\gg_\gamma),$ $\gamma\in R.$ For $a\in J_\nu$, we set
\begin{equation}
\sg_a^{\pm}:=(0,\ldots,0,\underbrace{\pm1}_{a \hbox{\tiny
-th}},0,\ldots,0)\andd\vp_a^{\pm}:=\vp_{\sg_a^{\pm}}\label{n4}.
\end{equation}
From (\ref{3.3}) we have
$\varphi_a^\pm(e_n)=e_{\sg_a^\pm}=(1/2)[h_{n,a}^\pm,e_n]$. Now
considering (\ref{fix1}), (\ref{*''}) and (R3) we have, for $i\in
J_\ell$ and $a\in J_\nu$,
\begin{equation}\label{ma3}
e_{i,a}^\pm=\vp_a^\pm(e_i)\andd h_{i,a}^\pm=\vp_a^\pm(h_i).
\end{equation}

 Now set, for $t\in J_n$, $a\in J_\nu$ and $\ell+1\leq j\leq n$,
\begin{equation}
\begin{array}{c}
e_{t,\sg}:=\vp_\sg(e_t),\;\;f_{t,\sg}:=\vp_\sg(f_t),\;\;h_{t,\sg}:=\vp_\sg(h_t),\vspace{2mm}\\
\;\;f_{t,a}^\pm:=\vp_a^\pm(f_t),\;\; e_{j,a}^\pm:=\vp_a^\pm(e_j ),\;\;h_{j,a}:=\vp_a^\pm(h_j)\\
\end{array}
\label{y9}
\end{equation}
and finally set
\begin{equation}\label{cen}
\z:=\hbox{span}_\bbbc\{[h_{i,\sg},h_{i,\tau}]\mid i\in
J_\ell,\;\sg,\tau\in\zn\}.
\end{equation}
%
We next note that as  $\vp_\sg,$ $\sg\in\zn,$  is a $\gg$-module
isomorphism and  $[e_i,f_j]=\d_{i,j}h_i,$  $i,j\in J_\ell,$ we get
for $a\in J_\nu$ and  $i,j\in J_\ell$,
\begin{equation}
[e_{i,\sg},f_j]=-[f_{i,\sg},e_j]= \d_{i,j}h_{i,\sg}.\label{m1}
\end{equation}


Now let $r\in J_m,$ then  Proposition \ref{y3.1}($ii$) guarantees
the existence of  a $\gg$-module isomorphism
\begin{equation}
\psi_\sg^r:\v\longrightarrow\v^r_\sg\;\;\hbox{such that}\;\;{\bf
v}=v_{n_s}\mapsto\vrs\hbox{ (see (\ref{3.3}))}. \label{iso}
\end{equation}
Set $ \psi_r:=\psi_0^r.$ Since $\psi_r({\bf v})=v^r,$
(\ref{my1}) and (\ref{my2}) imply that
\begin{equation}
v_{\pm t}^r=\psi_r(v_{\pm t});\;\;\;\;n_\ell+1\leq t\leq n_s. \label{3.33}
\end{equation}
Also as $\v^r_\sg$ is an irreducible finite dimensional $\gg$-module of
highest weight $\a_{n_s},$ we have
\begin{equation}
\begin{array}{c}
\v^r_\sg=(\v^r_\sg)_0\op\bigoplus_{t=n_\ell+1}^{n_s}(\v_\sg^r)_{\pm t} \hbox{ where }\vspace{2mm}\\
(\v_\sg^r)_{\pm t}:=(\v_\sg^r)_{\pm \a_t} =\bbbc v_{\sg,\pm t}^r \hbox{ with }
v_{\sg,\pm t}^r:=\psi_\sg^r(v_{\pm t});\;\;\;\;n_\ell+1\leq t\leq n_s.
\end{array}
\label{ma5}
\end{equation}


\begin{pro}
Let  $t\in J_n$ and $\a=\pm\a_t.$ Suppose for $j\in J_{|\sg|},$ $x_j\in\hh$ is such
that $\a(x_j)\neq 0.$  For $j\in J_{|\sg|}$, recall  Convention \ref{my8} and  set
$c_j:=(x_j)^{j,\sg}.$   Then for
$e\in\Lt_{\a}^0$ (see Proposition \ref{m8}), we have
$$\frac{1}{\prod_{j=1}^{|\sg|}(\pm\a(x_j))}[c_1,\ldots,c_{|\sg|},e]=\left\{\begin{array}{rl}
\vp_\sg(e)& \hbox{if } e\in\gg_{\a}\\
\psi_\sg^r(e)& \hbox{if } e\in\v^r_{\a}\hbox{ for some $r\in J_m.$}
\end{array}\right.$$
 \label{y3.17}
\end{pro}

\pf Using Proposition \ref{m8}($ii$), it is enough to show that the right hand side of the expression  equals to
$[y_1^{1,\sg},\ldots,y_{|\sg|}^{|\sg|,\sg},e]$
for some $y_1,\ldots,y_{|\sg|}\in \hh$ with $\a(y_i)=1,$ $i\in J_{|\sg|}.$ Take $$\left\{\begin{array}{ll} v:=e_n\andd \lam:=\a_n&\hbox{if } e\in\gg_{\a}\\
v:=v^r \andd \lam:=\a_{n_s}&\hbox{if } e\in\v^r_{\a}\hbox{ for
some $r\in J_m.$}
\end{array}\right.$$
Using Lemma \ref{rev1}($i$), we may suppose
$e=[f_{i_{p}},\ldots,f_{i_1},v]$ where  $p\in\bbbn,$
$i_1,\dots,i_{p}\in J_n$ and for  $q\in J_p,$ $\a_{i_q}\neq
\pm(\lam-\sum_{r=1}^{q-1}\a_{i_r}).$

For $1\leq q\leq p,$ $\a_{i_q}\neq
\pm(\lam-\sum_{r=1}^{q-1}\a_{i_r}),$ so there is $x_q\in\hh,$ such
that $\a_{i_q}(x_q)=0$ and
$(\lam-\sum_{r=1}^{q-1}\a_{i_r})(x_q)=1.$ Now we recall Convention
\ref{my8} and note that Proposition \ref{ma5}($i$) implies
that $[f_{i_q},x_q^{i,\sg}]=0,$ $i\in J_{|\sg|}.$ So contemplating
Proposition \ref{m8}($ii$) together with the Jacobi identity, we
get {\small\begin{eqnarray*}
[f_{i_q},x_q^{1,\sg},\ldots,x_q^{|\sg|,\sg},f_{i_{q-1}},\ldots,f_{i_1},v]&=&[x_q^{1,\sg},\ldots,x_q^{|\sg|,\sg},
f_{i_q},f_{i_{q-1}},\ldots,f_{i_1},v]\\
&\stackrel{\hbox{\tiny if $q\neq p
$}}{=}&[x_{q+1}^{1,\sg},\ldots,x_{q+1}^{|\sg|,\sg},f_{i_q},f_{i_{q-1}},\ldots,f_{i_1},v].
\end{eqnarray*}}
Using this repeatedly, one has
\begin{eqnarray*}
[x_{p}^{1,\sg},\ldots,x_{p}^{|\sg|,\sg},e]&=&[x_{p}^{1,\sg},\ldots,x_{p}^{|\sg|,\sg},f_{i_{p}},\ldots,f_{i_1},v]\\
&=&[f_{i_{p}},\ldots,f_{i_1},x_{1}^{1,\sg},\ldots,x_{1}^{|\sg|,\sg},v]\\
&=&\left\{\begin{array}{rl}
\;[f_{i_{p}},\ldots,f_{i_1},e_\sg]& \hbox{if } e\in\gg_{\a}\\
\;[f_{i_{p}},\ldots,f_{i_1},\vrs]& \hbox{if } e\in\v^r_{\a}\hbox{ for some $r\in J_m.$}
\end{array}\right.
\end{eqnarray*}
Now we are done as $\varphi_\sg$ and $\psi_\sg^r$ are $\gg$-module  homomorphisms.\qed

\bigskip

 Recall from Section
\ref{Terminology} that we have fixed an specific order for (simple)
roots, in terms of their lengths. Accordingly, for type $F_4$, we
consider the following fundamental system
\begin{equation}\label{fin7}
\{\a_1=\ve_2-\ve_3,\a_2=\ve_3-\ve_4,\a_3=\ve_4,\a_4=\frac{1}{2}(\ve_1-\ve_2-\ve_3-\ve_4)\},
\end{equation} where as usual $\ve_i$'s are the standard orthogonal basis for $\bbbr^4$.
Using the module theory for type $F_4$, we may find fix  complex numbers
$a,b,a',b',a_3'',b_3'',a_4'',b_4''$ satisfying
\begin{equation}
\begin{array}{c}
(a_3''a'+b_3''a)b=-(a_3''b'+b_3''b)a,(a_4''a'+b_4''a)b'=-(a_4''b'+b_4''b)a',\\
t(a'f_3\cdot v_3+a f_4\cdot v_{4},b' f_3\cdot v_{3}+b f_4\cdot v_{4})=1,\\
\;(a'f_3\cdot v_{3}+a f_4\cdot v_{4})*(b'f_3\cdot v_{3}+b,f_4\cdot v_{4})=0,\\
e_{i}\cdot(f_{3}\cdot v_{3})=a_i''v_{i},e_{i}\cdot(f_{4}\cdot
v_{4})=b_i''v_{i}.
\end{array}
\label{ch9}
\end{equation}
where $3\leq i\leq 4$ and $t$ is the trace form introduced in
Section \ref{Terminology} (see Subsection \ref{exceptional}). We then define
$$D^{r,s}_{\sg,\tau}:=\left\{\begin{array}{ll}
\;[a'[f_3,v^r_{\sg,3}]+a [f_4,v_{\sg,4}^r],b'[f_3,v^s_{\tau,3}]+
b[f_4,v_{\tau,4}^s]]& \hbox{for type $F_4$},\\
 \;[[f_\ell,v^r_{\sg,\ell}],[f_\ell,v^s_{\tau,\ell}]]&\hbox{otherwise,}
\end{array}\right.$$
where $r,s\in J_m$ and $\sg,\tau\in\zn.$ Now  set
\begin{equation}\label{def-d}
\dd:=\hbox{span}_{\bbbc}\{D_{\sg,\tau}^{r,s}\mid r,s\in
J_m,\;\;\sg,\tau\in\zn \}.
\end{equation}

In Propositions \ref{all2} and \ref{alltogether}, we prove several crucial
relations among the elements of $\Lt$,
which reveals the algebra structure between parts
$\gg_\sg$, $\v^r_\sg$, $\dd$ and $\mathcal Z$ (see (\ref{3.4}), (\ref{cen}) and
(\ref{def-d})) of $\Lt$. Since the proofs of these two propositions are quite
technical and a bit long, we found it more convenient to records their statements
here and postpone their proofs until Section \ref{postpone}.

\begin{pro}\label{all2}
Let  $r\in J_m,$ $a\in J_\nu$, $i\not=j\in J_\ell$ and
$\sg,\tau\in\zn.$ Set  $\tau^{\pm a}=\tau\pm\sg_a,$ then
considering (\ref{y9}), we have

(i) $[e_{i,\tau},e_{i,\sg}]=0=[f_{i,\tau},f_{i,\sg}],$

$(ii)$ $[e_{i,\sg},h_{j,\tau}]=-\a_i(h_j)e_{i,\sg+\tau}\andd
[f_{i,\sg},h_{j,\tau}]=\a_i(h_j)f_{i,\sg+\tau},$

$(iii)$
$[h_{i,\tau},e_{i,a}^\pm]=[h_{i,a}^\pm,e_{i,\tau}]=2e_{i,\tau^{\pm_a}}\hbox{
and }
[h_{i,\tau},f_{i,a}^\pm]=[h_{i,a}^\pm,f_{i,\tau}]=-2f_{i,\tau^{\pm_a}},$

$(iv)$
$[h_{j,a}^\pm,h_{i,\sg}]=({\a_i(h_j)}/{2})[h_{i,a}^\pm,h_{i,\sg}],$

(v) $[e_{i,a}^\pm,f_{i,\tau}]=\frac{1}{2}[h_{i,a}^\pm,h_{i,\tau}]+h_{i,\tau^{\pm a}},$

(vi) $[v^r,h_{i,\sg}]=-\a_{n_s}(h_i)v^r_\sg,$ in particular
$$[h_{i,\sg},h_{i,\tau},v^r]=[h_{i,\tau},h_{i,\sg},v^r].$$
 \end{pro}

\begin{pro}
\label{alltogether} For $\a\in \hh^\star\setminus \{0\},$
$\sg\in\zn$ and $r\in J_m,$ define $(\gs)_\a$ and $(\v^r_\sg)_\a$
to be zero if $\a$ is not a root or a short root respectively. Let
$\sg,\tau\in \zn,$ $r,s\in J_m$ and $t\in J_n,$ then  we have
\medskip

\noindent (i)
 $\displaystyle{[H,\sum_{t'=1}^{n}\sum_{\sg\in\zn}(\gs)_{\pm\a_{t'}}
]\subseteq \sum_{t'=1}^{n}\sum_{\sg\in\zn}(\gs)_{\pm\a_{t'}}}.$

\medskip

\noindent (ii) $[v^r,(\gs)_{\pm\a_t}]\subseteq
(\v^r_\sg)_{\pm\a_t+\a_{n_s}}$ and  $[v^r,(\gs)_{0}]\subseteq
(\v^r_\sg)_{\a_{n_s}}.$

\medskip

\noindent (iii)  $[H,(\v^r_\sg)_0]=\{0\}.$

\noindent (iv)
 $\displaystyle{[H,\sum_{p=1}^{m}\sum_{\d\in\zn}\v^p_\d
]\subseteq \sum_{p=1}^{m}\sum_{\d\in\zn}\v^p_\d}.$

\medskip

\noindent (v)  $[v^s,v^r_\sg]=0.$

\noindent (vi) For  $\a\in R_{sh}\cup\{0\}\setminus\{\pm\a_{n_s}\},$
we have
$$[v^s,(\v^r_\sg)_\a]\sub\sum_{\d\in\zn}(\gg_\d)_{\a+\a_{n_s}}+\sum_{1\leq
p\leq m}\sum_{\d\in\zn}(\v^p_\d)_{\a+\a_{n_s}}.$$  Moreover
only the first (second)
summation on the right hand side happens if $r=s$
($r\not=s$).

\medskip

\noindent (vii)
$[v^s,v^r_{\sg,-n_s}]\sub\sum_{\d\in\zn}(\gg_\d)_{0}+\sum_{1\leq
k\leq m}\sum_{\d\in\zn}(\v^k_\d)_{0}+\dd.$

\medskip

\noindent (viii)
$\displaystyle{[v^r,\dd]\sub\sum_{k=1}^m\sum_{\d\in\zn}\v^k_\d.}$

\noindent (ix) $[H,\dd]=\{0\}.$

\medskip

\noindent (x) $\dd$ is a trivial $\gg$-module.

\medskip

\noindent (xi)
Considering (\ref{cen}), we have
$\z\sub Z(\Lt).$
\end{pro}

\section{ Main theorem}\label{maintheorem}
In this section we state and prove the main theorem of this work. The notation and terminology will be as in the previous sections.

\begin{thm}[{\bf Main Theorem}]\label{main} Let $\LL$ be a centerless Lie torus of type $X\not= A, C, BC$ with the universal covering $\fA$, and associated presented Lie algebra $\Lt$. Then $\Lt\cong\fA$.
 In
particular  $\mathfrak{A}$ is a finitely presented Lie algebra.
\end{thm}

\pf We proceed with the proof in a few steps.
\medskip

\underline{{\bf Step1.}}\label{last1} Considering (\ref{ma3}) and
 Propositions \ref{y3.1}, \ref{all2}($iv$) together with  Proposition
\ref{alltogether} and using  the same argument as in \cite[Theorem
2.2 and  Proposition  2.11]{You}, we get that
\begin{equation*}\label{decom}\Lt=\sum_{\sg\in\zn}\gs+\sum_{\sg\in\zn}\sum_{r\in
J_m}\v^r_{\sg}+\dd+\z\andd Z(\Lt)\sub\dd+\z.
\end{equation*}
and that $\Lt$ is an $R$-graded Lie algebra.  Moreover
$\Lt=\op_{\a\in R}\Lt_\a$ where
\begin{equation*}\label{lequ1}
\Lt_\a=\left\{
\begin{array}{ll}
\displaystyle{\sum_{\sg\in\zn}(\gs)_0+
\sum_{\sg\in\zn}\sum_{r\in J_m}(\v^r_\sg)_{0}+\mathcal{D}+\z}& \hbox{if $\a=0,$}\\
\vspace{-2mm}\\
\displaystyle{\sum_{\sg\in\zn}(\gs)_{\a}+
\sum_{\sg\in\zn}\sum_{r\in J_m}(\v^r_\sg)_{\a}} &\hbox{if $\a\in R_{sh},$}\\
\vspace{-2mm}\\
\displaystyle{\sum_{\sg\in\zn}(\gs)_{\a}}&\hbox{if $\a\in
R_{lg}$}.
\end{array}\right.
\end{equation*}
Now we note that as $\Lt$ is an $R$-graded Lie algebra, so is
$\Lt/Z(LL)$ and  as $Z(\Lt)\sub \dd+\z,$ one can identify
$\Lt/Z(\Lt)$ as
\begin{equation}\label{lequ4}\Lt/Z(\Lt)=\sum_{\sg\in\zn}\gs+\sum_{\sg\in\zn}\sum_{r\in
J_m}\v^r_{\sg}+\dd'\;\;\hbox{where }
\dd':=\dd/Z(\Lt).\end{equation}

We shall keep the same  notation for the  images of
$e_i,f_i,h_i,h_{i,j}^\pm$ and $v^r$ in ${\Lt}/{Z(\Lt)}$ and use
$[\cdot,\cdot]^{\bar ~}$ for the Lie bracket on $\Lt/Z(\Lt).$
Now using (\ref{epi}), (\ref{structure}) and Theorem \ref{v1}, we
have an epimorphism
$$\begin{array}{c}
\psi':{\Lt}/{Z(\Lt)}\longrightarrow
\mathfrak{A}/Z(\mathfrak{A})=(\gg\ot A_{[\nu]})\op (\v\ot
A_{[\nu]}^{m})\op \dd_{A_{[\nu]}^{m},A_{[\nu]}^{m}}  \hbox{ with} \\
\vspace{-3mm}\\
e_i\mapsto e_i\ot1,\;f_i\mapsto f_i\ot1,\;h_i\mapsto
h_i\ot1,\;h_{i,j}^\pm\mapsto h_i\ot t_j^{\pm1},\; v^r\mapsto
v_{n_s}\ot w_r.
\end{array}$$
for $i\in J_\ell,\; j\in J_\nu\hbox{ and }r\in J_m.$ Next let
$\a,\b,\a+\b\in R^\times,\;\gamma,\a+\gamma\in
R_{sh},\;\sg,\tau\in\zn$ and $ r\in J_m.$ It  follows using
Proposition
 \ref{y3.17}
 and (\ref{ubracket}) that $\psi'((\gg_\sg)_\a)=\gg_\a\ot
t^\sg$ and $\psi'((\v^r_\sg)_\gamma)=\v_\gamma\ot t^\sg w_r.$
Therefore by (\ref{ubracket}), we have
\begin{equation}
\label{equ-gen1}
\begin{array}{l} \psi'([(\gs)_\a,(\gg_\tau)_\b]^{\bar
~})=[\psi'((\gs)_\a),\psi'((\gg_\tau)_\b)]=\gg_{\a+\b}\ot
t^{\sg+\tau}\neq\{0\},
\\
\\
\psi'([(\gs)_\a,(\v^r_\tau)_\gamma]^{\bar
~})=[\psi'((\gs)_\a),\psi'((\v^r_\tau)_\gamma)]=\v_{\a+\gamma}\ot
t^{\sg+\tau} w_r\neq\{0\},\\
\\
\psi'([(\v^s_\eta)_\zeta,(\v^r_\tau)_\gamma]^{\bar
~})=[\psi'((\v^s_\eta)_\zeta),\psi'((\v^r_\tau)_\gamma)]=[\v_\eta\ot
t^\zeta w_s,\v_\tau\ot t^\gamma w_r].
\end{array}
\end{equation}

\medskip

\underline{{\bf Step2 (Simply-laced types)}}. We define a
$\zn$-gradeding on $\Lt$  as follows. Recalling (\ref{n4}), we set
$$deg(e_i)=deg(f_i)=deg(h_i):=0,\;\; deg(h^\pm_{i,a}):=\pm2\sg_a$$
and note that this defines a $\zn-$grading on the free Lie algebra
generated by the generating set $\{e_i,f_i,h_i,h^\pm_{i,a}\mid
i\in J_\ell,a\in J_\nu\}$, then as the relations (R1)-(R9)   are
generated by homogenous elements, the defined grading induces a
grading on $\Lt$ which  naturally defines a grading on
$\Lt/Z(\Lt).$ We  set
\begin{equation}\label{lequ5}
(\Lt/Z(\Lt))_\a^\sg:=(\Lt/Z(\Lt))_\a\cap(\Lt/Z(\Lt))^\sg;\;\;\a\in
R,\;\sg\in\zn,
\end{equation}
so
\begin{equation}\label{lequ7}
(\Lt/Z(\Lt))_\a^\sg:=(\gs)_\a;\;\;\a\in R,\;\sg\in\zn.
\end{equation}

We next  note that as  $\Lt/Z(\Lt)$  is a centerless $R$-graded
Lie algebra, using Recognition Theorem results in the
existence of a unital commutative associative algebra $A$ such
that $\Lt/Z(\Lt)=\gg\ot A$ and so the multiplication on
$\Lt/Z(\Lt)$ obeys the Lie bracket on $\gg\ot A$ which implies
that
\begin{equation*}\label{lequ6}
(\Lt/Z(\Lt))_\a=\gg_\a\ot A;\;\a\in \rcross.
\end{equation*}
Using the same argument as in \cite[Subsection 2.6]{You}, one can
identify $A$ with $A_{[\nu]}$ such that $h_{i,a}^\pm$ is
identified with $h_i\ot t_a^{\pm 1}.$
This in particular implies that $\psi'$  is an
isomorphism from $\Lt/Z(\Lt)$  to $\fA/Z(\fA).$ Now taking
$\pi_1:\Lt\longrightarrow\Lt/Z(\Lt)$  and
$\pi_2:\fA\longrightarrow\fA/Z(\fA)$  to be the natural projection
maps and considering (\ref{epi}), we get that
$\psi'^{-1}\circ\pi_2\circ\psi=\pi_1.$ Now one concludes that
$\psi:\Lt\longrightarrow\fA$  is a central extension of $\fA,$ but
$\Lt$ is perfect, so by \cite[Proposition I.9.3]{MP}, we get that
$\psi:\Lt\longrightarrow\fA$ is an isomorphism and so
$\Lt\simeq\fA.$
\medskip

\underline{\bf Step3 (Type $G_2$).} We recall  that for type
$G_2,$ there is $0\leq p\leq 3$ such that the corresponding pair of
the $R$-graded Lie algebra $\fA/Z(\fA)$ is $(S,L)$ where
$S:=\bbbz\sg_1\op\cdots\op\bbbz\sg_\nu,$
$L:=3\bbbz\sg_1\op\cdots\op3\bbbz\sg_p\op\bbbz\sg_{p+1}\op\cdots\op\bbbz\sg_\nu$
and that $m=3^p-1.$  For $r\in J_{m},$ define
\begin{equation}\label{final2}
\begin{array}{c}
\sg^r:=\displaystyle{\sum_{i=1}^{p}s_i\sg_i} \;\;\hbox{ where }
r=\displaystyle{\sum_{i=1}^{p}3^{i-1}s_i} \hbox{ for
}s_1,s_2,s_3\in\{0,1,2\}.
\end{array}
\end{equation}
Let  $\sg=(n_1,\ldots,n_\nu)\in\zn,$ for $1\leq i\leq p,$ suppose
that $m_i,r_i\in\bbbz$ are such that $r_i\in\{0,1,2\}$ and
$n_i=3m_i+r_i,$    set
\begin{equation}\label{cormem}
\begin{array}{c}
r_\sg:=\displaystyle{\sum_{i=1}^p 3^{i-1}r_i},\;
\d_\sg:=\displaystyle{\sum_{i=1}^pm_i\sg_i+\sum_{i=p+1}^\nu
n_i\sg_i} \\
\theta_\sg:=\d_\sg+\sg^{r_\sg}=\displaystyle{\d_\sg+\sum_{i=1}^pr_i\sg_i}.
\end{array}
\end{equation}
Note that  $\d_\sg=\theta_\sg$ if $\sg\in L. $ Now we would like
to define a $\zn-$grading on $\Lt,$ for this, we recall
(\ref{n4}), (\ref{cormem}) and  set  for $i\in J_\ell,\;a\in
J_p,\;p+1\leq b\leq \nu$ and $r\in J_m,$
$$
\begin{array}{c}
deg(e_i)=deg(f_i)=deg(h_i):=0,\;deg(h_{i,a}^\pm):=\pm3\sg_{a},\\\;deg(h_{i,b}^\pm):=\pm\sg_b,\;\deg(v^r):=\sg^r.\\
\end{array}$$ This  defines  a $\zn-$grading
on the free  Lie algebra generated by
$$\{e_i,f_i,h_i,h_{i,a}^\pm,v^r\mid i\in J_\ell,\;a\in J_\nu,\;
r\in J_m\}.$$ As before relations (R1)-(R9) are generated by
homogeneous elements, so this grading is naturally transferred to
$\Lt$ and also to $\Lt/Z(\Lt).$ Next we note that if
$\sg=(n_1,\ldots,n_\nu)\in\zn$ and $r\in J_m$ be with
$\sg^r=\displaystyle{\sum_{i=1}^ps_i\sg_i},$ then
\begin{equation}\label{fin9}
\begin{array}{c}
\hbox{ $\gs$ is homogeneous of degree
$\displaystyle{\sum_{i=1}^p3n_i\sg_i+\sum_{i=p+1}^\nu n_i\sg_i}$ and}\\
 \;\hbox{$\v_\sg^r$ is
homogenous of degree
$\displaystyle{\sum_{i=1}^p(3n_i+s_i)\sg_i+\sum_{i=p+1}^\nu
n_i\sg_i.}$}
\end{array}
\end{equation}

As before set

$$(\Lt/Z(\Lt))_\a^\sg:=(\Lt/Z(\Lt))_\a\cap(\Lt/Z(\Lt))^\sg;\;\;\a\in R,\;\;\sg\in\zn.$$

Since $\Lt$ is a $G_2$-graded Lie algebra, Recognition Theorem
states that there are a unital commutative associative algebra $A$
and a unital Jordan algebra $\jj$ over $A$ having  a normalized
trace  satisfying the $ch_3$-identity (see Theorem \ref{recog}) such that $\Lt$ is centrally
isogenous with
$$(\gg\ot A)\op(\v\ot \jj_0)\op D_{\jj,\jj}$$ where $\v=\fc_0$ (see \S
1.1.2). Next set $B:=\jj_0$  and note that
$D_{\jj,\jj}=D_{\jj_0,\jj_0}.$
Now we have
\begin{equation}\label{weight-g}
(\Lt/Z(\Lt))_\a=\left\{\begin{array}{ll} (\gg_\a\ot A)\op(\v_\a\ot
B)&\hbox{if $\a\in R_{sh},$}\\
\gg_\a\ot A&\hbox{if $\a\in R_{lg}.$}
\end{array}\right.
\end{equation}
On the other hand,  using (\ref{lequ4}), (\ref{fin9}),
(\ref{final2}) and (\ref{cormem}),  we have for $\sg\in\zn$ and
$\a\in \rcross$ that
\begin{equation}\label{grad1}
(\Lt/Z(\Lt))_\a^\sg=\left\{
\begin{array}{ll}
(\v_{\d_\sg}^{r_\sg})_\a&\hbox{if $\a\in R_{sh}$ and $\sg\in
S\setminus L,$}\\
(\gg_{\d_\sg})_\a&\hbox{if $\a\in R^\times$ and $\sg\in L.$}
\end{array}
\right.
\end{equation}

Next let $\a\in R_{lg}$ and $\sg\in L,$ then using
(\ref{weight-g}), we get from the one dimensionality of
$(\Lt/Z(\Lt))_\a^{\sg}$  and $\gg_\a$  that there is a one
dimensional subspace $A_\a^{\d_\sg}$ of $A$ such that
$$(\Lt/Z(\Lt))_\a^{\sg}=\gg_\a\ot A_\a^{\d_\sg}.$$
Using the same argument as in simply laced types,we get that
we get $A_\a^{\d_\sg}=A_\b^{\d_\sg}$ for all
$\a,\b\in  R_{lg}$ and $\sg\in L.$ Set $$A^{\d_\sg}:=A_\a^{\d_\sg}
\hbox{ for }\sg\in L\hbox{ and any choice of }\a\in R_{lg}.$$

Now let $\sg\in L$ and  $\a$ be  a short root, take $\b$ to be a
long root such that $\a-\b\in \rcross,$ then recalling
(\ref{ubracket}), we  have
\begin{eqnarray*}
({\Lt}/{Z(\Lt)})_\a^{\sg}=(\gg_{\d_\sg})_\a=\varphi_{\d_\sg}([\gg_{\a-\b},\gg_\b])
=[\gg_{\a-\b},\varphi_{\d_\sg}(\gg_\b)]
&=&[\gg_{\a-\b},(\gg_{\d_\sg})_\b]\\
&=&[\gg_{\a-\b},(\gg_{\d_\sg})_\b]^-\\
&=&[\gg_{\a-\b},\gg_\b\ot A^{\d_\sg}]^-\\
& =&\gg_\a\ot A^{\d_\sg}.
\end{eqnarray*}
Therefore we have
\begin{equation}\label{equ-g1}
({\Lt}/{Z(\Lt)})_\a^{\sg}=\gg_\a\ot
A^{\d_\sg};\;\a\in\rcross,\;\sg\in L.
\end{equation}

Next suppose that $\a\in R_{sh}$ and  $\sg\in S\setminus L,$ then
the one dimensionality of $(\Lt/Z(\Lt))_\a^\sg$ and $\v_\a$
together with (\ref{weight-g}) and (\ref{grad1}) implies that
there is a one dimensional subspace $B_{\a}^{\theta_\sg}$ such
that
\begin{equation}\label{sumand1}
(\Lt/Z(\Lt))_\a^\sg=\v_\a\ot B_{\a}^{\theta_\sg}.
\end{equation}

If    $\sg\in S\setminus L$ and $\a,\b\in R_{sh}$ such that
$\a-\b\in \rcross$, then  we have
\begin{eqnarray*}
\v_\a\ot
B_\a^{\theta_\sg}=(\Lt/Z(\Lt))_\a^\sg=(\v_{\d_\sg}^{r_{\sg}})_\a=\psi_{\d_\sg}^{r_\sg}(\gg_{\a-\b}\cdot\v_\b)
&=& [\gg_{\a-\b},\psi_{\d_\sg}^{r_\sg}(\v_\b)]\\
&=& [\gg_{\a-\b},(\v_{\d_\sg}^{r_\sg})_\b]\\
&=&[\gg_{\a-\b},\v_\b\ot B_\b^{\theta_\sg}]^- \\
&=&\gg_{\a-\b}\cdot\v_\b\ot B_\b^{\theta_\sg}\\
&=&\v_\a\ot B_\b^{\theta_\sg},
\end{eqnarray*}
which implies that  $B_\a^{\theta_\sg}=B_\b^{\theta_\sg}$. Use
\cite[(5.11)]{AG} to conclude that
$B_\a^{\theta_\sg}=B_\b^{\theta_\sg}$ for all $\a,\b\in R_{sh}$
and define
\begin{equation}\label{sumand2}
B^{\theta_\sg}:=B_\a^{\theta_\sg}\hbox{ for }\sg\in S\setminus
L\hbox{ and any choice of } \a\in R_{sh}.\end{equation}

Now take $\sg,\tau\in L$ and  $\a,\b\in R^\times$ to be  such that
$\a+\b\in R^\times,$ then (\ref{equ-gen1}), (\ref{equ-g1}),
(\ref{ubracket}) and Remark \ref{jadid} imply that
\begin{equation*}\begin{array}{c}
0\neq
\;[(\frac{\Lt}{Z(\Lt)})^{\sg}_\a,(\frac{\Lt}{Z(\Lt)})^{\tau}_\b]^{\bar{~}}
\sub(\frac{\Lt}{Z(\Lt)})^{\sg+\tau}_{\a+\b}=\gg_{\a+\b}\ot
A^{\d_{\sg+\tau}}\andd\\
\\
\;[(\frac{\Lt}{Z(\Lt)})^{\sg}_\a,(\frac{\Lt}{Z(\Lt)})^{\tau}_\b]^{\bar{~}}=[\gg_\a\ot
A^{\d_\sg},\gg_\b\ot A^{\d_\tau}]^{\bar{~}}\sub \gg_{\a+\b}\ot
A^{\d_\sg}A^{\d_\tau}.
\end{array}
\end{equation*}

Also  for $\sg\in L, $ $\tau\in S\setminus L,$ $\a\in R_{lg}$ and
$\b\in R_{sh}$ such that $\a+\b\in R_{sh},$ (\ref{equ-gen1}),
(\ref{sumand1}), (\ref{equ-g1}), (\ref{ubracket}) and Remark
\ref{jadid} imply that
$$\begin{array}{l}
\;0\neq[(\frac{\Lt}{Z(\Lt)})_\a^{\sg},(\frac{\Lt}{Z(\Lt)})_\b^{\tau}]^-\sub
(\frac{\Lt}{Z(\Lt)})_{\a+\b}^{\sg+\tau}=\v_{\a+\b}\ot
B^{\theta_{\sg+\tau}}=\v_{\a+\b}\ot
B^{\d_\sg+\theta_{\tau}},\\
\vspace{-2mm}\\
\;[(\frac{\Lt}{Z(\Lt)})_\a^{\sg},(\frac{\Lt}{Z(\Lt)})_\b^{\tau}]^-\hspace{-1mm}=\hspace{-1mm}[\gg_\a\ot
A^{\d_\sg},\v_\b\ot
B^{\theta_\tau}]^-\hspace{-1mm}\sub\hspace{-1mm}\v_{\a+\b}\ot
A^{\d_\sg}\cdot B^{\theta_\tau}. \end{array}$$ Therefore  the one
dimensionality of the subspaces appearing on the right hand sides
implies that
\begin{equation}
\begin{array}{cl}
(i)&A^{\d_\sg }\cdot A^{\d_\tau}=A^{\d_\sg+\d_\tau}=A^{\d_{\sg+\tau}};\;\;\sg,\tau\in L,\\
\vspace{-2mm}\\
(ii)& A^{\d_\sg}\cdot
B^{\theta_\tau}=B^{\d_\sg+\theta_\tau}=B^{\theta_{\sg+\tau}};\;
\sg\in L,\tau\in S\setminus L. \label{equal-g2}
\end{array}
\end{equation}

One can see that there are short roots $\a,\b$ such that $\a+\b\in
R_{sh}$ and  $x\in \v_\a,$ $y\in\v_\b$ are such that $$D_{x,y}\neq
0\andd x*y\neq0$$ (see Subsection \ref{app1}). Now if $\sg,\gamma\in S\setminus L,$
considering (\ref{cormem}) and using (\ref{ubracket}),
(\ref{fin8}) and (\ref{newm1}), we have
$$[\v_\a\ot
t^{\d_\sg} w_{r_\sg},\v_\b\ot t^{\d_\gamma} w_{r_\gamma}]\neq 0$$
and so (\ref{grad1}) and (\ref{equ-gen1}) imply that  $0\neq
[(\Lt/Z(\Lt))_\a^\sg,(\Lt/Z(\Lt))_\b^\gamma]\sub(\Lt/Z(\Lt))_{\a+\b}^{\sg+\gamma}$
and so one  dimensionality of the subspaces  implies that
\begin{equation}\label{equ-g2}
[(\Lt/Z(\Lt))_\a^\sg,(\Lt/Z(\Lt))_\b^\gamma]=(\Lt/Z(\Lt))_{\a+\b}^{\sg+\gamma}.
\end{equation}
This together with  (\ref{grad1}), (\ref{equ-g1}), (\ref{sumand1})
and (\ref{sumand2}) implies that
\begin{equation*}\label{equ-g3}
[\v_\a\ot B^{\theta_\sg},\v_\b\ot B^{\theta_\gamma}]=\left\{
\begin{array}{ll}
\v_{\a+\b}\ot B^{\theta_{\sg+\gamma}}&\hbox{if $\sg+\gamma\in S\setminus L$}\\
\gg_{\a+\b}\ot A^{\d_{\sg+\gamma}}&\hbox{if $\sg+\gamma\in L.$}
\end{array}
\right.
\end{equation*}
Now if $0\neq a\in B^{\theta_\sg}$ and $0\neq b\in
B^{\theta_\gamma},$ then $(\Lt/Z(\Lt))_\a^\sg=\bbbc x\ot a$ and
$(\Lt/Z(\Lt))_\b^\tau=\bbbc y\ot b.$ Also using (\ref{ubracket}),
we have
$$[x\ot a,y\ot b]=D_{x,y}\ot t(a,b)+(x*y)\ot (a*b)+t(x,y)D_{a,b}.
$$
Therefore if $\sg+\gamma\in L,$ we have $0\neq t(a,b)\in
A^{\d_{\sg+\tau}}$ and $a*b=0$ which result in the fact that
$ab=a*b+t(a,b)\in A^{\d_{\sg+\gamma}}=A^{\theta_{\sg+\gamma}}.$
Also if $\sg+\gamma\in S\setminus L,$  then $0\neq a*b\in
B^{\theta_{\sg+\gamma}}$ and $t(a,b)=0$ which imply that
$ab=a*b+t(a,b)\in B^{\theta_{\sg+\gamma}}.$ Summarizing our
results, we have
\begin{equation}\label{equal3}
B^{\theta_\sg}\cdot B^{\theta_\gamma}=\left\{
\begin{array}{ll}
A^{\theta_{\sg+\gamma}}&\hbox{if $\sg+\gamma\in L,$}\\
B^{\theta_{\sg+\gamma}}&\hbox{if $\sg+\gamma\in S\setminus L.$}
\end{array}
\right.
\end{equation}

Next set $$\jj^{\sg}:=\left\{
\begin{array}{ll}
A^{\theta_{\sg}}&\hbox{if $\sg\in L,$}\\
B^{\theta_{\sg}}&\hbox{if $\sg\in S\setminus L.$}
\end{array}
\right.$$ Now as  $(\Lt/Z(\Lt))_\a=\sum_{\sg\in
\zn}(\Lt/Z(\Lt))_\a^\sg,$ $\a\in R^\times,$ one can use
(\ref{equal-g2}) and (\ref{equal3}) to get the fact that $\jj=A\op
B=\op_{\sg\in\zn}\jj^{\sg}$ is a graded Jordan algebra with one
dimensional summands satisfying
$\jj^{\sg}\cdot\jj^{\tau}=\jj^{\sg+\tau}.$ Therefore using the
same argument as in \cite[Proposition 5.58]{AG}, we get that  $A$
is graded isomorphic to $A_{[\nu]}$ and $\jj$ is graded isomorphic
to $\jj_p.$ Therefore  considering Remark \ref{jadid} and using
the same argument as in \cite[Subsection 2.6]{You}, we may
identify $A$ with $A_{[\nu]}$ such that
 \begin{equation} h_{i,a}^\pm=h_i\ot t_a^{\pm1};\;i\in
J_n,\;j\in J_\nu. \label{ja-g2} \end{equation}

Now let $r\in J_m=J_{p^3-1},$ considering (\ref{cormem}) and
(\ref{final2}), one finds that $r_{\sg^r}=r$ and $\d_{\sg^r}=0,$
so $\v_{\d_{\sg^r}}^{r_{\sg^r}}=\v_0^{r_{\sg^r}}=\v^r$ and
$\theta_{\sg^r}=\sg^r.$ Therefore using (\ref{grad1}),
(\ref{sumand1}),  (\ref{sumand2}) and (\ref{ma5}), we have
\begin{equation}\label{weight-g2}
(\v^r)_{n_s}=(\v_{0}^{r_{\sg^r}})_{n_s}=
(\v^{r_{\sg^r}}_{\d_{\sg^r}})_{n_s}=(\Lt/Z(\Lt))_{\a_{n_s}}^{\sg^r}=\v_{\a_{n_s}}\ot
B^{\theta_{\sg^r}}=\v_{\a_{n_s}}\ot B^{\sg^r}.
\end{equation}

Thus considering (\ref{my1}) for each $r\in J_{m},$ one finds
$\b^r\in B^{\sg^r} $ such that $v^r=v_{n_s}\ot \b^r.$ Now since
$\v,$ as a $\gg$-module, is generated by $v_{n_s},$ (\ref{m9})
implies that
\begin{equation} v^r_{\pm i}=v_{\pm i}\ot \b^r;\;\;n_\ell+1\leq i\leq n_s,\;r\in
J_m.\label{equal7}
\end{equation}

Now one can use Subsection \ref{app1} to  see that there is $n_\ell+1\leq i\leq n_s-1$ such that
$v_{n_s}*v_{-i}\neq 0$ and $D_{v_{n_s},v_{-i}}\neq 0.$ Using (R9),
we have   for $r,s\in J_{m}$ that

\begin{equation}\label{equal8}
[v^r,v^s_{-i}]=
m'_{-i}a'_{r,s}[(1/2)\fb^{t'_{-i}}_{\sg'_{r,s}},e_{t'_{-i}}]+m_{-i}a_{r,s}[(1/2)\fb^{t_i}_{\sg_{r,s}},v_{t_i}^{t_{r,s}}].
\end{equation}
Now if $q_1:=|\sg_{r,s}|$ and  $q_2:=|\sg'_{r,s}|,$ then there are
$j_1,\ldots,j_{q_1},i_1,\ldots,i_{q_2}\in J_\nu$ such that
$t^{\sg_{r,s}}=t_{j_1}^\pm\cdots t_{j_{q_1}}^\pm$ and
$t^{\sg'_{r,s}}=t_{i_1}^\pm\cdots t_{i_{q_2}}^\pm.$ Now
(\ref{fin3}), Convention \ref{my8}, (\ref{ja-g2}), (R4),
(\ref{equal7}), (\ref{ubracket}) and (\ref{newm2}) imply that
\begin{eqnarray}
m_{-i}a_{r,s}[(1/2)\fb^{t_i}_{\sg_{r,s}},v_{t_i}^{t_{r,s}}]&=&m_{-i}a_{r,s}[\frac{1}{2}h_{t_i}\ot
t_{j_{q_1}}^\pm,\ldots,\frac{1}{2}h_{t_i}\ot
t_{j_1}^\pm,v_{t_i}\ot \b^{t_{r,s}}]\nonumber\\
&=&m_{-i}a_{r,s}v_{t_i}\ot t_{j_1}^\pm\ldots
t_{j_{q_1}}^\pm\b^{t_{r,s}}\nonumber\\
&=&m_{-i}a_{r,s}v_{t_i}\ot t^{\sg_{r,s}}\b^{t_{r,s}}\label{equal4}\\
&=&v_{n_s}*v_{-i}\ot a_{r,s}t^{\sg_{r,s}}\b^{t_{r,s}},\nonumber
\end{eqnarray}
and
\begin{eqnarray}
m'_{-i}a'_{r,s}[(1/2)\fb^{t'_{-i}}_{\sg'_{r,s}},e_{t'_{-i}}]&=&m'_{-i}a'_{r,s}[\frac{1}{2}h_{t'_{-i}}\ot
t_{i_{q_2}}^\pm,\ldots,\frac{1}{2}h_{t'_{-i}}\ot
t_{i_1}^\pm,e_{t'_{-i}}\ot 1]\nonumber\\
&=&m'_{-i}a'_{r,s}e_{t'_{-i}}\ot t_{i_{q_2}}^\pm\ldots
t_{i_1}^\pm\nonumber\\
&=&d_{v_{n_s},v_{-i}}\ot a'_{r,s} t^{\sg'_{r,s}}.\label{equal9}
\end{eqnarray}
Also  (\ref{equal8}), (\ref{equal4}) and  (\ref{equal9}) imply
that
\begin{eqnarray*}
[v^r,v^s_{-i}]=d_{v_{n_s},v_{-i}}\ot a'_{r,s}
t^{\sg'_{r,s}}+(v_{n_s}*v_{-i})\ot
a_{r,s}t^{\sg_{r,s}}\b^{t_{r,s}},
\end{eqnarray*}
on the other hand  (\ref{equal7}) together with (\ref{ubracket})
implies that

\begin{eqnarray*}
[v^r,v^s_{-i}]&=&[v_{n_s}\ot \b^r,v_{-i}\ot \b^{s}]\\
&=&d_{v_{n_s},v_{-i}}\ot
t(\b^r,\b^s)+(v_{n_s}*v_{-i})\ot(\b^r*\b^s)+d_{\b^r,\b^s}.
\end{eqnarray*}

These together imply that
\begin{equation}\label{fequ1}
t(\b^r,\b^s)=a'_{r,s} t^{\sg'_{r,s}}\andd
\b^r*\b^s=a_{r,s}t^{\sg_{r,s}}\b^{t_{r,s}}.
\end{equation}
Now (\ref{newm1}) and (\ref{fequ1}) imply that we can identify
$\b^r$ with $w_r$ for $r\in J_m.$ Therefore we get that
$\Lt/Z(\Lt)$ is isomorphic to $\fA/Z(\fA).$ Now using the same
argument as in simply-laced types, one get that $\Lt\simeq\fA.$

\medskip\underline{\bf Step3 (Type $F_4$).} We recall  that for type
$F_4,$ there is $0\leq p\leq 3$ such  that the  corresponding pair of
the $R$-graded Lie algebra $\fA/Z(\fA)$ is $(S,L)$ where
$S:=\bbbz\sg_1\op\cdots\op\bbbz\sg_\nu$ and
$L:=2\bbbz\sg_1\op\cdots\op2\bbbz\sg_p\op\bbbz\sg_{p+1}\op\cdots\op\bbbz\sg_\nu$
and that $m=2^p-1.$  For $r\in J_{m},$ define
\begin{equation}
\begin{array}{c}
\sg^r:=\displaystyle{\sum_{i=1}^{p}s_i\sg_i} \;\;\hbox{ where }
r=\displaystyle{\sum_{i=1}^{p}2^{i-1}s_i} \hbox{ for
}s_1,s_2,s_3\in\{0,1\}.
\end{array}
\end{equation}
Let  $\sg=(n_1,\ldots,n_\nu)\in\zn,$ for $1\leq i\leq p,$ suppose
that $m_i,r_i\in\bbbz$ are such that $r_i\in\{0,1\}$ and
$n_i=2m_i+r_i,$ set
\begin{equation}\label{cormem-f4}
\begin{array}{c}
r_\sg:=\displaystyle{\sum_{i=1}^p 2^{i-1}r_i},\;
\d_\sg:=\displaystyle{\sum_{i=1}^pm_i\sg_i+\sum_{i=p+1}^\nu
n_i\sg_i} \\
\theta_\sg:=\d_\sg+\sg^{r_\sg}=\displaystyle{\d_\sg+\sum_{i=1}^pr_i\sg_i}.
\end{array}
\end{equation}
We mention that  $\d_\sg=\theta_\sg$ if $\sg\in L. $ Now we would
like to define a $\zn-$grading on $\Lt,$ for this, we note that
for $i\in J_\ell,\;a\in J_p,\;p+1\leq b\leq \nu$ and $r\in J_m,$
$$
\begin{array}{c}
deg(e_i)=deg(f_i)=deg(h_i):=0,\;deg(h_{i,a}^\pm):=\pm2\sg_a,\\\;deg(h_{i,b}^\pm):=\pm\sg_b,\;\deg(v^r):=\sg^r,\\
\end{array}$$ define  a $\zn-$grading
on the free  Lie algebra generated by
$\{e_i,f_i,h_i,h_{i,a}^\pm,v^r\mid i\in J_\ell,\;a\in J_\nu,\;
r\in J_m\}.$ As before relations (R1)-(R9) are generated by
homogeneous elements, so this grading is naturally transferred to
$\Lt$ and also on $Z(\Lt).$ Next we note that if
$\sg=(n_1,\ldots,n_\nu)\in\zn,$ and if $i\in J_p,$ $s_i\in\{0,1\}$
and $r:=\displaystyle{\sum_{i=1}^p2^{i-1}s_i},$ then
\begin{equation}\label{final3}
\begin{array}{c}
\hbox{ $\gs$ is homogeneous of degree
$\displaystyle{\sum_{i=1}^p2n_i\sg_i+\sum_{i=p+1}^\nu n_i\sg_i}$ and}\\
 \;\hbox{$\v_\sg^r$ is
homogenous of degree
$\displaystyle{\sum_{i=1}^p(2n_i+s_i)\sg_i+\sum_{i=p+1}^\nu
n_i\sg_i.}$}
\end{array}
\end{equation}
 As before set

$$(\Lt/Z(\Lt))_\a^\sg:=(\Lt/Z(\Lt))_\a\cap(\Lt/Z(\Lt))^\sg;\;\;\a\in R,\;\;\sg\in\zn.$$

Since $\Lt$ is a $F_4$-graded Lie algebra, Recognition Theorem
states that there is a unital commutative associative algebra $A$
and  an alternative  algebra $\cc$ over $A$ having  a normalized
trace $T$ satisfying the $ch_2$-identity (see Theorem \ref{recog}) such that $\Lt$ is
centrally isogenous with
$$(\gg\ot A)\op(\v\ot \cc_0)\op D_{\cc,\cc}$$ where
$\v=\mathfrak{J}_0$ (see \S 1.1.2). Next set $B:=\cc_0$ and note
that $D_{\cc,\cc}=D_{\cc_0,\cc_0}.$ Now we have
\begin{equation}\label{weight1-f4}
(\Lt/Z(\Lt))_\a=\left\{\begin{array}{ll} (\gg_\a\ot A)\op(\v_\a\ot
B)&\hbox{if $\a\in R_{sh},$}\\
\gg_\a\ot A&\hbox{if $\a\in R_{lg}.$}
\end{array}\right.
\end{equation}
Also for $\sg\in\zn$ and $\a\in \rcross$ by (\ref{final3}) and
(\ref{lequ4}),  we have
\begin{equation}\label{grad1-f4}
(\Lt/Z(\Lt))_\a^\sg=\left\{
\begin{array}{ll}
(\v_{\d_\sg}^{r_\sg})_\a&\hbox{if $\a\in R_{sh}$ and $\sg\in
S\setminus L,$}\\
(\gg_{\d_\sg})_\a&\hbox{if $\a\in R^\times$ and $\sg\in L.$}
\end{array}
\right.
\end{equation}

Using the  same argument as in type $G_2,$ one finds  one
dimensional subspaces $A^{\d_\sg},$ $\sg\in L$ and
$B^{\theta_\sg},$ $\sg\in S\setminus L,$ of $A$ and $B$
respectively satisfying
\begin{equation}
\begin{array}{cl}
(i)&A^{\d_\sg }\cdot A^{\d_\tau}=A^{\d_\sg+\d_\tau}=A^{\d_{\sg+\tau}};\;\;\sg,\tau\in L,\\
\vspace{-2mm}\\
(ii)& A^{\d_\sg}\cdot
B^{\theta_\tau}=B^{\d_\sg+\theta_\tau}=B^{\theta_{\sg+\tau}};\;
\sg\in L,\tau\in S\setminus L, \label{equal-f4}\\
(iii) &B^{\theta_\sg}\cdot B^{\theta_\tau}=\left\{
\begin{array}{ll}
A^{\theta_{\sg+\tau}}&\hbox{if $\sg,\tau\in S\setminus L, $ $\sg+\tau\in L,$}\\
B^{\theta_{\sg+\tau}}&\hbox{if $\sg,\tau\in S\setminus L, $
$\sg+\tau\in S\setminus L$}
\end{array}
\right.

\end{array}
\end{equation}
 such that
 \begin{equation}
(\Lt/Z(\Lt))_\a^\sg=\left\{
\begin{array}{ll}
\gg_\a\ot A^{\d_\sg} & \a\in R^\times,\;\; \sg\in L\\
\v_\a\ot B^{\theta_\sg} & \a\in R_{sh},\;\; \sg\in S\setminus L.
\end{array}
\right.
\end{equation}
Now for $\sg\in \zn,$ set $$\cc^{\sg}:=\left\{
\begin{array}{ll}
A^{\theta_{\sg}}&\hbox{if $\sg\in L,$}\\
B^{\theta_{\sg}}&\hbox{if $\sg\in S\setminus L.$}
\end{array}
\right.$$ Then  (\ref{equal3}) implies that
$\cc:=\op_{\sg\in\zn}\cc^{\sg}$ is a graded alternative algebra
with one dimensional summands satisfying
$\cc^{\sg}\cdot\cc^{\tau}=\cc^{\sg+\tau}.$ Therefore using the
same argument as in \cite[Proposition 5.45]{AG}, we have that $A$
is graded isomorphic to $A_{[\nu]}$ and $\cc$ is graded isomorphic
to $\aa_p.$ Therefore as before we may identify $A$ with
$A_{[\nu]}$ such that\begin{equation} h_{i,a}^\pm=h_i\ot
t_a^{\pm1};\;i\in J_n,\;j\in J_\nu. \label{ja-f4} \end{equation}
Using the same argument as in  type $G_2,$  one gets that
$\Lt/Z(\Lt)$  is isomorphic to $\fA/Z(\fA),$ and so as before, one
concludes that  $\Lt\simeq\fA.$

\section{Postponed proofs}\label{postpone}
In this section we present the proofs of Propositions \ref{all2}
and \ref{alltogether}. For this, we first need to prove the
following claim about type $F_4$.
\smallskip

\noindent
{\bf Claim.}
Suppose that $R$ is of type $F_4.$  For $\sg,\tau\in\zn,$ $ 1\leq
r,s \leq m$ with $r\neq s,$ we have

(i) $[v^s_{\sg,3},[f_3,v^r_{\tau,3}]]=0.$

(ii) $[v^s_{\sg,4},[f_4,v^r_{\tau,4}]]=0.$

(iii) $[v^s_{\sg,3},[f_4,v^r_{\tau,4}]]= -
[v^r_{\sg,3},[f_4,v^s_{\tau,4}]].$

(iv)
$[v^s_{\sg,4},[f_3,v^r_{\tau,3}]]=-[v^r_{\sg,4},[f_3,v^s_{\tau,3}]].$
\label{ch8}
\medskip

\noindent {\bf Proof of the  Claim.} $(i)$ One knows that
$\ve_1-\ve_4$ is a positive long root, take $i\in J_n$ to be such
that $\a_i=\ve_1-\ve_4.$ Now as $-\a_i+\a_{n_s}=\ve_4=\a_3$ (see
(\ref{fin7})), we have  $\bbbc[f_{i},v_\sg^s]=\bbbc v^s_{\sg,3}.$
Also  as $\ve_1-\ve_4$ is not a short root, $[f_3,v_\sg^s]=0$ and
as $\a_i$ is a long root  $[f_{i},f_3,v^r_{\tau,3}]=0.$ Therefore
using the Jacobi identity together with Propositions \ref{y3.17},
\ref{m8}($iii$) and (R9), we have
\begin{eqnarray*}
[v^s_{\sg,3},[f_3,v^r_{\tau,3}]]\in\bbbc
[[f_{i},v_\sg^s],[f_3,v^r_{\tau,3}]]&=&\bbbc([f_{i},v_\sg^s,f_3,v^r_{\tau,3}]-
[v_\sg^s,f_{i},f_3,v^r_{\tau,3}])\\
&=&\bbbc[f_i,v_\sg^s,f_3,v^r_{\tau,3}]-0\\
&=&\bbbc[f_{i},f_3,v_\sg^s,v^r_{\tau,3}]\\
&\sub&\bbbc[f_{i},f_3,\underbrace{H,\ldots,H}_{|\sg|+|\tau|},v^s,v^r_{3}]=0.
\end{eqnarray*}

$(ii)$ We know that $\a_4=\frac{1}{2}(\ve_1-\ve_2-\ve_3-\ve_4).$ Take $j,t\in J_n$ be such that
$\a_j=\frac{1}{2}(\ve_1+\ve_2+\ve_3-\ve_4)$
and $\a_t=\ve_2+\ve_3,$  which are positive roots. Now we have
$$\a_4+\a_j=\a_i=\ve_1-\ve_4 \andd\a_j-\a_4=\a_t
$$
which are long roots and so $[f_t,[e_4,\v_{\tau,-4}^r]]=0$ and $[e_4,\v_{\sg,j}^s]=0.$ Now the Jacobi identity together with  Propositions \ref{y3.17} and  \ref{m8}($iii$) implies that
\begin{eqnarray*}
[v^s_{\sg,4},[f_4,v^r_{\tau,4}]]&\in& [v^s_{\sg,4},[e_4,\v_{\tau,-4}^r]]\\&=&[[f_t,\v_{\sg,j}^s],[e_4,\v_{\tau,-4}^r]]\\
&=&[f_t,\v_{\sg,j}^s,[e_4,\v_{\tau,-4}^r]]+[\v_{\sg,j}^s,f_t,[e_4,\v_{\tau,-4}^r]]\\
&=&[f_t,[\v_{\sg,j}^s,e_4],\v_{\tau,-4}^r]+[f_t,e_4,\v^s_j,\v_{\tau,-4}^r]\\
&=&[f_t,e_4,\v_{\sg,j}^s,\v_{\tau,-4}^r]\\
&\sub&\bbbc[f_t,e_4,\underbrace{H,\ldots,H}_{|\sg|+|\tau|},\v^s_j,\v^r_{-4}]=0.
\end{eqnarray*}

(iii)
Consider Lemma \ref{rev7} and take  $i,j\in J_n$ be such that $\a_i=\ve_1-\ve_4$ and $-\a_j=\frac{1}{2}(-\ve_1-\ve_2-\ve_3+\ve_4)=-\a_i+\a_4.$ Next we note   that $-\a_4+\a_3$ is not a root and that $\a_i+\a_{n_s}=\ve_4=\a_3.$ Therefore
$[f_4, v^s_{\sg,3}]=0$ and  $f_i\cdot v_{n_{s}}=z v_{3} $ for some $z\in\bbbc$
which in turn implies that  $z[f_i, v_\sg^s]= v^s_{\sg,3}. $  We next note that as $\a_{n_s}\neq \a_j,$ there are $h_1,h_2\in\hh$ such that   $\a_{n_s}(h_1)=1=\a_j(h_2)$ and
$\a_{n_s}(h_2)=0=\a_j(h_1).$ We take $c_i:=(h_1)^{i,\sg},$
$c'_j:=(h_2)^{j,\tau}$ for $1\leq i\leq |\sg|,1\leq j\leq |\tau|$ (see Convention \ref{my8}). Also as $\a_{n_s}+\a_4$
is not a root, Propositions \ref{y3.17}, \ref{m8}($iii$)  together with (R9) imply that $[v_\sg^s,v^r_{\tau,4}]=0.$
So the Jacobi identity and  Proposition \ref{m8}($iii$) imply that
\begin{eqnarray*}
[v^s_{\sg,3},[f_4,v^r_{\tau,4}]]=[f_4,v^s_{\sg,3},v^r_{\tau,4}]&=&z[f_4,[f_i, v_\sg^s],v^r_{\tau,4}]\\
&=&z[f_4,f_i, [v_\sg^s,v^r_{\tau,4}]]-z[f_4,v_\sg^s,f_i, v^r_{\tau,4}]\\
&=&-z[f_4,v_\sg^s,f_i, v^r_{\tau,4}]\\
&=&-z[f_4,c_1,\ldots,c_{|\sg|},c'_1,\ldots,c'_{|\tau|},v^s,f_i,
v^r_{4}].
\end{eqnarray*}
Using (R4) together with the same argument as above, we conclude
that
$$[v^r_{\sg,3},[f_4,v^s_{\tau,4}]]=-z[f_4,c_1,\ldots,c_{|\sg|},c'_1,\ldots,c'_{|\tau|},v^r,f_\gamma, v^s_2].$$
But $F_4,$ $a_{r,s}=-a_{s,r}$ (see (\ref{newm1})), so
using (R9), we are done.

$(iv)$ Use the same  argument as in part $(i).$\qed

\vspace{5mm}
We are now ready to present proofs of the postponed propositions.

\vspace{3mm}
\noindent{\bf Proof of Proposition
\ref{all2}.}

($i$) We prove the first equality, as the argument for
the second equality is similar. We use induction on
$|\tau|.$  If $|\tau|=0,$ then $e_{i,\tau}=e_i,$ and
$[e_i,e_{i,\sg}]=\varphi_{\sg}([e_i,e_i])=0$. Now let the result
holds for all $\tau$ with $|\tau|=s.$ Suppose that $
\sg'=(m_1,\ldots,m_\nu)$ be such that $|\sg'|=s+1,$ then there is
$a\in J_\nu,$ such that $m_a\neq 0.$ Set
$\tau:=\sg'-\sg_a^{sgn(m_a)}$ and $\gamma:=\sg+\sg_a^{sgn(m_a)}$
(see (\ref{n4})), then $|\tau|=|\sg'|-1=s.$
 Now  Proposition
\ref{y3.17},  the Jacobi identity and probably (R4), (R6) together
with the induction hypothesis give
\begin{eqnarray*}
[e_{i,\sg'},e_{i,\sg}]&=&[[\frac{1}{2}h_{i,a}^{sgn(m_a)},e_{i,\tau}],e_{i,\sg}]\\
&=&[\frac{1}{2}h_{i,a}^{sgn(m_a)},e_{i,\tau},e_{i,\sg}]-[e_{i,\tau},\frac{1}{2}h_{i,a}^{sgn(m_a)},e_{i,\sg}]\\
&=&[\frac{1}{2}h_{i,a}^{sgn(m_a)},e_{i,\tau},e_{i,\sg}]-[e_{i,\tau},e_{i,\gamma}]
=0-0=0.
\end{eqnarray*}

\medskip
($ii$)   Since
$i\neq j,$ there are $x,y\in\hh$ such that
$$
\a_j(x)=1,\; \a_i(x)=0,\; \a_j(y)=0,\; \a_i(y)=1.
$$
Now using  (\ref{m1}), the Jacobi identity and Propositions
\ref{y3.17}, \ref{m8}($iii$),($ii$),
 we have
\begin{eqnarray*}
[e_{i,\sg},h_{j,\tau}]&=&-[e_{i,\sg},f_j,e_{j,\sg}]\\
&=&[f_j,e_{j,\sg},e_{i,\sg} ]\\
&=&[f_j,[x^{1,\tau},\ldots,x^{|\tau|,\tau},e_j],[y^{1,\sg},\ldots,y^{|\sg|,\sg},e_i] ]\\
&=&[f_j,x^{1,\sg},\ldots,x^{|\sg|,\sg},y^{1,\tau},\ldots,y^{|\tau|,\tau},e_j,e_i ]\\
&=&[f_j,y^{1,\sg},\ldots,y^{|\sg|,\sg},y^{1,\tau},\ldots,y^{|\tau|,\tau},e_{j},e_i ]\\
&=&[y^{1,\sg},\ldots,y^{|\sg|,\sg},y^{1,\tau},\ldots,y^{|\tau|,\tau},f_j,e_{j},e_i ]\\
&=&-\a_i(h_j)[y^{1,\sg},\ldots,y^{|\sg|,\sg},y^{1,\tau},\ldots,y^{|\tau|,\tau},e_{i}]=-\a_{i}(h_j)e_{i,\tau+\sg}.
\end{eqnarray*}
Note that to get the last equality, we may  use (R6), (R4)
(Canceling relations). The second assertion  is similarly proved.
\medskip

($iii$) We first
note that by (\ref{m1}), the Jacobi identity  and part $(i)$ of the
proposition, we have
\begin{equation}\label{y22}
[e_{i,\tau},h_{i,\sg}]=[e_{i,\tau},[e_{i,\sg},f_i]]=
[f_i,\stackrel{0} {\overbrace{e_{i,\sg},e_{i,\tau}}}]-[e_{i,\sg},
\varphi_\tau [f_i,e_i]] =[e_{i,\sg},h_{i,\tau}].
\end{equation}
Therefore using Proposition \ref{y3.17}, we have
\begin{eqnarray*}
[h_{i,\tau},e_{i,a}^\pm]=[h_{i,a}^\pm,e_{i,\tau}]=2e_{i,\tau^{\pm_a}}.
\end{eqnarray*}
The second claim   is similarly proved.

\medskip

($iv$) For $t\in
J_{|\sg|},$  using (R5), one finds   that
$(\a_i(h_j)/2)[h_{i,a}^\pm,f_i]=[h_{j,a}^\pm,f_i]$ and
$(\a_i(h_j)/2)[h_{i,a}^\pm,e_i]=[h_{j,a}^\pm,e_i].$ Also by
Proposition \ref{y3.17}, we get that
$[(\frac{-1}{2}h_i)^{1,\sg},\ldots,(\frac{-1}{2}h_i)^{|\sg|,\sg},f_i]=f_{i,\sg}.$
Now these together with  the Jacobi identity and (R4) imply that
{\small
\begin{eqnarray*} [h_{j,a}^\pm,e_i,f_{i,\sg}] &=&
[e_i,h_{j,a}^\pm,f_{i,\sg}]+[[h_{j,a}^\pm,e_i],f_{i,\sg}]\\
&=&
[e_i,(\frac{-1}{2}h_i)^{1,\sg},\ldots,(\frac{-1}{2}h_i)^{|\sg|,\sg},h_{j,a}^\pm,f_i]+[[h_{j,a}^\pm,e_i],f_{i,\sg}]\\
&=&
\frac{\a_i(h_j)}{2}([e_i,(\frac{-1}{2}h_i)^{1,\sg},\ldots,(\frac{-1}{2}h_i)^{|\sg|,\sg},h_{i,a}^\pm,f_i]+[[h_{i,a}^\pm,e_i],f_{i,\sg}])\\
&=&
\frac{\a_i(h_j)}{2}([e_i,h_{i,a}^\pm,f_{i,\sg}]+[[h_{i,a}^\pm,e_i],f_{i,\sg}])=
\frac{\a_i(h_j)}{2}[h_{i,a}^\pm,e_i,f_{i,\sg}].
\end{eqnarray*}}
Now we are done as by (\ref{m1}), $[e_i,f_{i,\sg}]=h_{i,\sg}.$

\medskip

($v$)
 Using Proposition \ref{y3.17} together with the Jacobi
identity and \ref{m1}, we have
\begin{eqnarray*}
[e_{i,a}^\pm,f_{i,\tau}]=[[\frac{1}{2}h_{i,a}^\pm,e_i],f_{i,\tau}]&=&[\frac{1}{2}h_{i,a}^\pm,e_i,f_{i,\tau}]-[e_i,\frac{1}{2}h_{i,a}^\pm,f_{i,\tau}]\\
&=&[\frac{1}{2}h_{i,a}^\pm,h_{i,\tau}]+[e_i,f_{i,\tau^{\pm_a}}]\\
&=&\frac{1}{2}[h_{i,a}^\pm,h_{i,\tau}]+h_{i,\tau^{\pm_a}}.
\end{eqnarray*}

($vi$)
 We first note that as $\a_{n_s}\neq \a_i,$  there
are  $x,y\in\hh$ such that $\a_{n_s}(y)=0,$ $\a_{i}(y)=1,$
$\a_{n_s}(x)=1,$ $\a_{i}(x)=0$ and so using Proposition
\ref{m8}($i$), we have  $[v^r,y^{i,\sg}]=0,$ $i\in J_{|\sg|},$
which together with Proposition \ref{y3.17} implies that
$[v^r,f_{i,\sg}]=[y^{1,\sg},\ldots,y^{|\sg|,\sg},v^r,f_i].$ Now
using this, (\ref{m1}), the Jacobi identity, (R2), Propositions
\ref{m8}($iii$) and \ref{y3.17}, we get that
\begin{eqnarray*}
[v^r,h_{i,\sg}]=[v^r,e_i,f_{i,\sg}]=[e_i,v^r,f_{i,\sg}]
&=&[x^{1,\sg},\ldots,x^{|\sg|,\sg},e_i,v^r,f_i]\\
&=&-\a_{n_s}(h_i)[x^{1,\sg},\ldots,x^{|\sg|,\sg},v^r]\\&=&-\a_{n_s}(h_i)v^r_\sg.
\end{eqnarray*}
This completes the proof of the  first assertion. Now one knows
that
$[v^r,h_{i,\sg}]=-\a_{n_s}(h_i)[x^{1,\sg},\ldots,x^{|\sg|,\sg},v^r]$
and
$[v^r,h_{i,\tau}]=-\a_{n_s}(h_i)[x^{1,\tau},\ldots,x^{|\tau|,\tau},v^r].$
Next we note that  using (\ref{m1}), Propositions \ref{y3.17},
\ref{m8}($i$), the Jacobi identity and (R4), we get that
$[h_{i,\tau},x^{j,\sg}]=0,$ $j\in J_{|\sg|},$ so
\begin{eqnarray*}
[h_{i,\tau},h_{i,\sg},v^r]&=&\a_{n_s}(h_i)[x^{1,\sg},\ldots,x^{|\sg|,\sg},h_{i,\tau},v^r]\\
&=&(\a_{n_s}(h_i))^2[x^{1,\sg},\ldots,x^{|\sg|,\sg},x^{1,\tau},\ldots,x^{|\tau|,\tau},v^r].
\end{eqnarray*}
Similarly
$
[h_{i,\tau},h_{i,\sg},v^r]=(\a_{n_s}(h_i))^2[x^{1,\sg},\ldots,x^{|\sg|,\sg},x^{1,\tau},\ldots,x^{|\tau|,\tau},v^r].
$
 Now  using (R4), we are done. This completes the proof.\qed

\bigskip

\noindent{\bf Proof of  Proposition \ref{alltogether}}

($i$) It is
immediate using Proposition \ref{y3.17} together with (R4), (R6).

\medskip

($ii$) We
first suppose that $t\neq n_s,$ then there exists $x\in\hh$ such
that $\a_{n_s}(x)=0$ and $\a_t(x)=1.$ Recall Convention \ref{my8},
and set $c_j:=x^{j,\sg}$ for $j\in J_{|\sg|}.$
 Then using  Propositions \ref{y3.17} and
\ref{m8}($i$), we have
\begin{eqnarray*}
[v^r,(\gs)_{\pm\a_t}]=[v^r,c_1,\ldots,c_{|\sg|},\gg_{\pm\a_t}]&=&[c_1,\ldots,c_{|\sg|},v^r,\gg_{\pm\a_t}]\\
&\sub&[c_1,\ldots,c_{|\sg|},(\v^r)_{\pm\a_t+\a_{n_s}}]\\
&=&(\v^r_\sg)_{\pm\a_t+\a_{n_s}}.
\end{eqnarray*}

Next  let $t=n_s,$ then as $\a_t$ is not a simple root, there is
$1\leq j\leq \ell$ and $1\leq i\leq n$ such that $\a_t-\a_j=\a_i.$
Now $(\gg_\sg)_{\a_t}=[e_j,(\gs)_{\a_i}]$ and so by (R2), the
Jacobi identity and the first part of the proof,  we have
\begin{eqnarray*}
[v^r,(\gs)_{\a_t}]=[v^r,e_j,(\gs)_{\a_i}]=[e_j,v^r,(\gs)_{\a_i}]&\sub&[e_j,(\v^r_\sg)_{\a_i+\a_{n_s}}]\\
&\sub&(\v^r_\sg)_{\a_t+\a_{n_s}}.
\end{eqnarray*}
Also as $\a_t=\a_{n_s}$ is a short  root, Lemma \ref{rev7} implies
that  there are $ i,j\in J_n$ such that $\a_{t}+\a_j=\a_i,$ so
$(\gg_\sg)_{-\a_t}=[e_j,(\gs)_{-\a_i}].$ Now as before   we have
\begin{eqnarray*}
[v^r,(\gs)_{-\a_t}]=[v^r,e_j,(\gs)_{-\a_i}]=[e_j,v^r,(\gs)_{-\a_i}]&\sub&[e_j,(\v^r_\sg)_{-\a_i+\a_{n_s}}]\\
&\sub&(\v^r_\sg)_{-\a_t+\a_{n_s}}.
\end{eqnarray*}
Next using (R2), the Jacobi identity and the first part of the
proof, we have
\begin{eqnarray*}
[v^r,(\gs)_0]=\sum_{i=1}^\ell[v^r,e_i,(\gs)_{-\a_i}]&=&\sum_{i=1}^\ell[e_i,v^r,(\gs)_{-\a_i}]\\
&\sub&\sum_{i=1}^\ell[e_i,(\v^r_\sg)_{-\a_i+\a_{n_s}}]\sub(\v^r_\sg)_{\a_{n_s}}.
\end{eqnarray*}

\medskip

($iii$)
Using Lemma \ref{rev1}($ii$), One knows that $(\v^r_\sg)_0$ is
spanned by $[f_{i}, v_{\sg,i}^r]$ for $n_\ell+1\leq i\leq\ell.$ So
it is enough to show that there is a basis $\{x_j\in\hh\mid j\in
J_\ell\}$ of $\hh$ such that $[(x_j)_a^\pm,f_{i}, v_{\sg,i}^r]=0$
for all $j\in J_\ell$ and $a\in J_\nu.$ Fix $n_\ell+1\leq
i\leq\ell.$ By Lemma  \ref{rev7}($i$), $\a_{i}=\b+\gamma$ where
$\b$ is a short root and $\gamma$ is a long root. Take $t,t'\in
J_n$ be such that $\b=\pm\a_t$ and $\gamma=\pm\a_{t'}$ and set

$$
e:=\left\{\begin{array}{ll} f_{t'}&\hbox{if $\gamma=\a_{t'}$}\\
e_{t'}&\hbox{if $\gamma=-\a_{t'},$}
\end{array}\right.\;
f:=\left\{\begin{array}{ll} f_{t}&\hbox{if $\b=\a_{t}$}\\
e_{t}&\hbox{if $\b=-\a_{t},$}
\end{array}\right.\;v:=\left\{\begin{array}{ll} v^r_{t}&\hbox{if $\b=\a_{t}$}\\
v^r_{-t}&\hbox{if $\b=-\a_{t},$}
\end{array}\right.
$$then $e\in\gg_{-\gamma},f\in\gg_{-\b},$ $[e,f]\in\gg_{-\a_i}$ and $v\in(\v^r)_\b$ and
so we have
\begin{equation}\label{rev4}
[e,v^r_{\sg,i}]\in\bbbc \psi^r_\sg(v)\andd [f,v^r_{\sg,i}]=0.
\end{equation}

Consider  a basis $\{x_i\mid 1\leq i\leq \ell-1\}$ for
$\hbox{ker}(\a_{i}).$ Also note that since $\a_{t'}\neq \a_{i},$
there is $x_\ell\in\hh$ such that $ \a_{t}(x_\ell)=0$ and
$\a_i(x_\ell)\neq0.$ Therefore  $\{x_i\mid 1\leq i\leq \ell\}$ is
a basis for $\hh.$ Now we note that using Proposition \ref{y3.17},
$v^r_{\sg,i}\in[\underbrace{H,\ldots,H}_{|\sg|-times},v^r_{i}]$
and $\psi^r_\sg(v)\in[\underbrace{H,\ldots,H}_{|\sg|{-times}},v].$
Now using  the Jacobi identity together with Proposition
\ref{m8}($i$), (R4) and (\ref{rev4}),
  we have for $1\leq j'\leq \ell-1,$
\begin{eqnarray*}
[(x_{j'})_a^\pm ,f_{i},v^r_{\sg,i}] \in\bbbc
[f_{i},(x_{j'})_a^\pm,\underbrace{H,\ldots,H}_{|\sg|{-times}},v^r_{i}]\sub
\bbbc[f_{i},\underbrace{H,\ldots,H}_{|\sg|{-times}},(x_{j'})_a^\pm,v^r_{i}]=0,
\end{eqnarray*}
and
\begin{eqnarray*}
[(x_\ell)_a^\pm,f_{i},v^r_{\sg,i}]\in
\bbbc[(x_\ell)_a^\pm,[e,f],v^r_{\sg,i}] &\sub&\bbbc
[(x_\ell)_a^\pm, [f,\psi^r_\sg(v)]]
\\
&\sub& \bbbc([f,(x_\ell)_a^\pm,\underbrace{H,\ldots,H}_{|\sg|{-times}},v]\\
&\sub&
\bbbc([f,\underbrace{H,\ldots,H}_{|\sg|{-times}},(x_\ell)_a^\pm,v]=0.
\end{eqnarray*}
 These all together complete the proof of this part.
 \medskip

($iv$) It
follows from part $(iii)$ of the proposition  together with
Proposition \ref{y3.17}, (R4) and (R7).
\medskip

($v$) Using
Remark \ref{rev7}$(i)$, one finds    $t,t'\in J_n$ such that
$\a_{t'}$ is a short root, $\a_{t}$ is a  long root and
$\a_{n_s}=\pm\a_{t'}+\a_{t}.$ We note that $v^r_\sg$ is a weight
vector of weight $\a_{n_s}$ and so
$v^r_\sg\in\bbbc[e_t,v^r_{\sg,\pm t'}].$ Now using (R2), we have
$$[v^s,v^r_\sg]\in\bbbc[v^s,[e_t,v^r_{\sg,\pm t'}]]=[[e_t,v^s],\v^r_{\sg,\pm t'}]+[e_t,v^s,v^r_{\sg,\pm t'}]=[e_t,v^s,v^r_{\sg,\pm t'}].$$

Since $t'\neq n_s,$ there is $h\in\hh$ such that $\a_{t'}(h)=1$
and $\a_{n_s}(h)=0.$ For
 $ i\in J_{|\sg|}$ set
$c_i:=h^{i,\sg}$ (Convention \ref{my8}). Now using Proposition
\ref{y3.17}, we have $\bbbc v^r_{\sg,\pm
t'}=\bbbc[c_1,\ldots,c_{|\sg|},v^r_{\pm t'}]$ and so by
Proposition \ref{m8}($i$),
\begin{eqnarray*}
[v^s,v^r_{\sg,\pm t'}]\in \bbbc[v^s,c_1,\ldots,c_{|\sg|},v^r_{\pm
t'}]=\bbbc[c_1,\ldots,c_{|\sg|},v^s,v^r_{\pm t'}].
\end{eqnarray*}
Now as $t'\neq n_s$ and $2\a_{n_s}$ is not a root, we are done
using (R9) together with Proposition \ref{y3.17}.

($vi$)
Suppose that $\a$ is a short root not equal to $\pm\a_{n_s},$ then
there is $n_\ell+1\leq t\leq n_s-1$ such that $\a=\pm\a_{t}.$
Since $\a\neq \a_{n_s},$ there is $h\in \hh$ such that $\a(h)=1$
and $\a_{n_s}(h)=0.$ Now using Proposition \ref{y3.17}, one gets
that
$(\v^r_\sg)_{\pm\a_t}=\bbbc[h^{1,\sg},\ldots,h^{|\sg|,\sg},v_{\pm
t}^r]$ and so Proposition \ref{m8}($iii$) implies that
\begin{eqnarray*} [v^s,(\v^r_\sg)_{\pm\a_t}]=\bbbc[v^s,h^{1,\sg},\ldots,h^{|\sg|,\sg},v^r_{\pm
t}]=\bbbc[h^{1,\sg},\ldots,h^{|\sg|,\sg},v^s,v^r_{\pm t}].
\end{eqnarray*}
Now we are done using (R9) together with Proposition \ref{y3.17}.
Next suppose $\a=0.$ Using Lemma \ref{rev1}($ii$), $(\v^r_\sg)_0$
is spanned by   $[e_{i}, v_{\sg,-i}^r],$  $n_\ell+1\leq
i\leq\ell.$ Fix $n_\ell+1\leq i\leq\ell,$ then we have
$$[v^s,[e_i,v_{\sg,-i}^r]]=[e_i,[v^s,v_{\sg,-i}^r]].$$
This together with  the previous step completes the proof.

($vii$)
 We first consider type $F_4.$ Take $i,j,{t'},p,q\in J_n $ to
be such that
$$\begin{array}{c}
\a_i=\frac{1}{2}(\ve_1-\ve_2-\ve_3-\ve_4),\;\a_j=\ve_2+\ve_3, \;
\a_{t'}=\frac{1}{2}(\ve_1-\ve_2-\ve_3+\ve_4),\\\a_p=\frac{1}{2}(\ve_1+\ve_2+\ve_3-\ve_4),\;
\a_q=\frac{1}{2}(\ve_1+\ve_2+\ve_3+\ve_4).
\end{array}$$
We recall $b,b'$ as in (\ref{ch9}), then one observes that there
are nonzero scalars $a_1,a_2,a_3,a_4,y,z\in\bbbc$ satisfying
$\frac{-za'}{a_2}+\frac{ya}{a_3a_1}\neq0$ and
\begin{equation}\label{fin4}\begin{array}{c}
\;[[f_j,f_i],v_{p}]=a_1[f_4,v_{4}],
\;[v_{n_s},f_{n_s}]=a_2[v_3,f_3],\;v_{n_s}=a_3[e_{t'},v_{p}],\\
\;[f_{t'},f_j,f_i,b'[f_3,v_3]+b [f_4,v_4]]=y v_{-{n_s}},\;
[f_{n_s},b'[f_3,v_3]+b [f_4,v_4]]=zv_{-{n_s}}.
\end{array}
\end{equation}
Now we recall that  $\v^s$ and $\v_\sg^r$ are $\gg-$modules whose
weights are short roots  and note that $\a_p-\a_i,\a_{t'}+\a_q$
are not short roots and $\a_i+\a_q=\ve_1=\a_{n_s}$. So the Jacobi
identity together with part $(vi)$ of the proposition implies that
\begin{eqnarray}
X:=y[e_{t'},v^s_p,v^r_{\sg,-n_s}]&\in&\bbbc[e_{t'},v^s_p,f_i,v^r_{\sg,-q}]\nonumber\\
&=&
\bbbc([e_{t'},\underbrace{[v^s_p,f_i]}_0,v^r_{\sg,-q}]+[e_{t'},f_i,v^s_p,v^r_{\sg,-q}])\nonumber\\
&=&\bbbc[e_{t'},f_i,[f_{t'},v^s],v^r_{\sg,-q}]\label{fin5}\\
&\sub&\bbbc[e_{t'},f_i,f_{t'},v^s,v^r_{\sg,-q}]+\bbbc[e_{t'},f_i,v^s,\underbrace{f_{t'},v^r_{\sg,-q}}_0]\nonumber\\
&\sub&[e_{t'},f_i,f_{t'},\sum_{\tau\in\zn}(\gg_\tau)_{\a_i}+\sum_{1\leq
p\leq m}\sum_{\tau\in\zn}(\v^p_\tau)_{\a_i}]
\nonumber\\
&\sub&\sum_{\tau\in\zn}(\gg_\tau)_{0}+\sum_{1\leq k\leq
m}\sum_{\tau\in\zn}(\v^k_\tau)_{0}.\nonumber
\end{eqnarray}
Next we note that $\a_j+\a_i=\a_p=\a_{n_s}-\a_{t'},$ so the Jacobi
identity, Proposition \ref{y3.1} and part $(vi)$ of the
proposition imply that
\begin{eqnarray}
[[f_j,f_i],v^s_p,(\v^r_\sg)_0]&\sub&\bbbc[f_p,v^s_p,(\v^r_\sg)_0]\nonumber\\
&\sub&\bbbc[f_p,[f_{t'},v^s],(\v^r_\sg)_0]\nonumber\\
&=&\bbbc([f_p,f_{t'},v^s,(\v^r_\sg)_0]-[f_p,v^s,f_{t'},(\v^r_\sg)_0])\label{fin6}\\
&=&\bbbc([f_p,f_{t'},v^s,(\v^r_\sg)_0]-[f_p,v^s,(\v^r_\sg)_{-\a_{t'}}])\nonumber\\
&\sub&\sum_{\tau\in\zn}(\gg_\tau)_{0}+\sum_{1\leq k\leq
m}\sum_{\tau\in\zn}(\v^k_\tau)_{0}.\nonumber
\end{eqnarray}
Now set $v_{0,1}:=[f_3,v^r_{\sg,3}],v_{0,2}:=[f_4,v^r_{\sg,4}],$
we mention that as $\psi^s$ and $\psi^r_\sg$ are $\gg-$module
homomorphisms, (\ref{fin4}) implies that
$$[[f_j,f_i],v^s_{p}]=a_1[f_4,v^s_{4}],
\;v^s_{n_s}=a_3[e_{t'},v^s_{p}],\;[f_{t'},f_j,f_i,b'v_{0,1}+b
v_{0,2}]=y v^r_{-\a_{n_s}}.$$ Next we note that
$b'v_{0,1}+bv_{0,2}\in(\v^r_\sg)_0,$ so  Therefore
$[f_j,f_i,b'v_{0,1}+bv_{0,2}]\in(\v^r_\sg)_{-\a_i-\a_j}$ and so
$[h_{t'},f_j,f_i,b'v_{0,1}+bv_{0,2}]=[f_j,f_i,b'v_{0,1}+bv_{0,2}].$
Also as $\a_j$ and $\a_i+\a_j-\a_{t'}$ are not short roots,
$[e_{t'},f_j,f_i,b'v_{0,1}+bv_{0,2}]=0$ and
$[f_j,b'v_{0,1}+bv_{0,2}]=0.$ These all together with (\ref{fin5})
and (\ref{fin6}) imply that

 {\small\begin{eqnarray*}
[v^s,v^r_{\sg,-n_s}]\hspace{-3mm}&=&\hspace{-3mm}\frac{a_3}{y}[[e_{t'},v^s_p],f_{t'},f_j,f_i,b'v_{0,1}+bv_{0,2}]\\
\hspace{-3mm}&=&\hspace{-3mm}\frac{a_3}{y}([e_{t'},v^s_p,f_{t'},f_j,f_i,b'v_{0,1}+bv_{0,2}]\hspace{-1mm}-\hspace{-1mm}[v^s_p,e_{t'},f_{t'},f_j,f_i,b'v_{0,1}+bv_{0,2}])\\
\hspace{-3mm}&=&\hspace{-3mm}\frac{a_3}{y}(X-[v^s_p,h_{t'},f_j,f_i,b'v_{0,1}+bv_{0,2}]
\hspace{-1mm}
-\hspace{-1mm}[v^s_p,f_{t'},e_{t'},f_j,f_i,b'v_{0,1}+bv_{0,2}])\\
\hspace{-2mm}&=&\hspace{-2mm}\frac{a_3}{y}(X-[v^s_p,f_j,f_i,b'v_{0,1}+bv_{0,2}])\\
\hspace{-2mm}&=&\hspace{-2mm}\frac{a_3}{y}(X-[v^s_p,[f_j,f_i],b'v_{0,1}+bv_{0,2}]-[v^s_p,f_i,f_j,b'v_{0,1}+bv_{0,2}])\\
\hspace{-2mm}&=&\hspace{-2mm}\frac{a_3}{y}(X-[[v^s_p,[f_j,f_i]],b'v_{0,1}+bv_{0,2}]-[[f_j,f_i],v^s_p,b'v_{0,1}+bv_{0,2}])\\
&\in&\frac{a_3}{y}(X+a_1[[f_4,v^s_{4}],b'v_{0,1}+bv_{0,2}]-[[f_j,f_i],v^s_p,(\v_\sg^r)_0]\\
&\equiv&\frac{a_1a_3}{y}[[f_4,v^s_{2}],b'v_{0,1}+bv_{0,2}]\;\;\;\hbox{(mod}\sum_{\tau\in\zn}(\gg_\tau)_{0}
+\sum_{1\leq k\leq m}\sum_{\tau\in\zn}(\v^k_\tau)_{0}).
\end{eqnarray*}}
On the other hand using part $(vi)$ of the Proposition, we have
$[f_{n_s},v^s,b'v_{0,1}+bv_{0,2}]\in\sum_{\tau\in\zn}(\gg_\tau)_{0}+\sum_{1\leq
k\leq m}\sum_{\tau\in\zn}(\v^k_\tau)_{0},$ so the Jacobi identity
implies that
\begin{eqnarray*}
[v^s,v^r_{\sg,-n_s}]&=&\frac{1}{z}[v^s,f_{n_s},b'v_{0,1}+bv_{0,2}]\\
&=&\frac{1}{z}([[v^s,f_{n_s}],b'v_{0,1}+bv_{0,2}]+[f_{n_s},v^s,b'v_{0,1}+bv_{0,2}])\\
&=&\frac{a_2}{z}([[v_3^s,f_{3}],b'v_{0,1}+bv_{0,2}]+[f_{n_s},v^s,b'v_{0,1}+bv_{0,2}])\\
&\equiv&\frac{a_2}{z}[[v_3^s,f_{3}],b'v_{0,1}+bv_{0,2}]\;\;\;\hbox{\small(mod
$\displaystyle{\sum_{\tau\in\zn}(\gg_\tau)_{0}+\sum_{1\leq k\leq
m}\sum_{\tau\in\zn}(\v^k_\tau)_{0}}).$}
\end{eqnarray*}
Now we have the following congruences   modulo
{\small$\displaystyle{\sum_{\tau\in\zn}(\gg_\tau)_{0}+\sum_{1\leq
k\leq m}\sum_{\tau\in\zn}(\v^k_\tau)_{0}}:$}
{\small\begin{eqnarray*}
(\frac{-za'}{a_2}+\frac{ya}{a_3a_1})[v^s,v^r_{\sg,-\a_{n_s}}]&\equiv&
\frac{-za'}{a_2}[v^s,v^r_{\sg,-\a_{n_s}}]
+\frac{ya}{a_3a_1}[v^s,v^r_{\sg,-\a_{n_s}}]\\
&\equiv&[a'[f_3,v^s_1]+a[f_4,v^s_2],b'[f_3,v_{\sg,3}^r]+b[f_4,v^r_{\sg,4}]]\in\dd.
\end{eqnarray*}}

For types other than $F_4$, we first note that
$\bbbc[f_\ell,v^s_{\ell}]=\bbbc[f_{n_s},v^s].$ Now we have
\begin{eqnarray*}
[v^s,v^r_{\sg,-{n_s}}]\in[v^s,[f_{n_s},(\v^r_\sg)_0]]&=&\bbbc[v^s,[f_{n_s},[f_\ell,v^r_{\ell,\sg}]]]\\
&\sub&\bbbc[[v^s,f_{n_s}],[f_\ell,v^r_{\ell,\sg}]]+\bbbc[f,v^s,[f_\ell,v^r_{\ell,\sg}]]\\
&\sub&\bbbc[[f_\ell,v^s_{\ell}],[f_\ell,v^r_{\ell,\sg}]]+\bbbc[f_{n_s},v^s,[f_\ell,v^r_{\ell,\sg}]].
\end{eqnarray*}
Now we are done using part $(vi)$ of the proposition.

\medskip

($viii$),($ix$) Use the same argument as in
[You, Lemma 2.5].

\medskip

($x$) We
fix  $p,q\in J_m$ and $\gamma,\d\in \bbbz$ and show that
$$[e_i,D^{p,q}_{\gamma,\d}]=[f_i,D^{p,q}_{\gamma,\d}]=0;\; 1\leq i\leq\ell.$$
We know that $\v^p_\gamma,$ $\v^q_\d$ are $\gg-$modules, so for
$n_\ell+1\leq j\leq \ell,$
$[f_j,v^p_{\gamma,j}]\in(\v^p_{\gamma})_0$ and
$[f_j,v^q_{\d,j}]\in(\v^q_{\d})_0,$ therefore
$[e_i,f_j,v^p_{\gamma,j}]=0$ and $[e_i,f_j,v^q_{\d,j}]=0$ for
$1\leq i\leq n_\ell$ which together with the Jacobi identity
implies that
$$[e_i,D^{p,q}_{\gamma,\d}]=[f_i,D^{p,q}_{\gamma,\d}]=0;\; 1\leq i\leq n_\ell.$$
Therefore it remains to prove that
$$[e_i,D^{p,q}_{\gamma,\d}]=[f_i,D^{p,q}_{\gamma,\d}]=0;\; n_\ell+1\leq i\leq \ell.$$
We first consider type $F_4.$ In this case the only simple short
roots appearing in our fixed basis are $\a_3,\a_4$. Therefore it
is enough to show that $[e_i,D^{p,q}_{\gamma,\d}]=0$ for $i=3,4.$
Using the Jacobi identity, the Claim  and
appealing (\ref{ch9}), we have
\begin{eqnarray*}
[e_3,D^{p,q}_{\gamma,\d}]&=&[e_3,[a'[f_3,v^p_{\gamma,3}]+a [f_4,v_{\gamma,4}^p],b'[f_3,v^q_{\d,3}]+b[f_4,v_{\d,4}^q]]]\\
&=&(a_3''a'+b_3''a)[v^p_{\gamma,3},b'[f_3,v^q_{\d,3}]+b[f_4,v_{\d,4}^q]]\\
&-&(a_3''b'+b_3''b)[v^q_{\gamma,3},a'[f_3,v^p_{\d,3}]+a[f_4,v_{\d,4}^p]]\\
&=&(a_3''a'+b_3''a)b[v^p_{\gamma,3},[f_4,v_{\d,4}^q]]-
(a_3''b'+b_3''b)a[v^q_{\gamma,3},[f_4,v_{\d,4}^p]]\\
&=&0.
\end{eqnarray*}
Also
\begin{eqnarray*}
[e_4,D^{p,q}_{\gamma,\d}]&=&[e_4,[a'[f_3,v^p_{\gamma,3}]+a [f_4,v_{\gamma,4}^p],b'[f_3,v^q_{\d,3}]+b[f_4,v_{\d,4}^q]]]\\
&=&(a_4''a'+b_4''a)[v^p_{\gamma,4},b'[f_3,v^q_{\d,3}]+b[f_4,v_{\d,4}^q]]\\
&-&(a_4''b'+b_4''b)[v^q_{\gamma,4},a'[f_3,v^p_{\d,3}]+a[f_4,v_{\d,4}^p]]\\
&=&(a_4''a'+b_4''a)b'[v^q_{\gamma,4},[f_3,v_{\d,3}^q]]-
(a_4''b'+b_4''b)a'[v^q_{\gamma,4},[f_3,v_{\d,3}^p]]]\\
&=&0.
\end{eqnarray*}
Using (R9) together with (\ref{fin8}), we get that
\begin{equation} [v^p,v^q_{\pm j}]=[v^q,v^p_{\pm j}]; \;n_{\ell+1}\leq j\leq n_s.\label{final}
\end{equation}  We note that in the cases under consideration, $(\v^q_\gamma)_0$ is a one dimensional subspace of
$\v^q_\gamma$. Let $n_\ell+1\leq i\leq \ell,$ then there is
$x\in\bbbc$ such that
$$x[f_i,v_{\gamma,i}^q]=[f_{\ell},v_{\gamma,\ell}^q]\andd x[f_i,v_{\d,i}^p]=[f_{\ell},v_{\d,\ell}^p].$$ Now using
Propositions \ref{y3.17}, \ref{m8}($iii$) and (\ref{final}), we
have
$[f_{n_s},v_\d^p,v_{\gamma,i}^q]=[f_{n_s},v_\gamma^q,v_{\d,i}^p].$
This together with the Jacobi identity implies that {{
{\small\begin{eqnarray}
[e_i,[f_{\ell},v_{\d,\ell}^p],[f_{\ell},v_{\gamma,\ell}^q]]\nonumber
\hspace{-2mm}&=&\hspace{-2mm}
[[f_{n_s},v_\d^p],[e_i,[f_{n_s},v_\gamma^q]]]-[[f_{n_s},v_\gamma^q],[e_i,[f_{n_s},v_\d^p]]]\nonumber\\
&=&\hspace{-2mm}
x[[f_{n_s},v_\d^p],[e_i,[f_i,v_{\gamma,i}^q]]-x[[f_{n_s},v_\gamma^q],[e_i,[f_i,v_{\d,i}^p]]\nonumber\\
&=&\hspace{-2mm}
x[[f_{n_s},v_\d^p],[[e_i,f_i],v_{\gamma,i}^q]]-x[[f_{n_s},v_\gamma^q],[[e_i,f_i],v_{\d,i}^p]]\label{cha}\\
&=&\hspace{-2mm}
2x[[f_{n_s},v_\d^p],v_{\gamma,i}^q]\hspace{-1mm}-\hspace{-1mm}2x[[f_{n_s},v_\gamma^q],v_{\d,i}^p]\nonumber\\
&=&\hspace{-2mm}
2x([f_{n_s},v_\d^p,v_{\gamma,i}^q]\hspace{-1mm}-\hspace{-1mm}[v_\d^p,f_{n_s},v_{\gamma,i}^q]\hspace{-1mm}-\hspace{-1mm}
[f_{n_s},v_\gamma^q,v_{\d,i}^p]\hspace{-1mm}+\hspace{-1mm}[v_\gamma^q,f_{n_s},v_{\d,i}^p])
\nonumber\\
&=&\hspace{-2mm}
2x(-[v_\d^p,f_{n_s},v_{\gamma,i}^q]+[v_\gamma^q,f_{n_s},v_{\d,i}^p]).
\nonumber
\end{eqnarray}}}}
 Now if $\a_{\ell}-\a_i$ is not a short root then
$[f_{n_s},v_{\d,i}^p]=[f_{n_s},v_{\gamma,i}^q]=0$ and so we are
done, otherwise there are $n_\ell+1\leq k\leq n_s$ and $y\in\bbbc$
such that
$$[f_{n_s},v_{\d,i}^p]=y v^p_{\d,\pm k}\andd [f_{n_s},v_{\gamma,i}^q]=y v^q_{\gamma,\pm k}.$$ This together with
(\ref{cha}), Propositions \ref{y3.17}, \ref{m8}($iii$) and
(\ref{final}) implies that
\begin{eqnarray*}
[e_i,D_{\gamma,\d}^{p,q}]\in
\bbbc[e_i,[f_{\ell},v_{\d,\ell}^p],[f_{\ell},v_{\gamma,\ell}^q]]=\bbbc(-[v_\d^p,v_{\gamma,\pm
k}^q]+[v_\gamma^q,v_{\d,\pm k}^p])=0.
\end{eqnarray*}
Using the same argument as above, one can show that
$$[f_i,D_{\gamma,\d}^{p,q}]=0.$$

\medskip

($xi$)  Let
$i\in J_\ell$ and $\gamma,\d\in\zn$ and fix $j\in J_\ell,$ $r\in
J_m$ and $a\in J_\nu.$ We need to prove
$$[e_j,[h_{i,\gamma},h_{i,\d}]]=[f_j,[h_{i,\gamma},h_{i,\d}]]=[h_{j,a}^\pm,[h_{i,\gamma},h_{i,\d}]]=
[v^p,[h_{i,\gamma},h_{i,\d}]]=0.$$ Using  the Jacobi identity
together with Proposition \ref{all2}($ii$), we have
{\small\begin{eqnarray*} [e_j,[h_{i,\gamma},h_{i,\d}]]=[h_{i,\d},
h_{i,\gamma},e_j]-[
h_{i,\gamma},h_{i,\d},e_j]&=&\a_j(h_i)([h_{i,\d},e_{j,\gamma}]-[h_{i,\gamma},e_{j,\d}])\\
&=&(\a_j(h_i))^2(e_{j,\gamma+\d}-e_{j,\gamma+\d})=0.
\end{eqnarray*}}
The proof for the second term in the statement is similar and the
last assertion in the statement is immediate using Proposition
\ref{all2}($vi$). Now it remains to prove
$[h_{j,a}^\pm,[h_{i,\gamma},h_{i,\d}]]=0,$ for this we first prove
$[h_{i,a}^\pm,h_{i,b}^\pm,h_{i,\d}]=0$ for all $b\in J_\nu.$ Fix
$b\in J_\nu$ and use Proposition \ref{all2}($v$),($iii$)
 to get
\begin{eqnarray*}
[f_{i,a}^\pm,h_{i,b}^\pm,h_{i,\d}]&=&2[f_{i,a}^\pm,[e_{i,b}^\pm,f_{i,\d}]-h_{i,\d+\gamma_b^\pm}]\\
&=&2[f_{i,a}^\pm,[e_{i,b}^\pm,f_{i,\d}]]-4f_{i,\d+\gamma_a^\pm+\gamma_b^\pm}\\
&=&2[f_{i,\d},e_{i,b}^\pm,f_{i,a}^\pm]-2[e_{i,b}^\pm,f_{i,\d},f_{i,a}^\pm]-4f_{i,\d+\gamma_b^\pm+\gamma_a^\pm}\\
&=&2[f_{i,\d},e_{i,b}^\pm,f_{i,a}^\pm]-4f_{i,\d+\gamma_b^\pm+\gamma_a^\pm}.
\end{eqnarray*}
But again using Proposition \ref{all2}($v$),($iii$), we obtain
\begin{eqnarray*}
2[f_{i\d},[e_{i,b}^\pm,f_{i,a}^\pm]]&=&2[f_{i\d},[e_{i,b}^\pm,f_{i,\gamma^\pm_a}]]\\
&=&[f_{i,\d},[h_{i,b}^\pm,h_{i,a}^\pm]]+2[f_{i,\d},h_{i,\gamma_a^\pm+\gamma_b^\pm}]\\
&=&[h_{i,a}^\pm,h_{i,b}^\pm,f_{i,\d}]-[h_{i,b}^\pm,h_{i,a}^\pm,f_{i,\d}]+4f_{i,\gamma_b^\pm+\gamma_b^\pm+\d}\\
&=&-2[h_{i,a}^\pm,f_{i,\d+\gamma_b^\pm}]+2[h_{i,b}^\pm,f_{i,\d+\gamma_a^\pm}]+4f_{i,\gamma_b^\pm+\gamma_j^\pm+\d}\\
&=&4f_{i,\d+\gamma_a^\pm+\gamma_b^\pm}-4f_{i,\d+\gamma_a^\pm+\gamma_b^\pm}+4f_{i,\gamma_b^\pm+\gamma_a^\pm+\d}\\
&=&4f_{i,\gamma_b^\pm+\gamma_a^\pm+\d}.
\end{eqnarray*}
Therefore $$[f_{i,a}^\pm,h_{i,b}^\pm,h_{i,\d}]=0$$ and so by
(\ref{m1}), the Jacobi identity and the first part of the proof,
we have
\begin{eqnarray*}
[h_{i,a}^\pm,h_{i,b}^\pm,h_{i,\d}]&=&[[e_i,f_{i,a}^\pm],h_{i,b}^\pm,h_{i,\d}]\\
&=&[e_i,f_{i,a}^\pm,h_{i,b}^\pm,h_{i,\d}]-[f_{i,a}^\pm,e_i,h_{i,b}^\pm,h_{i,\d}]\\
&=&0-0=0.
\end{eqnarray*}
Now using  Proposition \ref{all2}($iv$), our obtained information
and the fact that (\ref{3.1'}) is a generating set for the Lie
algebra, one can conclude  that $[h_{j,b}^\pm,h_{i,\d}]\in Z(\Lt)$
for all $b\in J_\nu.$ Similarly one can get
$[h_{j,b}^\pm,h_{i,\gamma}]\in Z(\Lt).$ Therefore the Jacobi
identity implies that
\begin{eqnarray*}
[h_{j,a}^\pm,[h_{i,\gamma},h_{i,\d}]]=[h_{i,\gamma},[h_{j,a}^\pm,h_{i,\d}]]-[h_{i,\d},[h_{j,a}^\pm,h_{i,\gamma}]]=0.
\end{eqnarray*}
This completes the proof.\qed
\medskip
\section{Appendix}\label{app}
\subsection{The exceptional simple Lie algebra of type $G_2$ }\label{app1}

Suppose that $\fc$ is the $8$-dimensional (non-associative) Octonion algebra over $\bbbc.$ It is a
celebrated theorem  that  all the derivations of $\fc$   are inner
and $\fg:=\hbox{Der}(\fc)$ is the $14$-dimensional simple Lie algebra  of type
$G_2$. Each element
of $\fc$ can be represented in a matrix form as
$\left(\begin{array}{cc}a&v\\w&b\end{array}\right)$ where
$a,b\in\bbbc$ and $v,w\in\bbbc^3.$ The add and the scalar product
on $\fc$ is as usual for matrices and the product on $\fc$ is
given by
$$\left(\begin{array}{cc}a&v\\w&b\end{array}\right)\left(\begin{array}{cc}a'&v'\\w'&b'\end{array}\right)
=\left(\begin{array}{cc}aa'-v\cdot w'&av'+b'v+w\times
w'\\a'w+bw'+v\times v'&bb'-w\cdot v'\end{array}\right)$$ where
$"\cdot"$ and $"\times"$ are respectively the  usual inner product
and cross product on $\bbbc^3.$ Take $\{e_i\mid 1\leq i\leq 3\}$
to be the standard basis for $\bbbc^3$ and fix the basis
$\{c_1,\ldots,c_8\}$ for $\fc$ as follows:

$$c_1:=\left(\begin{array}{cc}1&0\\0&0\end{array}\right),
c_2:=\left(\begin{array}{cc}0&0\\0&1\end{array}\right),c_{2+i}:=\left(\begin{array}{cc}0&e_i\\0&0\end{array}\right),
c_{5+i}:=\left(\begin{array}{cc}0&0\\e_i&0\end{array}\right)$$ For
$i=1,2,3.$ The normalized trace  on $\fc$ is
$$t:\fc\longrightarrow
\bbbc;\;\;\left(\begin{array}{cc}a&v\\w&b\end{array}\right)\mapsto
(1/2)(a+b).$$ Set
$$h_1:=D_{c_3,c_6}\andd h_2:=D_{c_4,c_7};\;\fh:=\bbbc h_1\op\bbbc h_2$$ and define $\a_1,\a_2\in\fh^\star$ by
$$\begin{array}{c}
\a_1:\;\;h_1\mapsto 2,\;h_2\mapsto -3\\
\a_2:\;\;h_1\mapsto -1,\;h_2\mapsto 3.
\end{array}
$$
One knows that $\fh$ is a Cartan subalgebra of $\fg$ and that $\fg$ has
a weight space decomposition $\fg=\op_{\a\in R}\fg_\a$ with respect to
$\fh$ where
$$R:=\{0,\pm\a_1,\pm\a_2,\pm(\a_1+\a_2),\pm(2\a_1+\a_2),\pm(3\a_1+\a_2),\pm(3\a_1+2\a_2)\}$$
and
\begin{equation*}\label{v6}
\begin{array}{lll}
\hbox{\fbox{$\fg_{\a_1}=\bbbc D_{c_5,c_4},$}}&\fg_{-\a_1}=\bbbc D_{c_1,c_3},&
\fg_{\a_2}=\bbbc D_{c_3,c_7},\\
\fg_{-\a_2}=\bbbc D_{c_4,c_6},&
\fg_{\a_1+\a_2}=\bbbc D_{c_1,c_7},&\fg_{-(\a_1+\a_2)}=\bbbc D_{c_1,c_4}\\
\fg_{2\a_1+\a_2}=\bbbc D_{c_1,c_5},&\fg_{-(2\a_1+\a_2)}=\bbbc
D_{c_1,c_8},& \fg_{3\a_1+\a_2}=\bbbc D_{c_5,c_6},\\
\fg_{-(3\a_1+\a_2)}=\bbbc D_{c_1,c_3},& \fg_{3\a_1+2\a_2}=\bbbc
D_{c_5,c_7},&\fg_{-(3\a_1+2\a_2)}=\bbbc D_{c_4,c_8}.
\end{array}
\end{equation*}

The subspace $\fc_0$  of $\fc$ consisting of all elements of trace
$0,$ is an irreducible $\fg-$module whose highest weight is a short root.
Indeed setting $\v:=\fc_0,$ we have the weight space decomposition
$\v=\op_{\a\in R_{sh}\cup\{0\}}\v_{\a}$ with $\v_0=c_1-c_2$ and
\begin{equation*}\label{v7}
\begin{array}{lll}
\v_{\a_1}=\bbbc c_6, &\v_{-\a_1}=\bbbc c_3,&\v_{\a_1+\a_2}=\bbbc c_7,\\
\v_{-(\a_1+\a_2)}=\bbbc c_4,&\v_{2\a_1+\a_2}=\bbbc
c_5,&\v_{-(2\a_1+\a_2)}=\bbbc c_8.
\end{array}
\end{equation*}

Now we would like to state some facts which we will use for our presentation.
 We note that $c_5$ is a maximal vector of weight
 $2\a_1+\a_2$ and  $c_8$
is a weight vector of weight $-(2\a_1+\a_2).$ We also have
\begin{equation}\label{fpre1}
\begin{array}{c}
t(c_5,c_t)=0\hbox{ for $3\leq t\leq 7$ with $t\neq5$}\andd\\
\begin{array}{llll}
c_5*c_3=c_7,&\hbox{\fbox{$c_5*c_4=-c_6,$}}&c_5*c_6=0,&c_5*c_7=0.
\end{array}\end{array}\end{equation}

\subsection{The exceptional simple Lie algebra of type $F_4$}\label{exceptional}
Take $\fc_{3\times 3}$ to be the algebra of all matrices with
entries in the Octonion algebra $\fc.$ Consider the conjugate transpose involution
$x\mapsto \bar{x}^t,\;x\in \fc_{3\times 3}$ on $\fc_{3\times 3}$
where $\bar{x}$ is  induced  from the usual involution on $\fc$
and take $\fcj$ to be the self-adjoint elements of $\fc_{3\times
3}$ with respect to the involution. Then $\fcj$ is an exceptional
simple Jordan algebra under the multiplication $x\cdot
y:=(1/2)(xy+yx)$ where the product by juxtaposition is the
multiplication in $\fc_{3\times 3}.$ It is well known fact that
all derivations of $\fcj$ are inner and $\fg:=\hbox{Der}(\fcj)$ is
a simple Lie algebra of type $F_4.$ Each element  of $\fcj$ is
of the form $\left(\begin{array}{ccc} \zeta_1&c& b\\
\bar c&\zeta_2& a\\
\bar b&\bar a&\zeta_3
 \end{array}\right)$
with $\zeta_1,\zeta_2,\zeta_3\in\bbbc$ and $a,b,c\in\fc.$ Such an
element can be denoted  as $(\zeta_1,\zeta_2,\zeta_3,a,b,c).$ One
can see that $(1,1,1,0,0,0)$ is the unit of the algebra $\fcj.$
The map
\begin{equation}\label{talb}
t:\fcj\longrightarrow\bbbc;\;
(\zeta_1,\zeta_2,\zeta_3,a,b,c)\mapsto(1/3)(\zeta_1+\zeta_2+\zeta_3).
\end{equation}
is a normalized trace form. We fix  a basis for $\fcj$ as follows:
\begin{equation}\label{albbas}
\begin{array}{ll}
\fj_1:=(1,0,0,{\bf0},{\bf0},{\bf0}),&\fj_2:=(0,1,0,{\bf0},{\bf0},{\bf0}),\\
\fj_3:=(0,0,1,{\bf0},{\bf0},{\bf0}),&
\fj_{3+i}:=(0,0,0,c_i,{\bf0},{\bf0}),\\
\fj_{11+i}:=(0,0,0,{\bf0},c_i,{\bf0}),&
\fj_{19+i}:=(0,0,0,{\bf0},{\bf0},c_i);\;1\leq i\leq 8.
\end{array}
\end{equation}

Set
\begin{equation}\label{caralb}
\dot h_1:=D_{\fj_4,\fj_5},\;\dot h_2:=D_{\fj_6,\fj_9}, \;\dot
h_3:=D_{\fj_7,\fj_10},\;\dot
h_4:=D_{\fj_8,\fj_{11}}\andd\fh:=\sum_{i=1}^4\bbbc\dot h_i,
\end{equation}
then $\fh$ is a Cartan subalgebra of $\fcj.$ Next define $\a_1,\a_2,\a_3,\a_4\in\fh^\star$ by
\begin{equation}\label{v8}
\begin{array}{l}
\a_1:\;\;\dot h_1\mapsto0, \dot
h_2\mapsto\frac{1}{2},\dot h_3\mapsto\frac{-1}{2},\dot
h_4\mapsto0,\\
\a_2:\;\;\dot h_1\mapsto0,\dot
h_2\mapsto0,\dot h_3\mapsto\frac{1}{2},\dot
h_4\mapsto\frac{-1}{2},\\
\a_3:\;\;\dot h_1\mapsto0,\dot
h_2\mapsto0,\dot h_3\mapsto0,\dot
h_4\mapsto\frac{1}{2},\\
\a_4:\;\;\dot h_1\mapsto\frac{1}{4},\dot
h_2\mapsto\frac{-1}{4},\dot h_3\mapsto\frac{-1}{4},\dot
h_4\mapsto\frac{-1}{4}.
\end{array}
\end{equation}
We have the weight space decomposition $\fg=\op_{\a\in R}\fg_{\a}$
with respect to $\fh$ where $R=\{0\}\cup R_{sh}\cup R_{lg}$ where
$$\begin{array}{rl}
R_{sh}=&\hspace{-2mm}\pm\{\a_3,\a_4,\a_2+\a_3,\a_3+\a_4,\a_1+\a_2+\a_3,\a_2+\a_3+\a_4,\\&
 \a_2+2\a_3+\a_4,\a_1+\a_2+\a_3+\a_4,\a_1+\a_2+2\a_3+\a_4,\\&
 \a_1+2\a_2+2\a_3+\a_4,\a_1+2\a_2+3\a_3+\a_4,\a_1+2\a_2+3\a_3+2\a_4
 \},\vspace{2mm}\\
  R_{lg}=&\hspace{-2mm}\pm
 \{0,\a_1,\a_2, \a_1+\a_2,
\a_2+2\a_3,\a_1+\a_2+2\a_3,\a_1+2\a_2+2\a_3,
\\&
\a_2+2\a_3+2\a_4,
\a_1+\a_2+2\a_3+2\a_4,
\a_1+2\a_2+2\a_3+2\a_4,\\&\a_1+2\a_2+4\a_3+2\a_4,
\a_1+3\a_2+4\a_3+2\a_4,
 2\a_1+3\a_2+4\a_3+2\a_4\},
\end{array}$$ and
$$
\begin{array}{ll}
\fg_{\a_1}=\bbbc D_{\fj_6,\fj_{10}},& \fg_{-\a_1}=\bbbc
D_{\fj_7,\fj_9},\\
\fg_{\a_2}=\bbbc D_{\fj_7,\fj_{11}},& \fg_{-\a_2}=\bbbc
D_{\fj_8,\fj_{10}},\\
\fg_{\a_3}=\bbbc D_{\fj_1,\fj_8},& \fg_{-\a_3}=\bbbc
D_{\fj_1,\fj_{11}},\\
\fg_{\a_4}=\bbbc D_{\fj_2,\fj_{21}},&
\fg_{-\a_4}=\bbbc
D_{\fj_2,\fj_{20}},\\
\fg_{\a_1+\a_2}=\bbbc D_{\fj_6,\fj_{11}},&
\fg_{-\a_1-\a_2}=\bbbc D_{\fj_8,\fj_9},\\
\fg_{\a_2+\a_3}=\bbbc D_{\fj_1,\fj_7},&
\fg_{-\a_2-\a_3}=\bbbc
D_{\fj_1,\fj_{10}},\\
\fg_{\a_2+2\a_3}=\bbbc D_{\fj_7,\fj_8},&
\fg_{-(\a_2+2\a_3)}=\bbbc
D_{\fj_{10},\fj_{11}},\\
\hbox{\fbox{$\fg_{\a_3+\a_4}=\bbbc D_{\fj_4,\fj_{24},}$}}&
\fg_{-(\a_3+\a_4)}=\bbbc
D_{\fj_1,\fj_{19}},\\
\fg_{\a_1+\a_2+\a_3}=\bbbc
D_{\fj_1,\fj_6},&\fg_{-(\a_1+\a_2+\a_3)}=\bbbc D_{\fj_1,\fj_2},\\
\fg_{\a_1+\a_2+2\a_3}=\bbbc
D_{\fj_6,\fj_8},&
\fg_{-(\a_1+\a_2+2\a_3)}=\bbbc D_{\fj_9,\fj_{11}},\\
\fg_{\a_1+2\a_2+2\a_3}=\bbbc
D_{\fj_6,\fj_7},&
\fg_{-(\a_1+2\a_2+2\a_3)}=\bbbc D_{\fj_9,\fj_{10}},\\
 \fg_{\a_2+\a_3+\a_4}=\bbbc D_{\fj_1,\fj_{15}},&
\fg_{-(\a_2+\a_3+\a_4)}=\bbbc D_{\fj_1,\fj_{18}},\\
\fg_{\a_2+2\a_3+\a_4}=\bbbc D_{\fj_2,\fj_{25}},&
\fg_{-(\a_2+2\a_3+\a_4)}=\bbbc
D_{\fj_2,\fj_{22}},\\
\fg_{\a_2+2\a_3+2\a_4}=\bbbc D_{\fj_4,\fj_9},&
\fg_{-(\a_2+2\a_3+2\a_4)}=\bbbc D_{\fj_5,\fj_6},
\\
\fg_{\a_1+\a_2+\a_3+\a_4}=\bbbc D_{\fj_1,\fj_{14}},&
\fg_{-(\a_1+\a_2+\a_3+\a_4)}=\bbbc
D_{\fj_1,\fj_{17}},\\
\fg_{\a_1+\a_2+2\a_3+\a_4}=\bbbc D_{\fj_2,\fj_{26}},&
\fg_{-(\a_1+\a_2+2\a_3+\a_4)}=\bbbc D_{\fj_2,\fj_{23}},\\
\fg_{\a_1+\a_2+2\a_3+2\a_4}=\bbbc D_{\fj_4,\fj_{10}},&
\fg_{-(\a_1+\a_2+2\a_3+2\a_4)}=\bbbc
D_{\fj_5,\fj_7},\\
\fg_{\a_1+2\a_2+2\a_3+\a_4 }=\bbbc D_{\fj_2,\fj_{27}},&
\fg_{-(\a_1+2\a_2+2\a_3+\a_4)}=\bbbc
D_{\fj_2,\fj_{24}},\\
\fg_{\a_1+2\a_2+2\a_3+2\a_4}=\bbbc D_{\fj_4,\fj_{11}},&
\fg_{-(\a_1+2\a_2+2\a_3+2\a_4)}=\bbbc
D_{\fj_5,\fj_8},\\
\fg_{\a_1+2\a_2+3\a_3+\a_4}=\bbbc
D_{\fj_1,\fj_{13}},&
\fg_{-(\a_1+2\a_2+3\a_3+\a_4)}=\bbbc D_{\fj_1,\fj_{12}},\\
\fg_{\a_1+2\a_2+3\a_3+2\a_4}=\bbbc
D_{\fj_1,\fj_4},&\fg_{-(\a_1+2\a_2+3\a_3+2\a_4)}=\bbbc
D_{\fj_1,\fj_5},\\
\fg_{\a_1+2\a_2+4\a_3+2\a_4}=\bbbc D_{\fj_4,\fj_8},&
\fg_{-(\a_1+2\a_2+4\a_3+2\a_4)}=\bbbc D_{\fj_5,\fj_{11}},\\
\fg_{\a_1+3\a_2+4\a_3+2\a_4}=\bbbc D_{\fj_4,\fj_7},&
\fg_{-(\a_1+3\a_2+4\a_3+2\a_4)}=\bbbc D_{\fj_5,\fj_{10}},\\
\fg_{2\a_1+3\a_2+4\a_3+2\a_4}=\bbbc
D_{\fj_4,\fj_6},&\fg_{-(2\a_1+3\a_2+4\a_3+2\a_4)}=\bbbc
D_{\fj_5,\fj_9}
\end{array}
$$

 The subspace $\v:=\fcj_0$ of trace zero elements of $\fcj$ is an irreducible  $\fg$-module whose highest weight is a short root, with the weight space decomposition, with respect to $\fh$,  $\v=\op_{\a\in
R_{sh}\cup\{0\}}\v_\a$, where $\v_0:=\bbbc(\fj_1-\fj_2)\op\bbbc(\fj_2-\fj_3)$ and
\begin{equation}\label{wsdal}
\begin{array}{ll}
\v_{\a_3}=\bbbc \fj_{8}& \v_{-\a_3}=\bbbc
\fj_{11}\\
\v_{\a_4}=\bbbc \fj_{21}& \v_{-\a_4}=\bbbc
\fj_{20}\\
\v_{\a_2+\a_3}=\bbbc \fj_{7}& \v_{-(\a_2+\a_3)}=\bbbc
\fj_{10}\\
\v_{\a_3+\a_4}=\bbbc \fj_{16}& \v_{-(\a_3+\a_4)}=\bbbc
\fj_{19}\\
\v_{\a_1+\a_2+\a_3}=\bbbc \fj_{6}& \v_{-(\a_1+\a_2+\a_3)}=\bbbc
\fj_{9}\\
\v_{\a_2+\a_3+\a_4}=\bbbc \fj_{15}&\v_{-(\a_2+\a_3+\a_4)}=\bbbc
\fj_{18}\\
\v_{\a_2+2\a_3+\a_4}=\bbbc \fj_{25}&\v_{-(\a_2+2\a_3+\a_4)}=\bbbc
\fj_{22}\\
\v_{\a_1+\a_2+\a_3+\a_4}=\bbbc
\fj_{14}&\v_{-(\a_1+\a_2+\a_3+\a_4)}=\bbbc
\fj_{17}\\
\v_{\a_1+\a_2+2\a_3+\a_4}=\bbbc
\fj_{26}&\v_{-(\a_1+\a_2+2\a_3+\a_4)}=\bbbc \fj_{23}\\
\v_{\a_1+2\a_2+2\a_3+\a_4}=\bbbc
\fj_{27}&\v_{-(\a_1+2\a_2+2\a_3+\a_4)}=\bbbc
\fj_{24}\\
\v_{\a_1+2\a_2+3\a_3+\a_4}=\bbbc
\fj_{13}&\v_{-(\a_1+2\a_2+3\a_3+\a_4)}=\bbbc
\fj_{12}\\\v_{\a_1+2\a_2+3\a_3+2\a_4}=\bbbc
\fj_{4}&\v_{-(\a_1+2\a_2+3\a_3+2\a_4)}=\bbbc \fj_5.
\end{array}
\end{equation}
Now we state some information suitable for our presentation.
 We note that $\fj_4$ is a maximal vector of weight
 $\a_1+2\a_2+3\a_3+2\a_4$ and  $\fj_5$
is a weight vector of weight $-(\a_1+2\a_2+3\a_3+2\a_4).$ We also
have
\begin{equation}\label{fpre2}
\begin{array}{c}
t(\fj_4,\fj_t)=0\hbox{ for $6\leq t\leq 27$ }\andd\\
\begin{array}{lllll}
\fj_4*\fj_6=0&\fj_4*\fj_7=0&\fj_4*\fj_8=0&\fj_4*\fj_9=0&\fj_4*\fj_{10}=0\\
\fj_4*\fj_{11}=0&\fj_4*\fj_{12}=\frac{1}{2}\fj_{20}&\fj_4*\fj_{13}=0&\fj_4*\fj_{14}=0&\fj_4*\fj_{15}=0\\
\fj_4*\fj_{16}=0&\fj_4*\fj_{17}=\frac{-1}{2}\fj_{24}&\fj_4*\fj_{18}=\frac{-1}{2}\fj_{25}&\fj_4*\fj_{19}=0&\fj_4*\fj_{20}=\frac{1}{2}\fj_{12}\\
\fj_4*\fj_{21}=0&\fj_4*\fj_{22}=\frac{-1}{2}\fj_{13}&
\fj_4*\fj_{23}=\frac{-1}{2}\fj_{15}&\hbox{\fbox{$\fj_4*\fj_{24}=\frac{-1}{2}\fj_{16}$}}&\fj_4*\fj_{25}=0\\
&\fj_4*\fj_{26}=0&\fj_4*\fj_{27}=0.
\end{array}
\end{array}
\end{equation}

Consider the normalized trace on $\fcj$ and take $\zeta=e^{2\pi i/3}$. We have
\begin{equation}\label{rel1}
\begin{array}{c}
t\big(4D_{\fj_1,\fj_{11}}(\fj_8)-4\zeta D_{\fj_2,\fj_{20}}(\fj_{21}),
4D_{\fj_1,\fj_{11}}(\fj_8)-4\zeta^2 D_{\fj_2,\fj_{20}}(\fj_{21})
\big)=1,\vspace{2mm}\\
 \big(4D_{\fj_1,\fj_{11}}(\fj_8)-4\zeta
D_{\fj_2,\fj_{20}}(\fj_{21})\big)*
\big(4D_{\fj_1,\fj_{11}}(\fj_8)-4\zeta^2
D_{\fj_2,\fj_{20}}(\fj_{21})\big)=0,
\end{array}
\end{equation}
and
\begin{equation}\label{rel2}
\begin{array}{ll}
D_{\fj_1,\fj_8}(D_{\fj_1,\fj_{11}}(\fj_8))=\frac{-1}{8}\fj_8&
D_{\fj_1,\fj_8}(D_{\fj_2,\fj_{20}}(\fj_{21}))=\frac{1}{16}\fj_8\\
D_{\fj_2,\fj_{21}}(D_{\fj_1,\fj_{11}}(\fj_8))=\frac{1}{16}\fj_{21}&
D_{\fj_2,\fj_{21}}(D_{\fj_2,\fj_{20}}(\fj_{21}))=\frac{-1}{8}\fj_{21}.
\end{array}
\end{equation}

Now setting
$$a:=-4\zeta,\;a':=4,\;b:=-4\zeta^2,\;b':=4,\;b''_4=a''_3:=\frac{-1}{8},\;a''_4=b''_3:=\frac{1}{16},$$
we have
\begin{equation}\label{rel3}
(a''_3a'+b''_3a)b=-(a''_3b'+b''_3b)a\andd
(a''_4a'+b''_4a)b'=-(a''_4b'+b''_4b)a'.
\end{equation}

\subsection{The algebras $\aa_p,\jj_p (\leq p\leq 3)$}\label{app3}
We  recall these  algebras from
\cite{BGKN} and \cite{AABGP}. Suppose that $\nu$ is a positive integer, for $1\leq p\leq 3$ with $p\leq
\nu,$ take $\aa_0$ to be $A_{[\nu]}$ and
$$\aa_p:=\aa_{p-1}\op \aa_{p-1}x_p$$
to be the algebra obtained from  $\aa_{p-1}$ using the {\it
Cayley-Dickson} process with $$x_p^2=t_p.$$ The last one is called
the  {\it Cayley torus}  (or {\it Octonion torus}). We mention
that the Cayley torus is alternative but not associative, also the
center and the nucleus of this algebra coincides with $A_{[\nu]}.$
We know that $\aa_{p-1}$ is a subalgebra of $\aa_p$ for $1\leq
p\leq 3.$ Moreover these four algebras are algebras over
$A_{[\nu]},$ in fact they are free $A_{[\nu]}$-modules. Next note
that for $1\leq i\leq 8,$ there exist unique
$s_1,s_2,s_3\in\{0,1\}$ such that $i=1+s_1+2s_2+4s_3.$ Take
$w_i:=(x_1^{s_1}x_2^{s_2}) x_3^{s_3}.$ Then $\{w_i\mid 1\leq i\leq
2^p\}$ is an $A_{[\nu]}$-basis for $\aa_p,$ $0\leq p\leq 3.$ We
consider the normalized trace \begin{equation}
\label{v2}T:\aa_p\longrightarrow
A_{[\nu]};\;\sum_{i=1}^{2^p}r_iw_i\mapsto r_1
\end{equation}
on $\aa_p,$ $0\leq p\leq 3$  and define $t$ and $*$ to be as
before.

Next take $\jj_o:=A_{[\nu]},$ $\jj_1$ to be the commutative
associative algebra over $A_{[\nu]}$ with generator $x_1$ subject
to the relation $x_1^3=t_1$ and  $\jj_2$ to be the plus algebra of
the associative  algebra over $A_{[\nu]}$ generated by $x_1,x_2$
subject to the relations $x_1^3=t_1,
x_2^3=t_2,x_1x_2=e^{\frac{2\pi i}{3}}x_2x_1$. Let $\jj_3:=\jj_2\op
(\jj_2\cdot x_3)\op(\jj_2\cdot x_3^2)$  be the Jordan algebra
obtained from $\jj_2$ using Tit's first Jordan algebra
construction \cite[Chapter IX]{J} with $x_3^3=t_3.$ We mention
that whenever we use  $\jj_p,$ $1\leq p\leq 3,$ we  assume
$p\leq \nu.$

One  sees that for $1\leq i\leq 27,$ there exist unique
$s_1,s_2,s_3\in\{0,1,2\}$ such that $i=1+s_1+3s_2+9s_3.$ Take
$w_i:=(x_1^{s_1}\cdot x_2^{s_2})\cdot x_3^{s_3}.$ Then $\{w_i\mid
1\leq i\leq 3^p\}$ is an $A_{[\nu]}$-basis for $\jj_p,$ $0\leq
p\leq 3.$ Define the following normalized  trace  on $\jj_p,$
$0\leq p\leq 3:$
\begin{equation}
\label{v22}T:\jj_p\longrightarrow
A_{[\nu]};\;\sum_{i=1}^{3^p}r_iw_i\mapsto r_1.
\end{equation}

Note that the algebra products on $\aa_p$ and $\jj_p,$ the trace form  on these algebras  and the corresponding $*$-operators  are $A_{[\nu]}$-bilinear. Therefore each of these are uniquely determined by its table of multiplication on the corresponding  $A_{[\nu]}$-basis. For the convenience of the reader, in what follows we accomplish the tables for $\aa_3$ and $\jj_3$ and we note that the tables for $\aa_1,\aa_2$
 and $\jj_1,\jj_2$
are sub-tables  of the ones for $\aa_3$ and $\jj_3$ respectively as  $\aa_1\sub\aa_2\sub\aa_3$ and $\jj_1\sub\jj_2\sub\jj_3.$

\newpage

\begin{table}\caption{the table of the products (for $\aa_p$)}
\begin{tabular}{l}
\begin{tabular}{|l||lp{1.2cm}|p{1.2cm}p{1.2cm}|p{1.2cm}p{1.2cm}p{1.2cm}p{1.2cm}|}
\hline
 &$w_1$&$w_2$&$w_3$&$w_4$&$w_5$&$w_6$&$w_7$&$w_8$\\
\hline
$w_1$&{\bf 1}&$w_2$&$w_3$&$w_4$&$w_5$&$w_6$&$w_7$&$w_8$\\
$w_2$&$w_2$&${\bf t_1}$&$w_4$&$t_1w_3$&$w_6$&$t_1w_5$&$-w_8$&$-t_1w_7$\\
\cline{1-3}
\end{tabular}\\
\begin{tabular}{|l||lp{1.2cm}p{1.2cm}p{1.2cm}|p{1.2cm}p{1.2cm}p{1.2cm}p{1.2cm}|}
$w_3$&$w_3$&$-w_4$&$\bf{t_2}$&$-t_2w_2$&$w_7$&$w_8$&$t_2w_5$&$t_2w_6$\\
$w_4$&$w_4$&$-t_1w_3$&$t_2w_2$&$\bf{-t_1t_2}$&$w_8$&$t_1w_7$&$-t_2w_6$&$-t_1t_2w_5$\\
\cline{1-5}
\end{tabular}\\
\begin{tabular}{|l||lp{1.2cm}p{1.2cm}p{1.2cm}p{1.2cm}p{1.2cm}p{1.2cm}p{1.2cm}|}
$w_5$&$w_5$&$-w_6$&$-w_7$&$-w_8$&$\bf{t_3}$&$-t_3w_2$&$-t_3w_3$&$-t_3w_4$\\
$w_6$&$w_6$&$-t_1w_5$&$-w_8$&$-t_1w_7$&$t_3w_2$&$\bf{-t_1t_3}$&$t_3w_4$&$t_1t_3w_3$\\
$w_7$&$w_7$&$w_8$&$-t_2w_5$&$t_2w_6$&$t_3w_3$&$-t_3w_4$&$\bf{-t_2t_3}$&$-t_2t_3w_2$\\
$w_8$&$w_8$&$t_1w_7$&$-t_2w_6$&$t_1t_2w_5$&$t_3w_4$&$-t_1t_3w_3$&$t_2t_3w_2$&$\bf{t_1t_2t_3}$\\
\hline
\end{tabular}
\end{tabular}
\end{table}

\begin{table}\caption{the table of the trace form (for $\aa_p$)}
\begin{tabular}{l}
\begin{tabular}{|p{1cm}||p{1cm}p{1cm}|p{1cm}p{1cm}|p{1cm}p{1cm}p{1cm}p{1cm}|}
\hline
$t\fm$ &$w_1$&$w_2$&$w_3$&$w_4$&$w_5$&$w_6$&$w_7$&$w_8$\\
\hline
$w_1$&{\bf 1}&0&0&0&0&0&0&0\\
$w_2$&0&${\bf t_1}$&0&0&0&0&0&0\\
\cline{1-3}
\end{tabular}\\
\begin{tabular}{|p{1cm}||p{1cm}p{1cm}p{1cm}p{1cm}|p{1cm}p{1cm}p{1cm}p{1cm}|}
$w_3$&0&0&${\bf t_2}$&0&0&0&0&0\\
$w_4$&0&0&0&$\bf{-t_1t_2}$&0&0&0&0\\
\cline{1-5}
\end{tabular}\\
\begin{tabular}{|p{1cm}||p{1cm}p{1cm}p{1cm}p{1cm}p{1cm}p{1cm}p{1cm}p{1cm}|}
$w_5$&0&0&0&0&$\bf{t_3}$&0&0&0\\
$w_6$&0&0&0&0&0&$\bf{-t_1t_3}$&0&0\\
$w_7$&0&0&0&0&0&0&$\bf{-t_2t_3}$&0\\
$w_8$&0&0&0&0&0&0&0&$\bf{t_1t_2t_3}$\\
\hline
\end{tabular}
\end{tabular}
\end{table}

\begin{table}\caption{the table of the operator $*$ (for $\aa_p$)}
\begin{tabular}{l}
\begin{tabular}{|l||lp{1.2cm}|p{1.2cm}p{1.2cm}|p{1.2cm}p{1.2cm}p{1.2cm}p{1.2cm}|}
\hline
 $*$&$w_1$&$w_2$&$w_3$&$w_4$&$w_5$&$w_6$&$w_7$&$w_8$\\
\hline
$w_1$&0&$w_2$&$w_3$&$w_4$&$w_5$&$w_6$&$w_7$&$w_8$\\
$w_2$&$w_2$&0&$w_4$&$t_1w_3$&$w_6$&$t_1w_5$&$-w_8$&$-t_1w_7$\\
\cline{1-3}
\end{tabular}\\
\begin{tabular}{|l||lp{1.2cm}p{1.2cm}p{1.2cm}|p{1.2cm}p{1.2cm}p{1.2cm}p{1.2cm}|}
$w_3$&$w_3$&$-w_4$&0&$-t_2w_2$&$w_7$&$w_8$&$t_2w_5$&$t_2w_6$\\
$w_4$&$w_4$&$-t_1w_3$&$t_2w_2$&0&$w_8$&$t_1w_7$&$-t_2w_6$&$-t_1t_2w_5$\\
\cline{1-5}
\end{tabular}\\
\begin{tabular}{|l||lp{1.2cm}p{1.2cm}p{1.2cm}p{1.2cm}p{1.2cm}p{1.2cm}p{1.2cm}|}
$w_5$&$w_5$&$-w_6$&$-w_7$&$-w_8$&0&$-t_3w_2$&$-t_3w_3$&$-t_3w_4$\\
$w_6$&$w_6$&$-t_1w_5$&$-w_8$&$-t_1w_7$&$t_3w_2$&0&$t_3w_4$&$t_1t_3w_3$\\
$w_7$&$w_7$&$w_8$&$-t_2w_5$&$t_2w_6$&$t_3w_3$&$-t_3w_4$&0&$-t_2t_3w_2$\\
$w_8$&$w_8$&$t_1w_7$&$-t_2w_6$&$t_1t_2w_5$&$t_3w_4$&$-t_1t_3w_3$&$t_2t_3w_2$&0\\
\hline
\end{tabular}
\end{tabular}
\end{table}

Setting
$$\begin{array}{ll}
b_{2,2}:=1&b_{2}^{3}=b_{3}^{2}:=t_1\\
b_{2}^{4}=b_{4}^{2}:=1& b_{2}^{5}=b_{5}^{2}:=\frac{(1+e^{\pi
i/3})^2}{2(1+e^{2\pi
i/3})}\\b_{2}^{6}=b_{6}^{2}:=\frac{t_1(1+e^{\pi i/3}+e^{2\pi
i/3}+e^{\pi
i})}{4}& b_{2}^{7}=b_{7}^{2}:=1\\
b_{2}^{8}=b_{8}^{2}:=\frac{(1+e^{2\pi i/3})^2}{2(1+e^{4\pi
i/3})}&b_{2}^{9}=b_{9}^{2}:=\frac{t_1(1+e^{2\pi i/3}+e^{4\pi
i/3}+e^{2\pi
i})}{4}\\
b_{2}^{10}=b_{10}^{2}:=\frac{-1}{2}&
b_{2}^{11}=b_{11}^{2}:=\frac{-1}{2}\\
b_{2}^{12}=b_{12}^{2}:=\frac{-t_1}{2}&
b_{2}^{13}=b_{13}^{2}:=\frac{-e^{\pi i/3}}{1+e^{\pi i/3}}\\
b_{2}^{14}=b_{14}^{2}:=\frac{-e^{\pi i/3}(1+e^{\pi
i/3})}{2(1+e^{2\pi i/3})}&b_{2}^{15}=b_{15}^{2}:=\frac{-t_1e^{\pi
i/3}}{4}\\
 b_{2}^{16}=b_{16}^{2}:=\frac{-e^{2\pi i/3}}{1+e^{2\pi i/3}}&
b_{2}^{17}=b_{17}^{2}:=\frac{-e^{2\pi i/3}(1+e^{2\pi
i/3})}{2(1+e^{4\pi i/3})}\\
b_{2}^{18}=b_{18}^{2}:=\frac{-t_1e^{2\pi i/3}(1+e^{4\pi i/3})}{4}&
b_{2}^{19}=b_{19}^{2}:=1\\
b_{2}^{20}=b_{20}^{2}:=\frac{-1}{2}&b_{2}^{21}=b_{21}^{2}:=\frac{-t_1}{2}\\
b_{2}^{22}=b_{22}^{2}:=\frac{-1}{1+e^{\pi i/3}}&
b_{2}^{23}=b_{23}^{2}:=\frac{-(1+e^{\pi i/3})}{2(1+e^{2\pi
i/3})}
\end{array}$$

$$\begin{array}{ll}
b_{2}^{24}=b_{24}^{2}:=\frac{-t_1(1+e^{2\pi i/3})}{4}&
 b_{2}^{25}=b_{25}^{2}:=\frac{-1}{1+e^{2\pi i/3}}\\
b_{2}^{26}=b_{26}^{2}:=\frac{-(1+e^{2\pi i/3})}{2(1+e^{4\pi
i/3})}&b_{2}^{27}=b_{27}^{2}:=\frac{-t_1(1+e^{4\pi
i/3})}{4}\\b_{3}^{3}:=t_1&
b_{3}^{4}=b_{4}^{3}:=1\\b_{3}^{5}=b_{5}^{3}:=\frac{(1+e^{\pi i/3}+
e^{\pi i/3}+e^{\pi
i})t_1}{4}&b_{3}^{6}=b_{6}^{3}:=\frac{(1+e^{2\pi
i/3})^2t_1}{2(1+e^{\pi
i/3})}\\b_{3}^{7}=b_{7}^{3}:=1&b_{3}^{8}=b_{8}^{3}:=\frac{(1+e^{2\pi
i/3}+e^{4\pi i/3}+e^{2\pi
i})t_1}{4}\\b_{3}^{9}=b_{9}^{3}:=\frac{(e^{4\pi
i/3})^2t_1}{2(1+e^{2\pi i/3})}&
b_{3}^{10}=b_{10}^{3}:=\frac{-1}{2}\\
b_{3}^{11}=b_{11}^{3}:=\frac{-t_1}{2}&b_{3}^{12}=b_{12}^{3}:=\frac{-t_1}{2}\\
b_{3}^{13}=b_{13}^{3}:=\frac{-1}{(1+e^{2\pi
i/3})}&b_{3}^{14}=b_{14}^{3}:=\frac{-t_1(1+e^{\pi i/3})}{4}\\
b_{3}^{15}=b_{15}^{3}:=\frac{-t_1(1+e^{2\pi i/3})}{2(1+e^{\pi
i/3})}& b_{3}^{16}=b_{16}^{3}:=\frac{-1}{1+e^{4\pi
i/3}}\\b_{3}^{17}=b_{17}^{3}:=\frac{-(1+e^{2\pi i/3})t_1}{4}&
b_{3}^{18}=b_{18}^{3}:=\frac{-t_1(1+e^{4\pi i/3})}{2(1+e^{2\pi i/3})}\\
b_{3}^{19}=b_{19}^{3}:=1&b_{3}^{20}=b_{20}^{3}=\frac{-t_1}{2}\\b_{3}^{21}=b_{21}^{3}:=\frac{-t_1}{2}&
b_{3}^{22}=b_{22}^{3}:=\frac{-e^{2\pi i/3}}{1+e^{2\pi i/3}}\\
b_{3}^{23}=b_{23}^{3}:=\frac{-e^{2\pi i/3}(1+e^{\pi i/3})t_1}{4}&
b_{3}^{24}=b_{24}^{3}:=\frac{-e^{2\pi i/3}(1+e^{2\pi
i/3})t_1}{2(1+e^{2\pi i/3})}\\
b_{3}^{25}=b_{25}^{3}:=\frac{-e^{4\pi i/3}}{1+e^{4\pi
i/3}}&b_{3}^{26}=b_{26}^{3}:=\frac{-e^{4\pi i/3}(1+e^{2\pi
i/3})t_1}{4}\\b_{3}^{27}=b_{27}^{3}:=\frac{-e^{4\pi i/3}(1+e^{4\pi
i/3})t_1}{2(1+e^{2\pi i/3})}&
b_{4}^{4}:=1\\b_{4}^{5}=b_{5}^{4}:=\frac{(1+e^{\pi
i/3})^2}{2(1+e^{2\pi i/3})}&b_{4}^{6}=b_{6}^{4}:=\frac{ (1+e^{2\pi
i/3}+e^{2\pi i/3})^2}{2(1+e^{4\pi i/3})}\\
b_{4}^{7}=b_{7}^{4}:=t_2&
 b_{4}^{8}=b_{8}^{4}:=
 \frac{(1+e^{\pi
i/3}+e^{2\pi i/3}+e^{\pi i})t_2}{4}\\
b_{4}^{9}=b_{9}^{4}:=\frac{(1+e^{2\pi i/3}+e^{4\pi i/3}+ e^{2\pi
i})t_2}{4}& b_{4}^{10}=b_{10}^{4}:=\frac{-1}{2}\\
b_{4}^{11}=b_{11}^{4}:=\frac{-e^{\pi i/3}}{1+e^{\pi i/3}}&
b_{4}^{12}=b_{12}^{4}:=\frac{-e^{2\pi i/3}}{1+e^{2\pi i/3}}\\
b_{4}^{13}=b_{13}^{4}:=\frac{-1}{2}&
b_{4}^{14}=b_{14}^{4}:=\frac{-e^{\pi i/3}(1+e^{\pi
i/3})}{2(1+e^{2\pi i/3})}\\ b_{4}^{15}= b_{15}^{4}:=\frac{-e^{2\pi
i/3}(1+e^{2\pi i/3})}{2(1+e^{4\pi i/3})}&
 b_{4}^{16}=b_{16}^{4}:=\frac{-t_2}{2}\\
b_{4}^{17}=b_{17}^{4}:=\frac{-e^{\pi i/3}(1+e^{2\pi i/3})t_2}{4})&
b_{4}^{18}=b_{18}^{4}:= \frac{-e^{2\pi i/3}(1+e^{4\pi
i/3})t_2}{4})\\
 b_{4}^{19}=b_{19}^{4}:=1&
b_{4}^{20}=b_{20}^{4}:=\frac{-1}{1+e^{\pi i/3}}\\
b_{4}^{21}=b_{21}^{4}:=\frac{-1}{1+e^{2\pi i/3}}&
b_{4}^{22}=b_{22}^{4}:=\frac{-1}{2}\\
b_{4}^{23}=b_{23}^{4}:=\frac{-(1+e^{\pi i/3})}{2(1+e^{2\pi
i/3})}&b_{4}^{24}=b_{24}^{4}:=\frac{-(1+e^{2\pi i/3})}{2(1+e^{4\pi i/3})}\\
 b_{4}^{25}=b_{25}^{4}:=\frac{-t_2}{2}&
b_{4}^{26}=b_{26}^{4}:=\frac{-(1+e^{2\pi i/3})t_2}{4}\\
b_{4}^{27}=b_{27}^{4}:=\frac{-(1+e^{4\pi i/3})t_2}{4}&
b_{5}^{5}:=\frac{e^{\pi i/3}(1+e^{\pi i/3})^2}{2(1+e^{4\pi
i/3})}\\
b_{5}^{6}=b_{6}^{5}:=\frac{t_1(1+e^{\pi i/3})^2(e^{\pi i/3}+e^{\pi
i})}{8}& b_{5}^{7}=b_{7}^{5}
:=\frac{(1+e^{\pi i/3}+e^{2\pi i/3}+e^{\pi i})}{4}\\
b_{5}^{8}=b_{8}^{5}:= \frac{t_2(1+e^{\pi i/3})^2(e^{\pi
i/3}+e^{\pi i})}{8}& b_{5}^{9}=b_{9}^{5} :=\frac{t_1t_2e^{2\pi
i/3}(1+e^{\pi i/3})(1+e^{4\pi i/3}}{4}\\
b_{5}^{10}=b_{10}^{5}:=\frac{-1}{2}&
b_{5}^{11}=b_{11}^{5}:=\frac{-e^{\pi i/3}(1+e^{\pi
i/3})}{2(1+e^{2\pi i/3})}\\
b_{5}^{12}=b_{12}^{5}:=\frac{-t_1e^{2\pi i/3}(1+e^{\pi i/3})}{4}&
b_{5}^{13}=b_{13}^{5}:=\frac{-(1+e^{\pi i/3})}{2(1+e^{2\pi i/3})}\\
b_{5}^{14}=b_{14}^{5}:=\frac{-e^{\pi i/3} (1+e^{\pi
i/3})^2}{4(1+e^{4\pi i/3})}&b_{5}^{15}=b_{15}^{5}
:=\frac{-t_1e^{2\pi i/3}(1+e^{\pi i/3})(1+e^{2\pi i/3})}{8}\\
 b_{5}^{16}=b_{16}^{5}
 := \frac{-t_2(1+e^{\pi i/3})}{4}&
b_{5}^{17}=b_{17}^{5} :=\frac{-t_2e^{\pi i/3}(1+e^{\pi
i/3})(1+e^{2\pi i/3})}{8}\\ b_{5}^{18}=b_{18}^{5}
:=\frac{-t_1t_2e^{2\pi i/3}(1+e^{\pi i/3})(1+e^{4\pi i/3})}{8}&
b_{5}^{19}=b_{19}^{5}:=1
\end{array}$$

$$\begin{array}{ll}
b_{5}^{20}=b_{20}^{5}:=\frac{-(1+e^{\pi i/3})}{2(1+e^{2\pi
i/3})}&b_{5}^{21}=b_{21}^{5}:=\frac{-t_1(1+e^{\pi/3})}{4}\\
b_{5}^{22}=b_{22}^{5}:=\frac{-e^{\pi i/3}(1+e^{\pi
i/3})}{2(1+e^{2\pi i/3})}& b_{5}^{23}=b_{23}^{5}:=\frac{-e^{\pi
i/3}(1+e^{\pi i/3})^2}{4(1+e^{4\pi i/3})}\\b_{5}^{24}=b_{24}^{5}
:=\frac{-t_1e^{\pi i/3}(1+e^{\pi i/3})(1+e^{2\pi i/3})}{8}&
 b_{5}^{25}=b_{25}^{5}:=\frac{-t_2e^{2\pi i/3}
 (1+e^{\pi i/3})}{4}\\
b_{5}^{26}=b_{26}^{5} :=\frac{-t_2e^{2\pi i/3}(1+e^{\pi
i/3})(1+e^{2\pi i/3})}{8}&b_{5}^{27}=b_{27}^{5}
:=\frac{-t_1t_2e^{2\pi i/3}(1+e^{\pi i/3})(1+e^{4\pi i/3})}{8}\\
b_{6}^{6} :=\frac{t_1e^{2\pi i/3}(1+e^{2\pi
i/3})}{2}&b_{6}^{7}=b_{7}^{6} :=\frac{(1+e^{2\pi i/3}+e^{4\pi
i/3}+e^{2\pi i})t_2}{4}\\
b_{6}^{8}=b_{8}^{6} :=\frac{e^{\pi
i/3}(1+e^{2\pi i/3})^2(1+e^{\pi i})t_1t_2}{8}& b_{6}^{9}=b_{9}^{6}
:=\frac{(e^{\pi i/3}+e^{\pi i})^2(1+e^{4\pi i/3})t_1t_2}{8}\\
b_{6}^{10}=b_{10}^{6}:=\frac{-1}{2}&
b_{6}^{11}=b_{11}^{6}:=\frac{-t_1e^{\pi i/3}(1+e^{2\pi i/3})}{4}\\
b_{6}^{12}=b_{12}^{6}:=\frac{-t_1e^{2\pi i/3}(1+e^{2\pi
i/3})}{2(1+e^{\pi i/3})}& b_{6}^{13}=b_{13}^{6}:=\frac{-(1+e^{2\pi
i/3})}{2(1+e^{4\pi
i/3})}\\
b_{6}^{14}=b_{14}^{6}:=\frac{-t_1e^{\pi i/3}(1+e^{\pi
i/3})(1+e^{4\pi i/3})}{8}&
b_{6}^{15}=b_{15}^{6}:=\frac{-t_1e^{2\pi i/3}(1+e^{2\pi i/3})}
{4}\\
b_{6}^{16}=b_{16}^{6}:=\frac{-t_2(1+e^{2\pi
i/3})}{4}&b_{6}^{17}=b_{17}^{6}:=\frac{-e^{\pi i/3}(1+e^{2\pi
i/3})^2t_1t_2}{8}\\
b_{6}^{18}=b_{18}^{6}:=\frac{-t_1t_2e^{2\pi i/3}(1+e^{2\pi
i/3})(1+e^{4\pi i/3})}{4}& b_{6}^{19}=b_{19}^{6}:=1\\
b_{6}^{20}=b_{20}^{6}: =\frac{-t_1(1+e^{2\pi
i/3})}{4}&b_{6}^{21}=b_{21}^{6}
:=\frac{-t_1(1+e^{2\pi i/3})}{2(1+e^{\pi i/3})}\\
b_{6}^{22}=b_{22}^{6} :=\frac{-e^{2\pi i/3}(1+e^{2\pi
i/3})}{2(1+e^{4\pi i/3})}& b_{6}^{23}=b_{23}^{6} :=\frac{-e^{2\pi
i/3}(1+e^{\pi i/3})(1+e^{2\pi i/3})t_1}{8}\\
 b_{6}^{24}=b_{24}^{6}
:=\frac{-e^{2\pi i/3}(1+e^{2\pi i/3})t_1}{4}&
b_{6}^{25}=b_{25}^{6}:=\frac{-e^{4\pi i/3}(1+e^{2\pi
i/3})t_2}{4}\\
b_{6}^{26}=b_{26}^{6} :=\frac{-e^{4\pi
i/3}(1+e^{2\pi i/3})^2t_1t_2}{8}&b_{6}^{27}=b_{27}^{6}
:=\frac{-e^{4\pi
i/3}(1+e^{2\pi i/3})(1+e^{4\pi i/3})t_1t_2}{8}\\
b_{7}^{7} :=t_2&b_{7}^{8}=b_{8}^{7} :=\frac{(1+e^{2\pi
i/3})^2t_2}{2(1+e^{\pi i/3})}\\ b_{7}^{9}=b_{9}^{7}
:=\frac{(1+e^{4\pi i/3})^2t_2}{2(1+e^{2\pi i/3})}&
b_{7}^{10}=b_{10}^{7}:=\frac{-1}{2}\\
b_{7}^{11}=b_{11}^{7}:=\frac{-e^{2\pi i/3}}{(1+e^{2\pi i/3})}&
b_{7}^{12}=b_{12}^{7}:=\frac{-e^{4\pi i/3}}{(1+e^{4\pi i/3})}\\
b_{7}^{13}=b_{13}^{7}:=\frac{-t_2}{2}&
b_{7}^{14}=b_{14}^{7}:=\frac{-t_2e^{2\pi i/3}(1+e^{\pi i/3})}{4}\\
b_{7}^{15}=b_{15}^{7}:=\frac{-t_2e^{4\pi i/3}(1+e^{2\pi i/3})}{4}
& b_{7}^{16}=b_{16}^{7}:=\frac{-t_2}{2}\\
b_{7}^{17}=b_{17}^{7} :=\frac{-e^{2\pi i/3}(1+e^{2\pi
i/3})t_2}{2(1+e^{\pi i/3})}&
b_{7}^{18}=b_{18}^{7}:=\frac{-e^{4\pi i/3}(1+e^{4\pi i/3})t_2}{2(1+e^{2\pi i/3})}\\
b_{7}^{19}=b_{19}^{7}:=1&b_{7}^{20}=b_{20}^{7}:
=\frac{1}{(1+e^{2\pi i/3})}\\b_{7}^{21}=b_{21}^{7}
:=\frac{1}{(1+e^{4\pi i/3})}& b_{7}^{22}=b_{22}^{7}
:=\frac{-t_2}{2}\\
 b_{7}^{23}=b_{23}^{7} :=\frac{-(1+e^{\pi
i/3}t_2)}{4}& b_{7}^{24}=b_{24}^{7}
:=\frac{-(1+e^{2\pi i/3}t_2)}{4}\\
b_{7}^{25}=b_{25}^{7}:=\frac{-t_2}{2}&b_{7}^{26}=b_{26}^{7}
:=\frac{-(1+e^{2\pi i/3})t_2}{2(1+e^{\pi
i/3})}\\b_{7}^{27}=b_{27}^{7} :=\frac{-(1+e^{4\pi
i/3})t_2}{2(1+e^{2\pi i/3})}& b_{8}^{8}:=\frac{e^{2\pi
i/3}(1+e^{2\pi i/3})t_2}{2}\\ b_{8}^{9}=b_{9}^{8} :=\frac{(e^{\pi
i/3}+e^{\pi i})^2(1+e^{4\pi i/3})t_1t_2}{8}&
b_{8}^{10}=b_{10}^{8}:=\frac{-1}{2}\\
b_{8}^{11}=b_{11}^{8}:=\frac{-e^{2\pi i/3}(1+e^{2\pi
i/3})}{2(1+e^{4\pi i/3})}& b_{8}^{12}=b_{12}^{8} :=\frac{-e^{4\pi
i/3}(1+e^{2\pi i/3})t_1}{4}\\
b_{8}^{13}=b_{13}^{8}:=\frac{-t_2(1+e^{2\pi i/3})}{4}&
b_{8}^{14}=b_{14}^{8} :=\frac{-t_2e^{2\pi i/3}(1+e^{\pi
i/3})(1+e^{2\pi i/3})}{8}\\
b_{8}^{15}=b_{15}^{8}:=\frac{-t_1t_2e^{4\pi i/3}(1+e^{2\pi
i/3})^2}{8} & b_{8}^{16}=b_{16}^{8}
:=\frac{-t_2(1+e^{2\pi i/3})}{2(1+e^{\pi i/3})}\\
b_{8}^{17}=b_{17}^{8} :=\frac{-t_2e^{2\pi i/3}(1+e^{2\pi
i/3})}{4}& b_{8}^{18}=b_{18}^{8} :=\frac{-e^{4\pi i/3}(1+e^{2\pi
i/3})(1+e^{4\pi i/3})t_1t_2}{8}\\
b_{8}^{19}=b_{19}^{8}:=1&b_{8}^{20}=b_{20}^{8}: =\frac{1+e^{2\pi
i/3}}{2(1+e^{4\pi i/3})}
\end{array}$$

$$\begin{array}{ll}
b_{8}^{21}=b_{21}^{8}
:=\frac{-t_1(1+e^{2\pi i/3})}{4}& b_{8}^{22}=b_{22}^{8}
:=\frac{-t_2e^{\pi i/3}(1+e^{2\pi i/3})}{4}\\
b_{8}^{23}=b_{23}^{8} :=\frac{-t_2e^{\pi i/3}(1+e^{\pi
i/3})(e^{2\pi i/3})t_2}{8}& b_{8}^{24}=b_{24}^{8} :=\frac{-e^{\pi
i/3}(1+e^{2\pi i/3})^2t_1t_2}{8}\\
b_{8}^{25}=b_{25}^{8}:=\frac{-t_2e^{2\pi i/3}(1+e^{2\pi
i/3})}{2(1+e^{\pi i/3})}&b_{8}^{26}=b_{26}^{8} :=\frac{-e^{2\pi
i/3}(1+e^{2\pi i/3})t_2}{4}\\b_{8}^{27}=b_{27}^{8}
:=\frac{-e^{2\pi i/3}(1+e^{2\pi i/3})(1+e^{4\pi
i/3})t_1t_2}{8}&b_{9}^9:=\frac{e^{4\pi i/3}(1+e^{4\pi
i/3})^2t_1t_2}{2(1+e^{\pi i/3})}\\
b_{9}^{10}=b_{10}^{9}:=\frac{-1}{2}& b_{9}^{11}=b_{11}^{9}
:=\frac{-t_1e^{2\pi i/3}(1+e^{4\pi i/3})}{4}\\
b_{9}^{12}=b_{12}^{9} :=\frac{-e^{4\pi i/3}(1+e^{4\pi i/3})t_1}
{2(1+e^{2\pi i/3})}& b_{9}^{13}=b_{13}^{9}:=\frac{-t_2(1+e^{4\pi
i/3}))}{4}\\ b_{9}^{14}=b_{14}^{9} :=\frac{-t_1t_2e^{2\pi
i/3}(1+e^{\pi i/3})(1+e^{4\pi i/3}))}{8}& b_{9}^{15}=b_{15}^{9}
:=\frac{-t_1t_2e^{4\pi i/3}(1+e^{2\pi i/3})(1+e^{4\pi i/3}))}{8}
\\b_{9}^{16}=b_{16}^{9}
:=\frac{-t_2(1+e^{4\pi i/3}))}{2(1+e^{2\pi i/3}))}&
b_{9}^{17}=b_{17}^{9} :=\frac{-e^{2\pi i/3}(1+e^{2\pi
i/3})(1+e^{4\pi i/3})t_1t_2}{8}\\ b_{9}^{18}=b_{18}^{9}
:=\frac{-e^{4\pi i/3}(1+e^{4\pi i/3})^2t_1t_2}{4(1+e^{\pi i/3})}&
b_{9}^{19}=b_{19}^{9}:=1\\ b_{9}^{20}=b_{20}^{9}:
=\frac{-t_1(1+e^{4\pi i/3})}{4}&b_{9}^{21}=b_{21}^{9}
:=\frac{(1+e^{4\pi i/3})t_1}{2(1+e^{2\pi i/3})}\\
b_{9}^{22}=b_{22}^{9} :=\frac{-t_2e^{2\pi i/3}(1+e^{4\pi
i/3})}{4}& b_{9}^{23}=b_{23}^{9} :=\frac{-e^{2\pi i/3}(1+e^{\pi
i/3})(1+e^{4\pi i/3})t_1t_2)}{8}\\ b_{9}^{24}=b_{24}^{9}
:=\frac{-e^{2\pi i/3}(1+e^{2\pi i/3})(1+e^{4\pi i/3})t_1t_2)}{8}&
b_{9}^{25}=b_{25}^{9}:=\frac{-t_2e^{4\pi i/3}(1+e^{4\pi
i/3})}{2(1+e^{2\pi i/3})}\\b_{9}^{26}=b_{26}^{9} :=\frac{-e^{4\pi
i/3}(1+e^{2\pi i/3})(1+e^{4\pi
i/3})t_1t_2}{8}&b_{9}^{27}=b_{27}^{9} :=\frac{-e^{4\pi
i/3}(1+e^{4\pi i/3})^2t_1t_2}{4(1+e^{\pi i/3})}
\\b_{10}^{10}:=\frac{-t_3}{2}& b_{10}^{11}=b_{11}^{10}
:=t_3\\
b_{10}^{12}=b_{12}^{10} :=t_3&
b_{10}^{13}=b_{13}^{10}:=t_3\\
b_{10}^{14}=b_{14}^{10} :=t_3& b_{10}^{15}=b_{15}^{10} :=t_3
\\b_{10}^{16}=b_{16}^{10}
:=t_3& b_{10}^{17}=b_{17}^{10} :=t_3\\
b_{10}^{18}=b_{18}^{10}:= t_3&
b_{10}^{19}=b_{19}^{10}:=1\\
b_{10}^{20}=b_{20}^{10}: =\frac{-1}{2}&
b_{10}^{21}=b_{21}^{10}
:=\frac{-1}{2}\\
b_{10}^{22}=b_{22}^{10} :=\frac{-1}{2}&
b_{10}^{23}=b_{23}^{10}
:=\frac{-1}{2}\\
b_{10}^{24}=b_{24}^{10}:= \frac{-1}{2}&
b_{10}^{25}=b_{25}^{10}:=\frac{-1}{2}\\b_{10}^{26}=b_{26}^{10}
:=\frac{-1}{2}&b_{10}^{27}=b_{27}^{10} :=\frac{-1}{2}\\
b_{11}^{11} :=t_3& b_{11}^{12}=b_{12}^{11} :=\frac{-t_1t_3}{2}\\
b_{11}^{13}=b_{13}^{11}:=t_3& b_{11}^{14}=b_{14}^{11}
:=\frac{t_3(1+e^{\pi i/3})^2}{2(1+e^{2\pi i/3})}\\
b_{11}^{15}=b_{15}^{11} := \frac{t_1t_2(1+e^{2\pi i/3}+e^{2\pi
i/3}+e^{\pi i}))}{4} &b_{11}^{16}=b_{16}^{11}
:=t_3\\
b_{11}^{17}=b_{17}^{11} :=\frac{(1+e^{2\pi
i/3})^2t_3)}{2(1+e^{4\pi i/3})}& b_{11}^{18}=b_{18}^{11}
:=\frac{(1+e^{2\pi i/3}+e^{4\pi i/3}+e^{2\pi i})t_1t_2}{4}\\
b_{11}^{19}=b_{19}^{11}:=\frac{-1}{2}&
b_{11}^{20}=b_{20}^{11}:=\frac{-1}{2}\\
b_{11}^{21}=b_{21}^{11} :=t_1& b_{11}^{22}=b_{22}^{11}
:=\frac{1}{(1+e^{\pi i/3})}\\
b_{11}^{23}=b_{23}^{11} :=\frac{-(1+e^{\pi i/3})}
{2(1+e^{2\pi
i/3})}& b_{11}^{24}=b_{24}^{11}
:=\frac{-t_1(1+e^{2\pi i/3})}{4}\\
b_{11}^{25}=b_{25}^{11}:=\frac{1}{1+e^{2\pi
i/3}}&b_{11}^{26}=b_{26}^{11} :=\frac{-(1+e^{2\pi
i/3})}{2(1+e^{4\pi i/3})}\\b_{11}^{27}=b_{27}^{11}
:=\frac{-t_1(1+e^{4\pi i/3})}{4}& b_{12}^{12}:=t_1t_3\\
b_{12}^{13}=b_{13}^{12}:=t_3& b_{12}^{14}=b_{14}^{12}
:=\frac{t_1t_3(1+e^{\pi i/3}+e^{2\pi i/3}+e^{\pi i})}{4}\\
b_{12}^{15}=b_{15}^{12} :=\frac{t_1t_3(1+e^{2\pi
i/3})^2}{2(1+e^{\pi i/3})} &b_{12}^{16}=b_{16}^{12}
:=t_3\\
b_{12}^{17}=b_{17}^{12}:= \frac{t_1t_3(1+e^{2\pi i/3}+e^{4\pi
i/3}+e^{2\pi i})}{4}& b_{12}^{18}=b_{18}^{12}
:=\frac{(1+e^{4\pi i/3})^2t_1t_3}{2(1+e^{2\pi i/3})}\\
b_{12}^{19}=b_{19}^{12}:=\frac{-1}{2}& b_{12}^{20}=b_{20}^{12}:
=t_1\\b_{12}^{21}=b_{21}^{12} :=\frac{-t_1}{2}&
b_{12}^{22}=b_{22}^{12} :=\frac{-1}{(1+e^{2\pi i/3})}
\end{array}$$
$$
\begin{array}{ll}
b_{12}^{23}=b_{23}^{12} :=\frac{-t_1(1+e^{\pi i/3})}{4}&
b_{12}^{24}=b_{24}^{12} :=\frac{-(1+e^{2\pi i/3})t_1}{2(1+e^{\pi
i/3})}\\ b_{12}^{25}=b_{25}^{12} :=\frac{1}{1+e^{4\pi
i/3})}&b_{12}^{26}=b_{26}^{12} :=\frac{-t_1(1+e^{2\pi
i/3}}{4}\\b_{12}^{27}=b_{27}^{12} :=\frac{-t_1(1+e^{4\pi
i/3)}}{2(1+e^{2\pi i/3})}& b_{13}^{13}:=t_3\\
b_{13}^{14}=b_{14}^{13} :=\frac{t_3(1+e^{\pi i/3})^2}{2(1+e^{2\pi
i/3})}& b_{13}^{15}=b_{15}^{13} :=\frac{t_3(1+e^{2\pi
i/3})^2}{2(1+e^{4\pi i/3})}\\b_{13}^{16}=b_{16}^{13}
:=\frac{-t_2t_3}{2}& b_{13}^{17}=b_{17}^{13}
:=\frac{t_2t_3(1+e^{\pi i/3}+e^{2\pi i/3}+e^{\pi
i})}{4}\\b_{13}^{18}=b_{18}^{13} :=\frac{t_2t_3(1+e^{2\pi
i/3}+e^{4\pi i/3}+e^{2\pi i})}{4}&
b_{13}^{19}=b_{19}^{13}:=\frac{-1}{2}\\ b_{13}^{20}=b_{20}^{13}:
=\frac{e^{\pi i/3}}{1+e^{\pi i/3}}&b_{13}^{21}=b_{21}^{13}
:=\frac{e^{2\pi i/3}}{1+e^{2\pi i/3}}\\ b_{13}^{22}=b_{22}^{13}
:=\frac{-1}{2}& b_{13}^{23}=b_{23}^{13} :=\frac{-e^{\pi
i/3}(1+e^{\pi i/3})}{2(1+e^{2\pi i/3})}\\ b_{13}^{24}=b_{24}^{13}
:\frac{-e^{2\pi i/3}(1+e^{2\pi i/3})}{2(1+e^{4\pi i/3})}&
b_{13}^{25}=b_{25}^{13}:=t_2\\b_{13}^{26}=b_{26}^{13}
:=\frac{-t_2e^{\pi i/3}(1+e^{2\pi
i/3})}{4}&b_{13}^{27}=b_{27}^{13} :=\frac{-t_2e^{2\pi
i/3}(1+e^{4\pi i/3})}{4}\\
b_{14}^{14} :=\frac{t_3e^{\pi i/3}(1+e^{\pi i/3})^2}{2(1+e^{4\pi
i/3})}&b_{14}^{15}=b_{15}^{14} :=\frac{t_1t_3(1+e^{\pi
i/3})^2(e^{\pi i/3}+e^{\pi  i})}{8}
\\b_{14}^{16}=b_{16}^{14}
:=\frac{t_2t_3(1+e^{\pi i/3}+e^{2\pi i/3}+e^{\pi i})}{4}&
b_{14}^{17}=b_{17}^{14} :=\frac{t_2t_3(1+e^{\pi i/3})^2
(e^{\pi
i/3}+e^{\pi  i})}{8}\\
b_{14}^{18}=b_{18}^{14} :=\frac{-t_1t_2t_3e^{2\pi i/3}(1+e^{\pi
i/3})(1+e^{4\pi i/3})}{8}& b_{14}^{19}=b_{19}^{14}:=\frac{-1}{2}\\
b_{14}^{20}=b_{20}^{14}: =\frac{-e^{\pi i/3}(1+e^{\pi
i/3})}{2(1+e^{2\pi i/3})}&b_{14}^{21}=b_{21}^{14}
:=\frac{-t_1e^{2\pi i/3}(1+e^{\pi i/3})}{4}\\
b_{14}^{22}=b_{22}^{14} :=\frac{-(1+e^{\pi i/3})}{2(1+e^{2\pi
i/3})}& b_{14}^{23}=b_{23}^{14} :=\frac{-e^{\pi i/3} (1+e^{\pi
i/3})^2}{4(1+e^{4\pi i/3})}\\ b_{14}^{24}=b_{24}^{14}
:=\frac{-e^{2\pi i/3}(1+e^{\pi i/3})(1+e^{2\pi i/3})t_1}{8}&
b_{14}^{25}=b_{25}^{14} :=\frac{-t_2(1+e^{\pi
i/3})}{4}\\b_{14}^{26}=b_{26}^{14} :=\frac{-t_2e^{\pi i/3}
(1+e^{\pi i/3})(1+e^{2\pi i/3})}{8}&b_{14}^{27}=b_{27}^{14}
:=\frac{e^{2\pi i/3}(1+e^{\pi i/3})(1+e^{4\pi i/3})t_1t_2}{4}\\
b_{15}^{15} :=\frac{t_1t_3e^{2\pi i/3}(1++e^{2\pi i/3})}{2}
&b_{15}^{16}=b_{16}^{15}
:=\frac{t_2t_3(1+e^{2\pi i/3}+e^{4\pi i/3}+e^{2\pi i})}{4}\\
b_{15}^{17}=b_{17}^{15} :=\frac{-t_1t_2t_3e^{\pi i/3}(1+e^{2\pi
i/3})^2 (1+e^{\pi i})}{16}& b_{15}^{18}=b_{18}^{15}
:=\frac{t_1t_2t_3(e^{\pi i/3}+e^{\pi i})^2(1+e^{4\pi i/3})}{8}\\
b_{15}^{19}=b_{19}^{15}:=\frac{-1}{2}& b_{15}^{20}=b_{20}^{15}:
=\frac{-t_1e^{\pi i/3}(1+e^{2\pi
i/3})}{4}\\b_{15}^{21}=b_{21}^{15} :=\frac{-t_1e^{2\pi
i/3}(1+e^{2\pi i/3})}{2(1+e^{\pi i/3})}& b_{15}^{22}=b_{22}^{15}
:=\frac{-(1+e^{2\pi i/3})}{2(1+e^{4\pi i/3})}\\
b_{15}^{23}=b_{23}^{15} :=\frac{-t_1e^{\pi i/3} (1+e^{\pi
i/3})(1+e^{2\pi i/3})}{8}& b_{15}^{24}=b_{24}^{15}
:=\frac{-e^{2\pi i/3}(1+e^{2\pi i/3})t_1}{4}\\
b_{15}^{25}=b_{25}^{15} :=\frac{-t_2(1+e^{2\pi
i/3})}{4}&b_{15}^{26}=b_{26}^{15} :=\frac{t_1t_2e^{\pi i/3}
(1+e^{2\pi i/3})}{4}\\b_{15}^{27}=b_{27}^{15}
 :=\frac{-e^{2\pi
i/3}(1+e^{2\pi i/3})(1+e^{4\pi i/3})t_1t_2}{8}&b_{16}^{16}
:=t_2t_3\\
b_{16}^{17}=b_{17}^{16} :=\frac{t_2t_3(1+e^{2\pi i/3})^2}
{2(1+e^{\pi/3})}& b_{16}^{18}=b_{18}^{16}
:=\frac{t_2t_3(1+e^{4\pi i/3})^2}{2(1+e^{2\pi/3})}\\
b_{16}^{19}=b_{19}^{16}:=\frac{-1}{2}& b_{16}^{20}=b_{20}^{16}:
=\frac{-e^{2\pi i/3}}{1+e^{2\pi i/3}}\\b_{16}^{21}=b_{21}^{16}
:=\frac{-e^{4\pi i/3}}{1+e^{4\pi i/3}}& b_{16}^{22}=b_{22}^{16}
:=t_2\\
b_{16}^{23}=b_{23}^{16} :=\frac{-t_2e^{2\pi i/3} (1+e^{\pi
i/3})}{4}& b_{16}^{24}=b_{24}^{16} := \frac{-t_2e^{4\pi i/3}
(1+e^{2\pi
i/3})}{4}\\
b_{16}^{25}=b_{25}^{16} :=\frac{-t_2}{2}& b_{16}^{26}=b_{26}^{16}
:=\frac{-t_2e^{2\pi i/3} (1+e^{2\pi i/3})}{2(1+e^{\pi
i/3})}\\b_{16}^{27}=b_{27}^{16}
 :=\frac{-t_2e^{4\pi i/3} (1+e^{4\pi i/3})}{2(1+e^{2\pi
i/3})}&b_{17}^{17}
:=\frac{t_2t_3e^{2\pi i/3}(1+e^{2\pi i/3})}{2}\\
 b_{17}^{18}=b_{18}^{17}
:=\frac{t_1t_2t_3(e^{\pi i/3}+e^{\pi i})^2(1+e^{4\pi i/3})}{8}&
b_{17}^{19}=b_{19}^{17}:=\frac{-1}{2}\\
b_{17}^{20}=b_{20}^{17}: =\frac{-e^{2\pi i/3}(1+e^{2\pi
i/3})}{2(1+e^{4\pi i/3})}&b_{17}^{21}=b_{21}^{17}
:=\frac{-t_1e^{4\pi i/3}(1+e^{2\pi i/3})}{4}\\
b_{17}^{22}=b_{22}^{17} :=\frac{-t_2(1+e^{2\pi i/3})}{4}&
b_{17}^{23}=b_{23}^{17} :=\frac{-t_2e^{2\pi i/3} (1+e^{\pi
i/3})(1+e^{2\pi i/3})}{8}\\ b_{17}^{24}=b_{24}^{17}
:=\frac{t_1t_2e^{4\pi i/3}(1+e^{2\pi i/3})^2}{4}&
b_{17}^{25}=b_{25}^{17} :=\frac{-t_2(1+e^{2\pi i/3})}{2(1+e^{\pi
i/3})}\\b_{17}^{26}=b_{26}^{17} :=\frac{-t_2e^{2\pi i/3}
(1+e^{2\pi i/3})}{4}&b_{17}^{27}=b_{27}^{17}
 :=\frac{-e^{4\pi
i/3}(1+e^{2\pi i/3})(1+e^{4\pi i/3})t_1t_2}{8}
\end{array}
$$

$$
\begin{array}{ll}
b_{18}^{18} :=\frac{t_1t_2t_3e^{4\pi i/3}(1+e^{4\pi
i/3})^2}{2(1+e^{\pi i/3}}& b_{18}^{19}=b_{19}^{18}:=\frac{-1}{2}\\
b_{18}^{20}=b_{20}^{18}: =
\frac{-t_1e^{2\pi i/3}(1+e^{4\pi
i/3})}{4}&b_{18}^{21}=b_{21}^{18} :=\frac{-t_1e^{4\pi
i/3}(1+e^{4\pi
i/3})}{2(1+e^{2\pi i/3})}\\
b_{18}^{22}=b_{22}^{18} :=\frac{-t_2(1+e^{4\pi i/3})}{4}&
b_{18}^{23}=b_{23}^{18} :=\frac{-t_1t_2e^{2\pi i/3} (1+e^{\pi
i/3})(1+e^{4\pi i/3})}{8}\\ b_{18}^{24}=b_{24}^{18}
:=\frac{-t_1t_2e^{4\pi i/3} (1+e^{2\pi i/3}) (1+e^{4\pi i/3})}{8}&
b_{18}^{25}=b_{25}^{18} :=\frac{-t_2(1+e^{4\pi i/3})}{2(1+e^{2\pi
i/3})}\\b_{18}^{26}=b_{26}^{18} :=\frac{-t_1t_2e^{2\pi i/3}
(1+e^{2\pi i/3})(1+e^{4\pi i/3})}{8}&b_{18}^{27}=b_{27}^{18}
 :=\frac{-e^{4\pi
i/3}(1+e^{4\pi i/3})^2t_1t_2}{4(1+e^{\pi i/3})}\\
b_{19}^{19}:=\frac{-t_3^{-1}}{2}& b_{19}^{20}=b_{20}^{19}:
=t_3^{-1}\\b_{19}^{21}=b_{21}^{19} t_3^{-1}&
b_{19}^{22}=b_{22}^{19}
:=t_3^{-1}\\
b_{19}^{23}=b_{23}^{19} :=t_3^{-1}& b_{19}^{24}=b_{24}^{19}
:=t_3^{-1}\\
b_{19}^{25}=b_{25}^{19} t_3^{-1}&b_{19}^{26}=b_{26}^{19}
:=t_3^{-1}\\b_{19}^{27}=b_{27}^{19}
 :=t_3^{-1}&b_{20}^{20}:=t_3^{-1}\\
 b_{20}^{21}=b_{21}^{20} :=\frac{-t_1t_3^{-1}}{2}&b_{20}^{22}=b_{22}^{20}
:=t_3^{-1}\\
b_{20}^{23}=b_{23}^{20} :=\frac{t_3^{-1}(1+e^{\pi i/3})^2}
{2(1+e^{2\pi i/3})}& b_{20}^{24}=b_{24}^{20}
:=\frac{t_1t_3^{-1}(1+e^{\pi i/3}+e^{2\pi i/3}+e^{\pi i})}{4}\\
b_{20}^{25}=b_{25}^{20} :=t_3^{-1}&b_{20}^{26}=b_{26}^{20}
:=\frac{t_3^{-1}(1+e^{2\pi i/3})^2 }{2(1+e^{4\pi
i/3})}\\b_{20}^{27}=b_{27}^{20}
 :=\frac{t_1t_3^{-1}(1+e^{2\pi i/3}+e^{4\pi i/3}+e^{2\pi i})}{4}
 &b_{21}^{21}:=t_1t_3^{-1}
 \\
 b_{21}^{22}=b_{22}^{21}
:=t_3^{-1}& b_{21}^{23}=b_{23}^{21} :=\frac{t_1t_3^{-1}(1+e^{\pi
i/3} +e^{2\pi i/3}+e^{\pi i})}{4}\\ b_{21}^{24}=b_{24}^{21}
:=\frac{-t_1t_3^{-1}(1+e^{2\pi i/3})^2}{2(1+e^{\pi i/3})}&
b_{21}^{25}=b_{25}^{21} :=t_3^{-1}\\b_{21}^{26}=b_{26}^{21}
:=\frac{t_1t_3^{-1} (1+e^{2\pi i/3}+e^{4\pi i/3}+e^{2\pi
i})}{4}&b_{21}^{27}=b_{27}^{21}
 :=\frac{t_1t_3^{-1}(1+e^{4\pi i/3})^2}{2(1+e^{2\pi i/3})}\\
b_{22}^{22} :=t_3^{-1}& b_{22}^{23}=b_{23}^{22}
:=\frac{t_3^{-1}(1+e^{\pi i/3})^2}{2(1+e^{2\pi i/3})}\\
b_{22}^{24}=b_{24}^{22} \frac{t_3^{-1}(1+e^{2\pi
i/3})^2}{2(1+e^{4\pi i/3})}& b_{22}^{25}=b_{25}^{22}
:=\frac{t_2t_3^{-1}}{2}\\b_{22}^{26}=b_{26}^{22}
:=\frac{t_2t_3^{-1} (1+e^{\pi i/3}+e^{2\pi i/3}+e^{\pi
i})}{4}&b_{22}^{27}=b_{27}^{22}
 :=\frac{t_2t_3^{-1} (1+e^{2\pi i/3}+e^{4\pi i/3}+e^{2\pi
i})}{4}\\
  b_{23}^{23} :=\frac{t_3^{-1}e^{\pi i/3}(1+e^{\pi
i/3} )^2}{2(1+e^{4\pi i/3})}& b_{23}^{24}=b_{24}^{23}
:=\frac{t_1t_3^{-1}(1+e^{2\pi i/3})^2(e^{\pi i/3}+e^{\pi i})}{8}\\
b_{23}^{25}=b_{25}^{23} :=\frac{t_2t_3^{-1}(1+e^{\pi i/3}+e^{2\pi
i/3}+e^{\pi i})}{4}&b_{23}^{26}=b_{26}^{23} :=\frac{t_2t_3^{-1}
(1+e^{\pi i/3})^2(e^{\pi i/3}+e^{\pi
i})}{8}\\b_{23}^{27}=b_{27}^{23}
 :=\frac{-t_1t_2t_3^{-1}e^{2\pi i/3}(1+e^{\pi i/3})(1+e^{4\pi i/3})}{8}
 &b_{24}^{24}:=\frac{t_1t_3^{-1}e^{2\pi i/3}(1+e^{2\pi i/3})}{2}\\
  b_{24}^{25}=b_{25}^{24}
  :=\frac{t_2t_3^{-1}(1+e^{2\pi i/3}+e^{4\pi
i/3}+e^{2\pi i})}{4}&b_{24}^{26}=b_{26}^{24}
:=\frac{-t_1t_2t_3^{-1} e^{\pi i/3}(1+e^{2\pi i/3})^2(1+e^{\pi
i})}{16}\\b_{24}^{27}=b_{27}^{24}
 :=\frac{-t_1t_2t_3^{-1}(e^{\pi i/3}+e^{\pi i})^2
 (1+e^{4\pi i/3})}{8}
 &b_{25}^{25}:=t_2t_3^{-1}\\
 b_{25}^{26}=b_{26}^{25}:=\frac{t_2t_3^{-1} (1+e^{2\pi i/3})^2}{2(1+e^{\pi i/3})}
 & b_{25}^{27}=b_{27}^{25}:=\frac{t_2t_3^{-1}
 (1+e^{4\pi i/3})^2}{2(1+e^{2\pi i/3})}\\
 b_{26}^{26}:=\frac{t_2t_3^{-1}e^{2\pi i/3}(1+e^{2\pi i/3})}{2}&
 b_{26}^{27}=b_{27}^{26}:=\frac{t_1t_2t_3^{-1}(e^{\pi i/3}+e^{\pi i})^2(1+e^{4\pi
 i/3})}{8}\\
 b_{27}^{27}:=
 \frac{t_1t_2t_3^{-1}e^{4\pi i/3}(1+e^{4\pi i/3})^2}{2(1+e^{\pi i/3})}
\end{array}
$$

\begin{center}
{\small\begin{table}\caption{the table of the products (for $\fcj_p$)}
\begin{tabular}{|l||p{.8cm}p{.8cm}p{.8cm}p{.8cm}p{.8cm}p{.8cm}p{.8cm}p{.8cm}p{.8cm}|}
\hline
 &$w_1$&$w_2$&$w_3$&$w_4$&$w_5$&$w_6$&$w_7$&$w_8$&$w_9$\\
\hline
\vspace{-1mm}&&&&&&&&&\\
$w_1$ &$b_{1}^{1}$&$b_{1}^{2}w_2$&$b_{1}^{3}w_3$&$b_{1}^{4}w_4$&$b_{1}^{5}w_5$&$b_{1}^{6}w_6$&$b_{1}^{7}w_7$&$b_{1}^{8}w_8$&$b_{1}^{9}w_9$\\
\vspace{-1mm}&&&&&&&&&\\
$w_2$ &$b_{2}^{1}w_2$&$b_{2}^{2}w_3$&$b_{2}^{3}$&$b_{2}^{4}w_5$&$b_{2}^{5}w_6$&$b_{2}^{6}w_4$&$b_{2}^{7}w_8$&$b_{2}^{8}w_9$&$b_{2}^{9}w_7$\\
\vspace{-1mm}&&&&&&&&&\\
$w_3$ &$b_{3}^{1}w_3$&$b_{3}^{2}$&$b_{3}^{3}w_2$&$b_{3}^{4}w_6$&$b_{3}^{5}w_4$&$b_{3}^{6}w_5$&$b_{3}^{7}w_9$&$b_{3}^{8}w_7$&$b_{3}^{9}w_8$\\
\vspace{-1mm}&&&&&&&&&\\
$w_4$ &$b_{4}^{1}w_4$&$b_{4}^{2}w_5$&$b_{4}^{3}w_6$&$b_{4}^{4}w_7$&$b_{4}^{5}w_8$&$b_{4}^{6}w_9$&$b_{4}^{7}$&$b_{4}^{8}w_2$&$b_{4}^{9}w_3$\\
\vspace{-1mm}&&&&&&&&&\\
$w_5$ &$b_{5}^{1}w_5$&$b_{5}^{2}w_6$&$b_{5}^{3}w_4$&$b_{5}^{4}w_8$&$b_{5}^{5}w_9$&$b_{5}^{6}w_7$&$b_{5}^{7}w_2$&$b_{5}^{8}w_3$&$b_{5}^{9}$\\
\vspace{-1mm}&&&&&&&&&\\
$w_6$ &$b_{6}^{1}w_6$&$b_{6}^{2}w_4$&$b_{6}^{3}w_5$&$b_{6}^{4}w_9$&$b_{6}^{5}w_7$&$b_{6}^{6}w_8$&$b_{6}^{7}w_3$&$b_{6}^{8}$&$b_{6}^{9}w_2$\\
\vspace{-1mm}&&&&&&&&&\\
$w_7$ &$b_{7}^{1}w_7$&$b_{7}^{2}w_8$&$b_{7}^{3}w_9$&$b_{7}^{4}$&$b_{7}^{5}w_2$&$b_{7}^{6}w_3$&$b_{7}^{7}w_4$&$b_{7}^{8}w_5$&$b_{7}^{9}w_6$\\
\vspace{-1mm}&&&&&&&&&\\
$w_8$ &$b_{8}^{1}w_8$&$b_{8}^{2}w_9$&$b_{8}^{3}w_7$&$b_{8}^{4}w_2$&$b_{8}^{5}w_3$&$b_{8}^{6}$&$b_{8}^{7}w_5$&$b_{8}^{8}w_6$&$b_{8}^{9}w_4$\\
\vspace{-1mm}&&&&&&&&&\\
$w_9$ &$b_{9}^{1}w_9$&$b_{9}^{2}w_7$&$b_{9}^{3}w_8$&$b_{9}^{4}w_3$&$b_{9}^{5}$&$b_{9}^{6}w_2$&$b_{9}^{7}w_6$&$b_{9}^{8}w_4$&$b_{9}^{9}w_5$\\
\vspace{-1mm}&&&&&&&&&\\
$w_{10}$ &$b_{10}^{1}w_{10}$&$b_{10}^{2}w_2$&$b_{10}^{3}w_3$&$b_{10}^{4}w_4$&$b_{10}^{5}w_5$&$b_{10}^{6}w_6$&$b_{10}^{7}w_7$&$b_{10}^{8}w_8$&$b_{10}^{9}w_9$\\
\vspace{-1mm}&&&&&&&&&\\
$w_{11}$ &$b_{11}^{1}w_{11}$&$b_{11}^{2}w_2$&$b_{11}^{3}w_3$&$b_{11}^{4}w_4$&$b_{11}^{5}w_5$&$b_{11}^{6}w_6$&$b_{11}^{7}w_7$&$b_{11}^{8}w_8$&$b_{11}^{9}w_9$\\
\vspace{-1mm}&&&&&&&&&\\
$w_{12}$ &$b_{12}^{1}w_{12}$&$b_{12}^{2}w_2$&$b_{12}^{3}w_3$&$b_{12}^{4}w_4$&$b_{12}^{5}w_5$&$b_{12}^{6}w_6$&$b_{12}^{7}w_7$&$b_{12}^{8}w_8$&$b_{12}^{9}w_9$\\
\vspace{-1mm}&&&&&&&&&\\
$w_{13}$ &$b_{13}^{1}w_{13}$&$b_{13}^{2}w_2$&$b_{13}^{3}w_3$&$b_{13}^{4}w_4$&$b_{13}^{5}w_5$&$b_{13}^{6}w_6$&$b_{13}^{7}w_7$&$b_{13}^{8}w_8$&$b_{13}^{9}w_9$\\
\vspace{-1mm}&&&&&&&&&\\
$w_{14}$ &$b_{14}^{1}w_{14}$&$b_{14}^{2}w_2$&$b_{14}^{3}w_3$&$b_{14}^{4}w_4$&$b_{14}^{5}w_5$&$b_{14}^{6}w_6$&$b_{14}^{7}w_7$&$b_{14}^{8}w_8$&$b_{14}^{9}w_9$\\
\vspace{-1mm}&&&&&&&&&\\
$w_{15}$ &$b_{15}^{1}w_{15}$&$b_{15}^{2}w_2$&$b_{15}^{3}w_3$&$b_{15}^{4}w_4$&$b_{15}^{5}w_5$&$b_{15}^{6}w_6$&$b_{15}^{7}w_7$&$b_{15}^{8}w_8$&$b_{15}^{9}w_9$\\
\vspace{-1mm}&&&&&&&&&\\
$w_{16}$ &$b_{16}^{1}w_{16}$&$b_{16}^{2}w_2$&$b_{16}^{3}w_3$&$b_{16}^{4}w_4$&$b_{16}^{5}w_5$&$b_{16}^{6}w_6$&$b_{16}^{7}w_7$&$b_{16}^{8}w_8$&$b_{16}^{9}w_9$\\
\vspace{-1mm}&&&&&&&&&\\
$w_{17}$ &$b_{17}^{1}w_{17}$&$b_{17}^{2}w_2$&$b_{17}^{3}w_3$&$b_{17}^{4}w_4$&$b_{17}^{5}w_5$&$b_{17}^{6}w_6$&$b_{17}^{7}w_7$&$b_{17}^{8}w_8$&$b_{17}^{9}w_9$\\
\vspace{-1mm}&&&&&&&&&\\
$w_{18}$ &$b_{18}^{1}w_{18}$&$b_{18}^{2}w_2$&$b_{18}^{3}w_3$&$b_{18}^{4}w_4$&$b_{18}^{5}w_5$&$b_{18}^{6}w_6$&$b_{18}^{7}w_7$&$b_{18}^{8}w_8$&$b_{18}^{9}w_9$\\
\vspace{-1mm}&&&&&&&&&\\
$w_{19}$ &$b_{19}^{1}w_{19}$&$b_{19}^{2}w_2$&$b_{19}^{3}w_3$&$b_{19}^{4}w_4$&$b_{19}^{5}w_5$&$b_{19}^{6}w_6$&$b_{19}^{7}w_7$&$b_{19}^{8}w_8$&$b_{19}^{9}w_9$\\
\vspace{-1mm}&&&&&&&&&\\
$w_{20}$ &$b_{20}^{1}w_{20}$&$b_{20}^{2}w_2$&$b_{20}^{3}w_3$&$b_{20}^{4}w_4$&$b_{20}^{5}w_5$&$b_{20}^{6}w_6$&$b_{20}^{7}w_7$&$b_{20}^{8}w_8$&$b_{20}^{9}w_9$\\
\vspace{-1mm}&&&&&&&&&\\
$w_{21}$ &$b_{21}^{1}w_{21}$&$b_{21}^{2}w_2$&$b_{21}^{3}w_3$&$b_{21}^{4}w_4$&$b_{21}^{5}w_5$&$b_{21}^{6}w_6$&$b_{21}^{7}w_7$&$b_{21}^{8}w_8$&$b_{21}^{9}w_9$\\
\vspace{-1mm}&&&&&&&&&\\
$w_{22}$ &$b_{22}^{1}w_{22}$&$b_{22}^{2}w_2$&$b_{22}^{3}w_3$&$b_{22}^{4}w_4$&$b_{22}^{5}w_5$&$b_{22}^{6}w_6$&$b_{22}^{7}w_7$&$b_{22}^{8}w_8$&$b_{22}^{9}w_9$\\
\vspace{-1mm}&&&&&&&&&\\
$w_{23}$ &$b_{23}^{1}w_{23}$&$b_{23}^{2}w_2$&$b_{23}^{3}w_3$&$b_{23}^{4}w_4$&$b_{23}^{5}w_5$&$b_{23}^{6}w_6$&$b_{23}^{7}w_7$&$b_{23}^{8}w_8$&$b_{23}^{9}w_9$\\
\vspace{-1mm}&&&&&&&&&\\
$w_{24}$ &$b_{24}^{1}w_{24}$&$b_{24}^{2}w_2$&$b_{24}^{3}w_3$&$b_{24}^{4}w_4$&$b_{24}^{5}w_5$&$b_{24}^{6}w_6$&$b_{24}^{7}w_7$&$b_{24}^{8}w_8$&$b_{24}^{9}w_9$\\
\vspace{-1mm}&&&&&&&&&\\
$w_{25}$ &$b_{25}^{1}w_{25}$&$b_{25}^{2}w_2$&$b_{25}^{3}w_3$&$b_{25}^{4}w_4$&$b_{25}^{5}w_5$&$b_{25}^{6}w_6$&$b_{25}^{7}w_7$&$b_{25}^{8}w_8$&$b_{25}^{9}w_9$\\
\vspace{-1mm}&&&&&&&&&\\
$w_{27}$ &$b_{26}^{1}w_{26}$&$b_{26}^{2}w_2$&$b_{26}^{3}w_3$&$b_{26}^{4}w_4$&$b_{26}^{5}w_5$&$b_{26}^{6}w_6$&$b_{26}^{7}w_7$&$b_{26}^{8}w_8$&$b_{26}^{9}w_9$\\
\vspace{-1mm}&&&&&&&&&\\
$w_{27}$ &$b_{27}^{1}w_{27}$&$b_{27}^{2}w_2$&$b_{27}^{3}w_3$&$b_{27}^{4}w_4$&$b_{27}^{5}w_5$&$b_{27}^{6}w_6$&$b_{27}^{7}w_7$&$b_{27}^{8}w_8$&$b_{27}^{9}w_9$\\
\hline
\end{tabular}
\end{table}}
\end{center}

\begin{center}

{\small\begin{table}\caption{the table of the products (for $\fcj_p$)}
\begin{tabular}{|l||p{.8cm}p{.8cm}p{.8cm}p{.8cm}p{.8cm}p{.8cm}p{.8cm}p{.8cm}p{.8cm}|}
\hline
 &$w_{10}$&$w_{11}$&$w_{12}$&$w_{13}$&$w_{14}$&$w_{15}$&$w_{16}$&$w_{17}$&$w_{18}$\\
\hline
$w_1$ &$b_{1}^{10}$&$b_{1}^{11}w_2$&$b_{1}^{12}w_2$&$b_{1}^{13}w_3$&$b_{1}^{14}w_4$&$b_{1}^{15}w_5$&$b_{1}^{16}w_6$&$b_{1}^{17}w_7$&$b_{1}^{18}w_8$\\
\vspace{-1mm}&&&&&&&&&\\
$w_2$ &$b_{2}^{10}w_2$&$b_{2}^{11}w_3$&$b_{2}^{12}$&$b_{2}^{13}w_5$&$b_{2}^{14}w_6$&$b_{2}^{15}w_4$&$b_{2}^{16}w_8$&$b_{2}^{17}w_9$&$b_{2}^{18}w_7$\\
\vspace{-1mm}&&&&&&&&&\\
$w_3$ &$b_{3}^{10}w_3$&$b_{3}^{11}$&$b_{3}^{12}w_2$&$b_{3}^{13}w_6$&$b_{3}^{14}w_4$&$b_{3}^{15}w_5$&$b_{3}^{16}w_9$&$b_{3}^{17}w_7$&$b_{3}^{18}w_8$\\
\vspace{-1mm}&&&&&&&&&\\
$w_4$ &$b_{4}^{10}$&$b_{4}^{11}w_4$&$b_{4}^{12}w_5$&$b_{4}^{13}w_6$&$b_{4}^{14}w_7$&$b_{4}^{15}w_8$&$b_{4}^{16}w_9$&$b_{4}^{17}w_2$&$b_{4}^{18}w_3$\\
\vspace{-1mm}&&&&&&&&&\\
$w_5$ &$b_{5}^{10}w_5$&$b_{5}^{11}w_5$&$b_{5}^{12}w_6$&$b_{5}^{13}w_4$&$b_{5}^{14}w_8$&$b_{5}^{15}w_9$&$b_{5}^{16}w_7$&$b_{5}^{17}w_2$&$b_{5}^{18}w_3$\\
\vspace{-1mm}&&&&&&&&&\\
$w_6$ &$b_{6}^{10}w_2$&$b_{6}^{11}w_6$&$b_{6}^{12}w_4$&$b_{6}^{13}w_5$&$b_{6}^{14}w_9$&$b_{6}^{15}w_7$&$b_{6}^{16}w_8$&$b_{6}^{17}w_3$&$b_{6}^{18}$\\
\vspace{-1mm}&&&&&&&&&\\
$w_7$ &$b_{7}^{10}w_6$&$b_{7}^{11}w_7$&$b_{7}^{12}w_8$&$b_{7}^{13}w_9$&$b_{7}^{14}$&$b_{7}^{15}w_2$&$b_{7}^{16}w_3$&$b_{7}^{17}w_4$&$b_{7}^{18}w_5$\\
\vspace{-1mm}&&&&&&&&&\\
$w_8$ &$b_{8}^{10}w_4$&$b_{8}^{11}w_8$&$b_{8}^{12}w_9$&$b_{8}^{13}w_7$&$b_{8}^{14}w_2$&$b_{8}^{15}w_3$&$b_{8}^{16}$&$b_{8}^{17}w_5$&$b_{8}^{18}w_6$\\
\vspace{-1mm}&&&&&&&&&\\
$w_9$ &$b_{9}^{10}w_5$&$b_{9}^{11}w_9$&$b_{9}^{12}w_7$&$b_{9}^{13}w_8$&$b_{9}^{14}w_3$&$b_{9}^{15}$&$b_{9}^{16}w_2$&$b_{9}^{17}w_6$&$b_{9}^{18}w_4$\\
\vspace{-1mm}&&&&&&&&&\\
$w_{10}$ &$b_{10}^{10}w_9$&$b_{10}^{11}w_{10}$&$b_{10}^{12}w_2$&$b_{10}^{13}w_3$&$b_{10}^{14}w_4$&$b_{10}^{15}w_5$&$b_{10}^{16}w_6$&$b_{10}^{17}w_7$&$b_{10}^{18}w_8$\\
\vspace{-1mm}&&&&&&&&&\\
$w_{11}$ &$b_{11}^{10}w_9$&$b_{11}^{11}w_{11}$&$b_{11}^{12}w_2$&$b_{11}^{13}w_3$&$b_{11}^{14}w_4$&$b_{11}^{15}w_5$&$b_{11}^{16}w_6$&$b_{11}^{17}w_7$&$b_{11}^{18}w_8$\\
\vspace{-1mm}&&&&&&&&&\\
$w_{12}$ &$b_{12}^{10}w_9$&$b_{12}^{11}w_{12}$&$b_{12}^{12}w_2$&$b_{12}^{13}w_3$&$b_{12}^{14}w_4$&$b_{12}^{15}w_5$&$b_{12}^{16}w_6$&$b_{12}^{17}w_7$&$b_{12}^{18}w_8$\\
\vspace{-1mm}&&&&&&&&&\\
$w_{13}$ &$b_{13}^{10}w_9$&$b_{13}^{11}w_{13}$&$b_{13}^{12}w_2$&$b_{13}^{13}w_3$&$b_{13}^{14}w_4$&$b_{13}^{15}w_5$&$b_{13}^{16}w_6$&$b_{13}^{17}w_7$&$b_{13}^{18}w_8$\\
\vspace{-1mm}&&&&&&&&&\\
$w_{14}$ &$b_{14}^{10}w_9$&$b_{14}^{11}w_{14}$&$b_{14}^{12}w_2$&$b_{14}^{13}w_3$&$b_{14}^{14}w_4$&$b_{14}^{15}w_5$&$b_{14}^{16}w_6$&$b_{14}^{17}w_7$&$b_{14}^{18}w_8$\\
\vspace{-1mm}&&&&&&&&&\\
$w_{15}$ &$b_{15}^{10}w_9$&$b_{15}^{11}w_{15}$&$b_{15}^{12}w_2$&$b_{15}^{13}w_3$&$b_{15}^{14}w_4$&$b_{15}^{15}w_5$&$b_{15}^{16}w_6$&$b_{15}^{17}w_7$&$b_{15}^{18}w_8$\\
\vspace{-1mm}&&&&&&&&&\\
$w_{16}$ &$b_{16}^{10}w_9$&$b_{16}^{11}w_{16}$&$b_{16}^{12}w_2$&$b_{16}^{13}w_3$&$b_{16}^{14}w_4$&$b_{16}^{15}w_5$&$b_{16}^{16}w_6$&$b_{16}^{17}w_7$&$b_{16}^{18}w_8$\\
\vspace{-1mm}&&&&&&&&&\\
$w_{17}$ &$b_{17}^{10}w_9$&$b_{17}^{11}w_{17}$&$b_{17}^{12}w_2$&$b_{17}^{13}w_3$&$b_{17}^{14}w_4$&$b_{17}^{15}w_5$&$b_{17}^{16}w_6$&$b_{17}^{17}w_7$&$b_{17}^{18}w_8$\\
\vspace{-1mm}&&&&&&&&&\\
$w_{18}$ &$b_{18}^{10}w_9$&$b_{18}^{11}w_{18}$&$b_{18}^{12}w_2$&$b_{18}^{13}w_3$&$b_{18}^{14}w_4$&$b_{18}^{15}w_5$&$b_{18}^{16}w_6$&$b_{18}^{17}w_7$&$b_{18}^{18}w_8$\\
\vspace{-1mm}&&&&&&&&&\\
$w_{19}$ &$b_{19}^{10}w_9$&$b_{19}^{11}w_{19}$&$b_{19}^{12}w_2$&$b_{19}^{13}w_3$&$b_{19}^{14}w_4$&$b_{19}^{15}w_5$&$b_{19}^{16}w_6$&$b_{19}^{17}w_7$&$b_{19}^{18}w_8$\\
\vspace{-1mm}&&&&&&&&&\\
$w_{20}$ &$b_{20}^{10}w_9$&$b_{20}^{11}w_{20}$&$b_{20}^{12}w_2$&$b_{20}^{13}w_3$&$b_{20}^{14}w_4$&$b_{20}^{15}w_5$&$b_{20}^{16}w_6$&$b_{20}^{17}w_7$&$b_{20}^{18}w_8$\\
\vspace{-1mm}&&&&&&&&&\\
$w_{21}$ &$b_{21}^{10}w_9$&$b_{21}^{11}w_{21}$&$b_{21}^{12}w_2$&$b_{21}^{13}w_3$&$b_{21}^{14}w_4$&$b_{21}^{15}w_5$&$b_{21}^{16}w_6$&$b_{21}^{17}w_7$&$b_{21}^{18}w_8$\\
\vspace{-1mm}&&&&&&&&&\\
$w_{22}$ &$b_{22}^{10}w_9$&$b_{22}^{11}w_{22}$&$b_{22}^{12}w_2$&$b_{22}^{13}w_3$&$b_{22}^{14}w_4$&$b_{22}^{15}w_5$&$b_{22}^{16}w_6$&$b_{22}^{17}w_7$&$b_{22}^{18}w_8$\\
\vspace{-1mm}&&&&&&&&&\\
$w_{23}$ &$b_{23}^{10}w_9$&$b_{23}^{11}w_{23}$&$b_{23}^{12}w_2$&$b_{23}^{13}w_3$&$b_{23}^{14}w_4$&$b_{23}^{15}w_5$&$b_{23}^{16}w_6$&$b_{23}^{17}w_7$&$b_{23}^{18}w_8$\\
\vspace{-1mm}&&&&&&&&&\\
$w_{24}$ &$b_{24}^{10}w_9$&$b_{24}^{11}w_{24}$&$b_{24}^{12}w_2$&$b_{24}^{13}w_3$&$b_{24}^{14}w_4$&$b_{24}^{15}w_5$&$b_{24}^{16}w_6$&$b_{24}^{17}w_7$&$b_{24}^{18}w_8$\\
\vspace{-1mm}&&&&&&&&&\\
$w_{25}$ &$b_{25}^{10}w_9$&$b_{25}^{11}w_{25}$&$b_{25}^{12}w_2$&$b_{25}^{13}w_3$&$b_{25}^{14}w_4$&$b_{25}^{15}w_5$&$b_{25}^{16}w_6$&$b_{25}^{17}w_7$&$b_{25}^{18}w_8$\\
\vspace{-1mm}&&&&&&&&&\\
$w_{27}$ &$b_{26}^{10}w_9$&$b_{26}^{11}w_{26}$&$b_{26}^{12}w_2$&$b_{26}^{13}w_3$&$b_{26}^{14}w_4$&$b_{26}^{15}w_5$&$b_{26}^{16}w_6$&$b_{26}^{17}w_7$&$b_{26}^{18}w_8$\\
\vspace{-1mm}&&&&&&&&&\\
$w_{27}$ &$b_{27}^{10}w_9$&$b_{27}^{11}w_{27}$&$b_{27}^{12}w_2$&$b_{27}^{13}w_3$&$b_{27}^{14}w_4$&$b_{27}^{15}w_5$&$b_{27}^{16}w_6$&$b_{27}^{17}w_7$&$b_{27}^{18}w_8$\\
\hline
\end{tabular}
\end{table}}
\end{center}

\begin{center}
{\small\begin{table}\caption{the table of the products (for $\fcj_p$)}
\begin{tabular}{|l||p{.8cm}p{.8cm}p{.8cm}p{.8cm}p{.8cm}p{.8cm}p{.8cm}p{.8cm}p{.8cm}|}
\hline
 &$w_{19}$&$w_{20}$&$w_{21}$&$w_{22}$&$w_{23}$&$w_{24}$&$w_{25}$&$w_{26}$&$w_{27}$\\
\hline
$w_1$
&$b_{1}^{19}w_{19}$&$b_{1}^{20}w_{20}$&$b_{1}^{21}w_{21}$&$b_{1}^{22}w_{22}$&$b_{1}^{23}w_{23}$&$b_{1}^{24}w_{24}$&$b_{1}^{25}w_{25}$&$b_{1}^{26}w_{26}$&$b_{1}^{27}w_{27}$\\
\vspace{-1mm}&&&&&&&&&\\
$w_2$ &$b_{2}^{19}w_{20}$&$b_{2}^{20}w_{21}$&$b_{2}^{21}w_{19}$&$b_{2}^{22}w_{23}$&$b_{2}^{23}w_{24}$&$b_{2}^{24}w_{22}$&$b_{2}^{25}w_{26}$&$b_{2}^{26}w_{27}$&$b_{2}^{27}w_{25}$\\
\vspace{-1mm}&&&&&&&&&\\
$w_3$ &$b_{3}^{19}w_{21}$&$b_{3}^{20}w_{19}$&$b_{3}^{21}w_{20}$&$b_{3}^{22}w_{24}$&$b_{3}^{23}w_{22}$&$b_{3}^{24}w_{23}$&$b_{3}^{25}w_{27}$&$b_{3}^{26}w_{25}$&$b_{3}^{27}w_{26}$\\
\vspace{-1mm}&&&&&&&&&\\
$w_4$ &$b_{4}^{19}w_{22}$&$b_{4}^{20}w_{23}$&$b_{4}^{21}w_{24}$&$b_{4}^{22}w_{25}$&$b_{4}^{23}w_{26}$&$b_{4}^{24}w_{27}$&$b_{4}^{25}w_{19}$&$b_{4}^{26}w_{20}$&$b_{4}^{27}w_{21}$\\
\vspace{-1mm}&&&&&&&&&\\
$w_5$ &$b_{5}^{19}w_{23}$&$b_{5}^{20}w_{24}$&$b_{5}^{21}w_{22}$&$b_{5}^{22}w_{26}$&$b_{5}^{23}w_{27}$&$b_{5}^{24}w_{25}$&$b_{5}^{25}w_{20}$&$b_{5}^{26}w_{21}$&$b_{5}^{27}w_{19}$\\
\vspace{-1mm}&&&&&&&&&\\
$w_6$ &$b_{6}^{19}w_{24}$&$b_{6}^{20}w_{22}$&$b_{6}^{21}w_{23}$&$b_{6}^{22}w_{27}$&$b_{6}^{23}w_{25}$&$b_{6}^{24}w_{26}$&$b_{6}^{25}w_{21}$&$b_{6}^{26}w_{19}$&$b_{6}^{27}w_{20}$\\
\vspace{-1mm}&&&&&&&&&\\
$w_7$ &$b_{7}^{19}w_{25}$&$b_{7}^{20}w_{26}$&$b_{7}^{21}w_{27}$&$b_{7}^{22}w_{19}$&$b_{7}^{23}w_{20}$&$b_{7}^{24}w_{21}$&$b_{7}^{25}w_{22}$&$b_{7}^{26}w_{23}$&$b_{7}^{27}w_{24}$\\
\vspace{-1mm}&&&&&&&&&\\
$w_8$ &$b_{8}^{19}w_{26}$&$b_{8}^{20}w_{27}$&$b_{8}^{21}w_{25}$&$b_{8}^{22}w_{20}$&$b_{8}^{23}w_{21}$&$b_{8}^{24}w_{19}$&$b_{8}^{25}w_{23}$&$b_{8}^{26}w_{24}$&$b_{8}^{27}w_{22}$\\
\vspace{-1mm}&&&&&&&&&\\
$w_9$ &$b_{9}^{19}w_{27}$&$b_{9}^{20}w_{25}$&$b_{9}^{21}w_{26}$&$b_{9}^{22}w_{21}$&$b_{9}^{23}w_{19}$&$b_{9}^{24}w_{20}$&$b_{9}^{25}w_{24}$&$b_{9}^{26}w_{22}$&$b_{9}^{27}w_{23}$\\
\vspace{-1mm}&&&&&&&&&\\
$w_{10}$ &$b_{10}^{19}$&$b_{10}^{20}w_2$&$b_{10}^{21}w_{3}$&$b_{10}^{22}w_4$&$b_{10}^{23}w_5$&$b_{10}^{24}w_6$&$b_{10}^{25}w_8$&$b_{10}^{26}w_8$&$b_{10}^{27}w_9$\\
\vspace{-1mm}&&&&&&&&&\\
$w_{11}$ &$b_{11}^{19}w_2$&$b_{11}^{20}w_3$&$b_{11}^{21}$&$b_{11}^{22}w_5$&$b_{11}^{23}w_6$&$b_{11}^{24}w_4$&$b_{11}^{25}w_8$&$b_{11}^{26}w_9$&$b_{11}^{27}w_7$\\
\vspace{-1mm}&&&&&&&&&\\
$w_{12}$ &$b_{12}^{19}w_3$&$b_{12}^{20}$&$b_{12}^{21}w_{2}$&$b_{12}^{22}w_6$&$b_{12}^{23}w_4$&$b_{12}^{24}w_5$&$b_{12}^{25}w_9$&$b_{12}^{26}w_7$&$b_{12}^{27}w_8$\\
\vspace{-1mm}&&&&&&&&&\\
$w_{13}$ &$b_{13}^{19}w_4$&$b_{13}^{20}w_5$&$b_{13}^{21}w_{6}$&$b_{13}^{22}w_7$&$b_{13}^{23}w_8$&$b_{13}^{24}w_9$&$b_{13}^{25}$&$b_{13}^{26}w_2$&$b_{13}^{27}w_3$\\
\vspace{-1mm}&&&&&&&&&\\
$w_{14}$ &$b_{14}^{19}w_5$&$b_{14}^{20}w_6$&$b_{14}^{21}w_{4}$&$b_{14}^{22}w_8$&$b_{14}^{23}w_9$&$b_{14}^{24}w_7$&$b_{14}^{25}w_2$&$b_{14}^{26}w_3$&$b_{14}^{27}$\\
\vspace{-1mm}&&&&&&&&&\\
$w_{15}$ &$b_{15}^{19}w_6$&$b_{15}^{20}w_4$&$b_{15}^{21}w_{5}$&$b_{15}^{22}w_9$&$b_{15}^{23}w_7$&$b_{15}^{24}w_8$&$b_{15}^{25}w_3$&$b_{15}^{26}$&$b_{15}^{27}w_2$\\
\vspace{-1mm}&&&&&&&&&\\
$w_{16}$ &$b_{16}^{19}w_7$&$b_{16}^{20}w_8$&$b_{16}^{21}w_{9}$&$b_{16}^{22}$&$b_{16}^{23}w_2$&$b_{16}^{24}w_3$&$b_{16}^{25}w_4$&$b_{16}^{26}w_5$&$b_{16}^{27}w_6$\\
\vspace{-1mm}&&&&&&&&&\\
$w_{17}$ &$b_{17}^{19}w_8$&$b_{17}^{20}w_9$&$b_{17}^{21}w_{7}$&$b_{17}^{22}w_2$&$b_{17}^{23}w_3$&$b_{17}^{24}$&$b_{17}^{25}w_5$&$b_{17}^{26}w_6$&$b_{17}^{27}w_4$\\
\vspace{-1mm}&&&&&&&&&\\
$w_{18}$ &$b_{18}^{19}w_9$&$b_{18}^{20}w_7$&$b_{18}^{21}w_{8}$&$b_{18}^{22}w_3$&$b_{18}^{23}$&$b_{18}^{24}w_2$&$b_{18}^{25}w_6$&$b_{18}^{26}w_4$&$b_{18}^{27}w_5$\\
\vspace{-1mm}&&&&&&&&&\\
$w_{19}$ &$b_{19}^{19}w_{10}$&$b_{19}^{20}w_{11}$&$b_{19}^{21}w_{12}$&$b_{19}^{22}w_{13}$&$b_{19}^{23}w_{14}$&$b_{19}^{24}w_{15}$&$b_{19}^{25}w_{16}$&$b_{19}^{26}w_{17}$&$b_{19}^{27}w_{18}$\\
\vspace{-1mm}&&&&&&&&&\\
$w_{20}$ &$b_{20}^{19}w_{11}$&$b_{20}^{20}w_{12}$&$b_{20}^{21}w_{10}$&$b_{20}^{22}w_{14}$&$b_{20}^{23}w_{15}$&$b_{20}^{24}w_{13}$&$b_{20}^{25}w_{17}$&$b_{20}^{26}w_{18}$&$b_{20}^{27}w_{16}$\\
\vspace{-1mm}&&&&&&&&&\\
$w_{21}$ &$b_{21}^{19}w_{12}$&$b_{21}^{20}w_{10}$&$b_{21}^{21}w_{11}$&$b_{21}^{22}w_{15}$&$b_{21}^{23}w_{13}$&$b_{21}^{24}w_{14}$&$b_{21}^{25}w_{18}$&$b_{21}^{26}w_{16}$&$b_{21}^{27}w_{17}$\\
\vspace{-1mm}&&&&&&&&&\\
$w_{22}$ &$b_{22}^{19}w_{13}$&$b_{22}^{20}w_{14}$&$b_{22}^{21}w_{15}$&$b_{22}^{22}w_{15}$&$b_{22}^{23}w_{17}$&$b_{22}^{24}w_{18}$&$b_{22}^{25}w_{10}$&$b_{22}^{26}w_{11}$&$b_{22}^{27}w_{12}$\\
\vspace{-1mm}&&&&&&&&&\\
$w_{23}$ &$b_{23}^{19}w_{14}$&$b_{23}^{20}w_{15}$&$b_{23}^{21}w_{13}$&$b_{23}^{22}w_{17}$&$b_{23}^{23}w_{18}$&$b_{23}^{24}w_{16}$&$b_{23}^{25}w_{11}$&$b_{23}^{26}w_{12}$&$b_{23}^{27}w_{10}$\\
\vspace{-1mm}&&&&&&&&&\\
$w_{24}$ &$b_{24}^{19}w_{15}$&$b_{24}^{20}w_{13}$&$b_{24}^{21}w_{14}$&$b_{24}^{22}w_{18}$&$b_{24}^{23}w_{16}$&$b_{24}^{24}w_{17}$&$b_{24}^{25}w_{12}$&$b_{24}^{26}w_{10}$&$b_{24}^{27}w_{11}$\\
\vspace{-1mm}&&&&&&&&&\\
$w_{25}$ &$b_{25}^{19}w_{16}$&$b_{25}^{20}w_{17}$&$b_{25}^{21}w_{18}$&$b_{25}^{22}w_{10}$&$b_{25}^{23}w_{11}$&$b_{25}^{24}w_{12}$&$b_{25}^{25}w_{13}$&$b_{25}^{26}w_{14}$&$b_{25}^{27}w_{15}$\\
\vspace{-1mm}&&&&&&&&&\\
$w_{26}$ &$b_{26}^{19}w_{17}$&$b_{26}^{20}w_{18}$&$b_{26}^{21}w_{16}$&$b_{26}^{22}w_{11}$&$b_{26}^{23}w_{12}$&$b_{26}^{24}w_{10}$&$b_{26}^{25}w_{14}$&$b_{26}^{26}w_{15}$&$b_{26}^{27}w_{13}$\\
\vspace{-1mm}&&&&&&&&&\\
$w_{27}$ &$b_{27}^{19}w_{18}$&$b_{27}^{20}w_{16}$&$b_{27}^{21}w_{17}$&$b_{27}^{22}w_{12}$&$b_{27}^{23}w_{10}$&$b_{27}^{24}w_{11}$&$b_{27}^{25}w_{15}$&$b_{27}^{26}w_{13}$&$b_{27}^{27}w_{14}$\\
\hline
\end{tabular}
\end{table}}
\end{center}

 \end{document}